\documentclass[11pt]{amsart}

\usepackage[a4paper,margin=1in]{geometry}
\usepackage[T1]{fontenc}
\usepackage{lmodern}
\usepackage{amsmath,amssymb,mathtools,mathrsfs}
\usepackage[nopatch=footnote]{microtype}
\usepackage{xcolor}
\usepackage{hyperref}
\usepackage{needspace}

\newtheorem{theorem}{Theorem}

\newtheorem{maintheorem}{Theorem}

\newtheorem{lemma}{Lemma}
\newtheorem{citedlemma}[lemma]{Lemma}
\newtheorem{proposition}{Proposition}
\newtheorem{corollary}{Corollary}
\newtheorem{claim}{Claim}[lemma]

\theoremstyle{remark}
\newtheorem{remark}{Remark}
\theoremstyle{definition}
\newtheorem{question}{Question}
\newtheorem{example}{Example}
\theoremstyle{plain}
\newcommand{\N}{\mathbb{N}}
\newcommand{\Z}{\mathbb{Z}}
\newcommand{\cL}{\mathcal{L}}
\newcommand{\cM}{\mathcal{M}}
\newcommand{\htop}{h_{\mathrm{top}}}
\newcommand{\hM}{h_{\mathcal{M}}}
\newcommand{\hG}{h_{\mathcal{G}}}
\newcommand{\TV}{\operatorname{TV}}

\title[Generating Forbidden Words and \mbox{Measures of Maximal Entropy}]{Generating Forbidden Words and
Measures of Maximal Entropy for Subshifts}
\author[X. Wu]{Xinxing Wu}
\address{School of Mathematics and Statistics, Guizhou University of Finance and Economics,
Guiyang, Guizhou 550025, China}
\email{wuxinxing5201314@163.com}
\author[X. Zhang]{Xu Zhang}
\address{Department of Mathematics,
Shandong University, Weihai, Shandong, 264209, China}
\email{xu\_zhang\_sdu@mail.sdu.edu.cn}
\date{\today}

\begin{document}

\begin{abstract}
We prove that if the entropy of generating forbidden words of a subshift
over a finite alphabet is strictly smaller than its topological entropy,
then the subshift has only finitely many ergodic measures of maximal
entropy, each of which is Bernoulli up to a finite period, and admits
a characteristic measure of maximal entropy.
In particular, these conclusions hold when the entropy of minimal
forbidden words is strictly smaller than the topological entropy,
giving a negative answer to a question of Buzzi~\cite{Buzzi2005}.
\end{abstract}
\maketitle

\section{Introduction}
\label{sec:introduction}

\subsection{Background and main result}
\label{sec:background-main-result}
Measures of maximal entropy provide a fundamental link
between topological dynamics and ergodic theory by realizing the
topological entropy of a dynamical system as the measure-theoretic
entropy of an invariant probability measure. Every subshift over a finite alphabet admits
such a measure (\cite[Chapter~8]{Walters1982}).
For these systems, we focus on the number and measure-theoretic
structure of ergodic measures of maximal entropy.

Buzzi~\cite{Buzzi2005} gave an entropy-based
approach to these questions by introducing subshifts
of quasi-finite type, defined by the condition
\(h_{\mathcal{SC}}(\Sigma)<\htop(\Sigma)\), where
\(h_{\mathcal{SC}}(\Sigma)\) denotes the symmetric constraint entropy
and \(\htop(\Sigma)\) the topological entropy of \(\Sigma\).
Using Markov diagrams, he~\cite[Theorem~1~(1)]{Buzzi2005}
proved that every
subshift of quasi-finite type has only finitely many ergodic measures of maximal entropy,
each of which is Bernoulli up to a finite period.
Buzzi~\mbox{\cite{Buzzi2010}} later extended this approach to
puzzles of quasi-finite type, proving existence and finiteness of ergodic
measures of maximal entropy, whose natural extensions are Bernoulli up to
a finite period.

Another description of subshifts uses minimal forbidden words,
studied in symbolic dynamics by
B\'eal, Mignosi, and Restivo~\cite{BealMignosiRestivo1996}, and
subsequently by B\'eal et al.~\cite{BealEtAl2000}.
A minimal forbidden word is a word absent from the language
whose proper contiguous subwords all belong to the language.
These words uniquely determine the subshift
(\cite[p.~166]{BealEtAl2000}), and their finiteness characterizes
subshifts of finite type (\cite[p.~168]{BealEtAl2000}).
Rather than requiring finitely many minimal forbidden
words, one can compare their entropy \(\hM(\Sigma)\), which records
their exponential growth rate, with the topological entropy
\(\htop(\Sigma)\).
The resulting condition \(\hM(\Sigma)<\htop(\Sigma)\) permits
infinitely many minimal forbidden words.
By~\mbox{\eqref{eq:entropy-class-inclusions}}
and Example~\mbox{\ref{ex:strict-entropy-classes}~(1)},
the class of subshifts satisfying
\(\hM(\Sigma)<\htop(\Sigma)\) strictly contains the class of
subshifts of quasi-finite type.
For this larger class, Buzzi asked whether the finiteness conclusion
for ergodic measures of maximal entropy still holds:
\begin{question}[{\cite[footnote~15]{Buzzi2005}}]
\label{ques:buzzi-maximal-measures}
Does there exist a subshift \(\Sigma\) with
\(\hM(\Sigma)<\htop(\Sigma)\) and infinitely many ergodic measures of maximal entropy\footnote{Buzzi~\mbox{\cite[p.~372]{Buzzi2005}} calls an ergodic measure of maximal entropy a \emph{maximal measure}.}?
\end{question}

A third description uses sets of forbidden words generating a subshift.
A set of forbidden words \(\mathcal{F}\) generates \(\Sigma\)
if a sequence belongs to \(\Sigma\) exactly when none of its finite
contiguous subwords belongs to \(\mathcal{F}\).
Using this description, we introduce the entropy of generating forbidden words,
denoted by \(\hG(\Sigma)\), as the infimum of the exponential growth
rates of all such sets generating \(\Sigma\); see
Section~\mbox{\ref{sec:preliminaries}} for the precise definition.
By~\mbox{\eqref{eq:entropy-class-inclusions}}
and Example~\mbox{\ref{ex:strict-entropy-classes}~(2)},
the class of subshifts satisfying
\(\hG(\Sigma)<\htop(\Sigma)\) strictly contains the class satisfying
\(\hM(\Sigma)<\htop(\Sigma)\).
For this larger class, we prove the following main result:

\begin{maintheorem}
\label{thm:finiteness}
If \(\Sigma\) is a subshift over a finite alphabet and
\(\hG(\Sigma)<\htop(\Sigma)\),
then \(\Sigma\) has only finitely many ergodic measures of maximal entropy,
each of which is Bernoulli up to a finite period.
Moreover, \(\Sigma\) admits a characteristic measure of maximal entropy.
\end{maintheorem}

Theorem~\ref{thm:finiteness} strictly extends Buzzi's conclusions on the
finiteness of ergodic measures of maximal entropy and the fact that each is
Bernoulli up to a finite period~\mbox{\cite[Theorem~1~(1)]{Buzzi2005}}
beyond subshifts of quasi-finite type. It also extends the corresponding
conclusions of~\mbox{\cite[Theorem~1.11]{Buzzi2010}} in the subshift setting.
Moreover, our proof does not rely on Markov diagrams, thereby
giving a partial answer to his question about direct proofs
\cite[Section~1.4, p.~380]{Buzzi2005}
(Remark~\ref{rem:buzzi-comparison}).
In particular, Theorem~\mbox{\ref{thm:finiteness}}
applies whenever \(\hM(\Sigma)<\htop(\Sigma)\).
This consequence, stated as Corollary~\mbox{\ref{cor:minimal-forbidden-entropy}}
below, answers Question~\mbox{\ref{ques:buzzi-maximal-measures}} negatively.

As an application, we prove that every power-free subshift \(\Sigma\)
with positive topological entropy satisfies
\(\hG(\Sigma)<\htop(\Sigma)\)
(Proposition~\ref{prop:power-free-entropy})
and is not of quasi-finite type
(Corollary~\ref{cor:power-free-measures}).
In particular, these include the classical power-free subshifts
\(\Sigma_2^{(3)}\) and \(\Sigma_3^{(2)}\),
defined in Section~\ref{sec:applications}.

\subsection{Related work}
\label{sec:related-work}
Pavlov~\mbox{\cite{Pavlov2022}} obtained a unique measure of
maximal entropy with the \(K\)-property under an explicit weighted
bound on the numbers of forbidden words of each length.
Climenhaga~\mbox{\cite[Theorem~1.1]{ClimenhagaPowerFree2025}}
used a weak form of specification to prove uniqueness of the measure
of maximal entropy for the power-free subshifts
\(\Sigma_\beta^{(d)}\) and \(\Sigma_{\beta^+}^{(d)}\)
when \mbox{\(d\ge2\)} and \mbox{\(\beta>12\)}.
The lengths of minimal forbidden words also determine the language
stability condition of Cyr and Kra~\cite{CyrKra2023}, under which
they obtained characteristic measures of maximal entropy.
Cyr, Kra, and
Petite~\cite{CyrKraPetite2025} further developed approximation
criteria for the existence of characteristic measures and gave
additional conditions ensuring maximal entropy.
Continuing the study of finite-type approximations, Coles et al.~\mbox{\cite{ColesCyrKraPavlov2026}}
introduced period-stable subshifts and proved that every such
subshift admits a characteristic measure.

Measures of maximal entropy are precisely the equilibrium states for the zero potential.
More generally, Climenhaga and Thompson~\mbox{\cite{ClimenhagaThompson2013}} combined specification and the
Bowen property on suitable word collections with a pressure gap for
prefixes and suffixes to prove uniqueness of equilibrium states.
Their obstruction-entropy criterion~\mbox{\cite{ClimenhagaThompson2014}}
gives a unique measure of maximal entropy when the entropies of
obstructions to specification and positive expansivity are both
strictly smaller than the topological entropy.
For almost specification, Pavlov~\cite[Theorem~1.2]{Pavlov2016}
constructed a subshift with mistake function \(g(n)=4\) and exactly
two ergodic measures of maximal entropy with disjoint supports.
For coded shifts, the entropy-gap criterion of Kucherenko, Schmoll,
and Wolf~\mbox{\cite{KucherenkoSchmollWolf2024}} gives a unique
measure of maximal entropy when the generating set uniquely
represents the concatenation set and the concatenation entropy
exceeds the residual entropy.

On topologically transitive one-sided countable Markov shifts,
Buzzi and Sarig~\mbox{\cite{BuzziSarig2003}} proved that
bounded-above potentials with summable variations and finite
Gurevich pressure have at most one equilibrium state.
Complementing this uniqueness result, Sarig~\mbox{\cite{Sarig2011}}
proved the Bernoulli property for equilibrium states of bounded-above
weakly H\"older potentials with finite Gurevich pressure on
topologically mixing two-sided countable Markov shifts.
Using strongly positive recurrent Markov codings,
Climenhaga~\mbox{\cite{Climenhaga2018}} obtained unique
equilibrium states, Bernoulli up to a finite period, for H\"older
potentials on two-sided subshifts under specification, overlap,
and pressure-gap hypotheses.
A geometric approach to equilibrium states was presented by
Climenhaga, Pesin, and Zelerowicz~\mbox{\cite{ClimenhagaPesinZelerowicz2019}},
who used Carath\'eodory reference measures on local unstable leaves to
construct the unique equilibrium state for each H\"older potential on a
topologically transitive locally maximal hyperbolic set of a
\mbox{\(C^{1+\alpha}\)} diffeomorphism.

For surface systems, Buzzi~\mbox{\cite{Buzzi2009}} proved
finiteness of ergodic measures of maximal entropy for piecewise affine homeomorphisms
of compact affine surfaces with positive topological entropy.
For diffeomorphisms of closed surfaces with positive topological
entropy, Sarig~\mbox{\cite{Sarig2013}} proved that there are at most
countably many ergodic measures of maximal entropy in the \(C^{1+\beta}\) class,
\(0<\beta<1\).
In the \(C^\infty\) case, Buzzi, Crovisier, and
Sarig~\cite{BuzziCrovisierSarig2022} strengthened this to finiteness,
and to uniqueness under topological transitivity, using homoclinic
classes and symbolic codings.
In subsequent work, Buzzi, Crovisier, and Sarig~\mbox{\cite{BuzziCrovisierSarig2025}}
proved that positive-entropy \mbox{\(C^\infty\)}
diffeomorphisms of closed surfaces are strongly positive recurrent.
In the topologically mixing case, they obtained exponential decay of
correlations for H\"older observables with respect to the unique measure
of maximal entropy.
In a recent survey, Buzzi~\mbox{\cite{Buzzi2026Survey}} discusses
strong positive recurrence and its statistical consequences for smooth
surface diffeomorphisms.

Beyond the surface setting, Lima, Obata, and Poletti~\mbox{\cite{LimaObataPoletti2026}}
considered possibly noninvertible maps that are
\mbox{\(C^{1+}\)} outside their singular sets.
Under their assumptions~(A1)--(A7), they proved that each homoclinic class
of an ergodic adapted hyperbolic measure carries at most one adapted
hyperbolic measure of maximal entropy. When such a measure exists, its natural
extension is Bernoulli up to a finite period.
For finite-horizon Sinai billiard maps, Baladi and
Demers~\mbox{\cite[condition~(1.5)]{BaladiDemers2020}}
obtained a unique measure of maximal entropy with full support and the Bernoulli
property under sparse recurrence to singularities.
More recently, Climenhaga and Day~\mbox{\cite{ClimenhagaDay2026}}
removed the sparse recurrence assumption and proved that every
finite-horizon Sinai billiard map has a unique measure of maximal entropy,
which is multiply mixing and has full support and local product structure.

For flows, Bowen and Ruelle~\mbox{\cite{BowenRuelle1975}}
established uniqueness of equilibrium states for H\"older potentials
on basic sets of Axiom~A flows using symbolic models.
Beyond uniform hyperbolicity, Burns et al.~\mbox{\cite{BurnsClimenhagaFisherThompson2018}}
proved that a H\"older potential for the geodesic flow on a closed
rank one manifold of nonpositive curvature has a unique equilibrium
state when its pressure on the singular set is strictly smaller than its
pressure on the whole system.

\subsection{Proof strategy}
\label{sec:proof-strategy}
The proof treats the zero- and positive-entropy cases in
Sections~\ref{sec:zero-entropy} and~\ref{sec:positive-entropy},
respectively.
In the zero-entropy case, the entropy gap forces all invariant
probability measures to be supported on a common finite union of
periodic orbits.  The ergodic measures of maximal entropy are therefore precisely
the uniform probability measures on these orbits.

In the positive-entropy case, the argument uses equal-length word
replacements.  Such replacements also appear in the comparison of
cylinder measures through extender sets in Garc\'ia-Ramos and
Pavlov~\cite{GarciaRamosPavlov2019}.
The entropy gap and the uniform estimate in
Lemma~\ref{lem:uniform-coverage} imply that these replacements
preserve every invariant measure of maximal entropy on their
admissible domains
(Lemma~\ref{lem:finite-replacement-invariance}).
Applied to the product subshift (Proposition~\mbox{\ref{prop:product-forbidden-entropy}}),
this invariance gives block exchanges preserving products of
measures of maximal entropy.  Together with the uniform estimate, these
exchanges yield uniform separation of the cylinder distributions
of distinct ergodic measures of maximal entropy at a fixed block length
(Lemma~\ref{lem:positive-entropy-separation}).
Finiteness then follows from compactness of the corresponding
finite-dimensional probability simplex.

For a single ergodic measure of maximal entropy, the same exchanges, applied to two
independent copies, show that the total variation distance between its
joint past--future law and the product of its marginals
is strictly smaller than \(1\).
Lemma~\ref{lem:berbee-symbolic} then gives a finite cyclic
decomposition.  The estimates in this lemma imply that the
finite block partitions on each component are weak Bernoulli for
the corresponding power of the shift, so
Lemma~\ref{lem:friedman-ornstein} yields the Bernoulli conclusion.

\section{Preliminaries}
\label{sec:preliminaries}

Throughout this paper, let
\(\Z=\{\ldots,-1,0,1,\ldots\}\),
\(\N=\{1,2,3,\ldots\}\), and \(\N_{0}=\{0,1,2,\ldots\}\).
A \emph{subshift} is a nonempty closed \(\sigma\)-invariant
subset \(\Sigma\subseteq\mathcal{A}^{\Z}\) (\(\mathcal{A}\) is some finite
set, \(\sigma:(x_{n})_{n\in\Z}\mapsto(x_{n+1})_{n\in\Z}\) is the
\emph{left-shift} on \(\mathcal{A}^{\Z}\) endowed with the product topology).
All measures are Borel measures.
For measures on \(\Sigma\), invariance and ergodicity refer to
\(\sigma|_\Sigma:\Sigma\to\Sigma\).
Write \(\mathscr{E}(\Sigma,\sigma)\) for the set of ergodic invariant
probability measures on \(\Sigma\).

An \emph{automorphism} of
\(\Sigma\) is a homeomorphism
\(\varphi:\Sigma\to\Sigma\) commuting with
\(\sigma|_\Sigma\).
Following~\mbox{\cite{CyrKra2023}}, a probability measure
\(\mu\) on \(\Sigma\)
is \emph{characteristic} if
\(\mu(\varphi^{-1}(B))=\mu(B)\)
for any automorphism \(\varphi\)
of \(\Sigma\) and any Borel set
\(B\subseteq\Sigma\).

For any \(x\in\mathcal{A}^{\Z}\) and any \(i,j\in\Z\) with \(i<j\),
the notation \(x_{[i,j)}\) denotes the word
\(x_{i}\cdots x_{j-1}\).
Let \(\mathbf{e}\) denote the empty word.
For any \(i\in\Z\), \(x_{[i,i)}=\mathbf{e}\).
We write \(|w|\) for the length of a finite
word \(w\), and \(uv\) for the concatenation of two finite words.
For any \(w\in\mathcal{A}^{n}\) (\(n\in\N\)),
let \([w]=\{x\in\mathcal{A}^{\Z}:x_{0}\cdots x_{n-1}=w\}\)
and \([w]_{\Sigma}=[w]\cap\Sigma\).
The \emph{language} of length \(n\) is
\(\cL(\Sigma,n)=\{w\in\mathcal{A}^{n}:
[w]_{\Sigma}\neq\varnothing\}\).
We set \([\mathbf{e}]=\mathcal{A}^{\Z}\),
\([\mathbf{e}]_{\Sigma}=\Sigma\), and
\(\cL(\Sigma,0)=\{\mathbf{e}\}\).

A finite word is \emph{forbidden} for
\(\Sigma\) if it does not occur in any point of
\(\Sigma\).
Thus, for any \(n\in\N\), the set of forbidden words of length \(n\)
for \(\Sigma\) is
\(\mathcal{A}^{n}\setminus\cL(\Sigma,n)
 =\{w\in\mathcal{A}^{n}:[w]_{\Sigma}=\varnothing\}.\)
The \emph{minimal forbidden words} of \(\Sigma\) are
the forbidden words with no proper forbidden factor,
where a factor is a contiguous subword.
Let \(\cM(\Sigma,n)\) denote the set of minimal
forbidden words of length \(n\in\N\).
By \cite[Remark~3]{BealEtAl2000}, for any
\(w=w_{1}\cdots w_{n}\in\mathcal{A}^{n}\),
\begin{equation}\label{eq:minimal-forbidden-characterization}
 w\in\cM(\Sigma,n)
 \ \Longleftrightarrow\
 \begin{cases}
 [w_{1}\cdots w_{n}]_{\Sigma}=\varnothing,\\
 [w_{1}\cdots w_{n-1}]_{\Sigma}\neq\varnothing,\\
 [w_{2}\cdots w_{n}]_{\Sigma}\neq\varnothing.
 \end{cases}
\end{equation}

For any set \(\mathcal{F}\subseteq\bigcup_{n\in\N}\mathcal{A}^n\),
define the closed \(\sigma\)-invariant set
\mbox{(\cite[Definition~2.6]{Pavlov2022})}
\[
 \Sigma_{\mathcal{F}}=\{x\in\mathcal{A}^{\Z}:
 x_{[i,i+n)}\notin\mathcal{F}
 \text{ for all }i\in\Z\text{ and }n\in\N\}.
\]
We call \(\mathcal{F}\)
a \emph{set of forbidden words generating} \(\Sigma\)
if \(\Sigma_{\mathcal{F}}=\Sigma\).
Such a set need not contain all forbidden words of
\(\Sigma\).
For any such set \(\mathcal{F}\) and any \(n\in\N\), write
\(\mathcal{F}_n=\mathcal{F}\cap\mathcal{A}^n\).
Then \(\mathcal{F}_n\subseteq\mathcal{A}^n\setminus\cL(\Sigma,n)\).
Every subshift admits such a set:
\(\bigcup_{n\in\N}\cM(\Sigma,n)\)
generates \(\Sigma\), since every forbidden word contains a minimal forbidden word.
A subshift \(\Sigma\) is
\emph{of finite type} if there exists a finite set of forbidden
words generating \(\Sigma\). Equivalently,
\(\bigcup_{n\in\N}\cM(\Sigma,n)\) is finite
\mbox{(\cite[p.~168]{BealEtAl2000}).}

Recall that the \emph{topological entropy} of $\Sigma$ is defined by
\begin{equation}
\label{eq:topological-entropy}
 \htop(\Sigma)=\limsup_{n\rightarrow\infty}
                  \frac{1}{n}\log\#\cL(\Sigma,n).
\end{equation}
Moreover, by subadditivity and \cite[Theorem~4.9]{Walters1982},
\begin{equation}
\label{eq:language-limit}
 \htop(\Sigma)=\lim_{n\rightarrow\infty}
 \frac{1}{n}\log\#\cL(\Sigma,n)
 =\inf_{n\in\N}\frac{\log\#\cL(\Sigma,n)}{n}.
\end{equation}

Following Buzzi~\mbox{\cite{Buzzi2005}}, we
recall the following definitions.
Let \(\Sigma^+=\{(x_i)_{i\geq0}:x\in\Sigma\}\).
For any word \(w=w_1\cdots w_n\in\mathcal{A}^n\)
(\(n\in\N\)),
the \emph{follower set} of \(w\) is defined by
\[
 \operatorname{Fol}(w)
 =\{(x_i)_{i\geq0}:x\in\Sigma,\ x_{1-n}\cdots x_0=w\},
\]
with \(\operatorname{Fol}(\mathbf{e})=\Sigma^+\).
A word \(w=w_1\cdots w_n\in\mathcal{A}^n\)
(\(n\in\N\)) is a
\emph{left constraint} if
\[
 \varnothing\neq\operatorname{Fol}(w)
 \subsetneq\operatorname{Fol}(w_2\cdots w_n),
\]
where \(w_2\cdots w_n=\mathbf{e}\) when \(n=1\).
For any \(n\in\N\), let \(\mathcal{C}(\Sigma,n)\) denote the set of
left constraints of length \(n\).
The \emph{left constraint entropy} of \(\Sigma\) is defined by
\[
 h_{\mathcal{C}}(\Sigma)
 =\limsup_{n\to\infty}\frac{1}{n}\log\#\mathcal{C}(\Sigma,n),
\]
where \(\log0=-\infty\).
Let \(\bar{\Sigma}=\{(x_{-i})_{i\in\Z}:x\in\Sigma\}\) denote the
time-reversed subshift.
The \emph{symmetric constraint entropy} of \(\Sigma\) is defined by
\[
 h_{\mathcal{SC}}(\Sigma)
 =\min\{h_{\mathcal{C}}(\Sigma),h_{\mathcal{C}}(\bar{\Sigma})\}.
\]
The subshift \(\Sigma\) is of
\emph{quasi-finite type} if
\(h_{\mathcal{SC}}(\Sigma)<\htop(\Sigma)\).

The \emph{entropy of minimal forbidden words} is defined by
(\cite[Lemma~15]{Buzzi2005})\footnote{The
formula in \cite[Lemma~15]{Buzzi2005} omits the logarithm.
In \cite[p.~389, proof of Lemma~15]{Buzzi2005},
the family \((a|b)(0|1|2)^{n}(a|b)^{n^{2}+1}\) has length
\(n^{2}+n+2\) and cardinality \(3^{n}2^{n^{2}+2}\).
Thus \(\frac{\log(3^{n}2^{n^{2}+2})}{n^{2}+n+2}\rightarrow\log 2\), whereas
\(\frac{3^{n}2^{n^{2}+2}}{n^{2}+n+2}\rightarrow+\infty\).
The value \(\hM(\Sigma_{1})=\log 2\) computed there therefore requires
the logarithm.}
\begin{equation}
\label{eq:entropy-definitions}
 \hM(\Sigma)=\limsup_{n\rightarrow\infty}
                  \frac{1}{n}\log\#\cM(\Sigma,n).
\end{equation}

By \cite[Lemma~15]{Buzzi2005},
\(\hM(\Sigma)\leq h_{\mathcal{C}}(\Sigma)\).
For any \(n\in\N\), reversal of words gives a bijection between
\(\cM(\Sigma,n)\) and \(\cM(\bar{\Sigma},n)\), so
\(\hM(\bar{\Sigma})=\hM(\Sigma)\).
By these facts and
\(\mathcal{C}(\Sigma,n)\subseteq\cL(\Sigma,n)\) for all \(n\in\N\),
\begin{equation}\label{eq:minimal-constraint-entropy}
 \hM(\Sigma)\leq h_{\mathcal{SC}}(\Sigma)\leq\htop(\Sigma).
\end{equation}
In particular, every subshift of quasi-finite type satisfies
\(\hM(\Sigma)<\htop(\Sigma)\).
The converse fails, as shown in Example~\mbox{\ref{ex:strict-entropy-classes}}~(1).

We define the \emph{entropy of generating forbidden words} of
\(\Sigma\) by
\begin{equation}\label{eq:generating-forbidden-entropy}
 \hG(\Sigma)
 =\inf_{\substack{
 \mathcal{F}\subseteq\bigcup_{n\in\N}\mathcal{A}^n\\
 \Sigma_{\mathcal{F}}=\Sigma}}
 \limsup_{n\to\infty}
 \frac{1}{n}\log\#(\mathcal{F}\cap\mathcal{A}^n).
\end{equation}
Since the set of all minimal forbidden words
\(\bigcup_{n\in\N}\cM(\Sigma,n)\)
generates \(\Sigma\),
\[
 \hG(\Sigma)\leq\hM(\Sigma).
\]
This, together with~\mbox{\eqref{eq:minimal-constraint-entropy}}, implies that
\begin{equation}\label{eq:entropy-comparison}
 \hG(\Sigma)\leq\hM(\Sigma)\leq h_{\mathcal{SC}}(\Sigma)\leq\htop(\Sigma).
\end{equation}
It is not difficult to check that
the entropy of generating forbidden words \(h_{\mathcal{G}}\)
is invariant under topological conjugacy.

Let \(\mu\) be an invariant probability measure on \(\Sigma\).
For \(n\in\N\), the partition \(\{[w]_\Sigma:w\in\cL(\Sigma,n)\}\)
has entropy (\cite[Definition~4.6]{Walters1982})
\begin{equation}
\label{eq:block-entropy}
 H_{\mu}(n)=-\sum_{w\in\cL(\Sigma,n)}
                           \mu([w]_{\Sigma})\log\mu([w]_{\Sigma}),
\end{equation}
where \(0\log 0=0\).
The \emph{measure-theoretic entropy} of \(\sigma|_\Sigma\)
with respect to \(\mu\) is given by
(\cite[Theorems~4.9 and~4.17]{Walters1982})
\begin{equation}
\label{eq:measure-entropy}
 h_{\mu}(\sigma|_\Sigma)=\lim_{n\rightarrow\infty}\frac{H_{\mu}(n)}{n}
                    =\inf_{n\in\N}\frac{H_{\mu}(n)}{n}.
\end{equation}
By \eqref{eq:measure-entropy} and
\cite[Corollary~4.2.1]{Walters1982},
\begin{equation}
\label{eq:block-entropy-bounds}
 n h_{\mu}(\sigma|_\Sigma)\leq H_{\mu}(n)\leq\log\#\cL(\Sigma,n).
\end{equation}
A \emph{measure of maximal entropy} is an invariant probability measure
whose measure-theoretic entropy equals the topological entropy.
An ergodic measure of maximal entropy exists for every subshift over a finite alphabet
(\cite[pp.~191--192]{Walters1982}).

For any nonempty finite alphabet \(\mathcal{A}\), let
\(\mathscr{C}_{\mathcal{SC}}^{(\mathcal{A})}\),
\(\mathscr{C}_{\mathcal{M}}^{(\mathcal{A})}\), and
\(\mathscr{C}_{\mathcal{G}}^{(\mathcal{A})}\)
denote the classes of subshifts \(\Sigma\subseteq\mathcal{A}^{\Z}\)
satisfying
\[
 h_{\mathcal{SC}}(\Sigma)<\htop(\Sigma),\quad
 \hM(\Sigma)<\htop(\Sigma),\quad\text{and}\quad
 \hG(\Sigma)<\htop(\Sigma),
\]
respectively. Thus \(\mathscr{C}_{\mathcal{SC}}^{(\mathcal{A})}\)
is the class of subshifts of quasi-finite type over
\(\mathcal{A}\). By~\mbox{\eqref{eq:entropy-comparison}},
\begin{equation}\label{eq:entropy-class-inclusions}
 \mathscr{C}_{\mathcal{SC}}^{(\mathcal{A})}
 \subseteq\mathscr{C}_{\mathcal{M}}^{(\mathcal{A})}
 \subseteq\mathscr{C}_{\mathcal{G}}^{(\mathcal{A})}.
\end{equation}
The following example shows that both inclusions are strict when
\(\#\mathcal{A}\geq4\).

\begin{example}
\label{ex:strict-entropy-classes}
Let \(\mathcal{A}=\{a_{0},a_{1},\ldots,a_{\#\mathcal{A}-1}\}\)
be a finite alphabet with \(\#\mathcal{A}\geq4\).

(1) By~\mbox{\cite[Lemma~16]{Buzzi2005}}, there exists a subshift
\(\Sigma\subseteq\{a_{0},a_{1}\}^{\Z}\) such that
\(\hM(\Sigma)=0<\htop(\Sigma)\) and
\(\Sigma\) has no periodic points.
Writing \(\cM_{\mathcal{B}}(\Sigma,n)\)
for the set of minimal forbidden words of length \(n\)
over the alphabet \(\mathcal{B}\),
\[
 \cM_{\mathcal{A}}(\Sigma,n)=
 \begin{cases}
 \cM_{\{a_{0},a_{1}\}}(\Sigma,1)\cup
       \{a_{2},\ldots,a_{\#\mathcal{A}-1}\},&n=1,\\
 \cM_{\{a_{0},a_{1}\}}(\Sigma,n),&n\geq2.
 \end{cases}
\]
Hence \(\hM(\Sigma)=0<\htop(\Sigma)\)
also over \(\mathcal{A}\). Moreover,
\(\Sigma\) is not of quasi-finite type, by
\mbox{\cite[Theorem~1~(3)]{Buzzi2005}}. Hence
\(
 \Sigma\in\mathscr{C}_{\mathcal{M}}^{(\mathcal{A})}
 \setminus\mathscr{C}_{\mathcal{SC}}^{(\mathcal{A})}.
\)

(2) Let
\[
 \begin{aligned}
 \mathcal{F}
 &=(\mathcal{A}\setminus\{a_{0},a_{1},a_{2},a_{3}\})\\
 &\quad\cup\{a_{0}a_{2},a_{1}a_{2},a_{2}a_{3},a_{3}a_{0},a_{3}a_{1},a_{3}a_{2}\}\\
 &\quad\cup\bigcup_{n\in\N}
   \{a_{2}u a_{3}^{n^2+1}:u\in\{a_{0},a_{1}\}^n\},
 \end{aligned}
\]
and put \(\Sigma=\Sigma_{\mathcal{F}}\). It can be verified that
\(\{a_{0},a_{1}\}^{\Z}\subseteq\Sigma\), and thus \(\Sigma\neq\varnothing\).

Fix any \(N\in\N\).
For any \(x=(x_i)_{i\in\Z}\in\Sigma\), the forbidden words in
\(\mathcal{A}\setminus\{a_{0},a_{1},a_{2},a_{3}\}\) give
\(x_i\in\{a_{0},a_{1},a_{2},a_{3}\}\) for all \(i\in\Z\).
Since \(a_{3}a_{0},a_{3}a_{1},a_{3}a_{2}\in\mathcal{F}\),
\(x_i=a_{3}\) implies \(x_{i+1}=a_{3}\), and hence
\(x_j=a_{3}\) for all \(j\geq i\), implying that
\(a_{2}u a_{3}\) does not occur in \(x\) for any
\(u\in\bigcup_{n\in\N}\{a_{0},a_{1}\}^{n}\)
(otherwise, \(x\) contains \(a_{2}u a_{3}^{|u|^{2}+1}\in\mathcal{F}\),
contradicting \(x\in\Sigma_{\mathcal{F}}\)).
Similarly, since \(a_{0}a_{2},a_{1}a_{2},a_{3}a_{2}\in\mathcal{F}\),
\(x_i=a_{2}\) implies \(x_{i-1}=a_{2}\), and hence
\(x_j=a_{2}\) for all \(j\leq i\).
To obtain \(\#\cL(\Sigma,N)\), we consider the following four mutually exclusive
cases for \(w=w_0\cdots w_{N-1}\in\cL(\Sigma,N)\):
\begin{enumerate}
\renewcommand{\theenumi}{\roman{enumi}}
\renewcommand{\labelenumi}{(\theenumi)}
\item Neither \(a_{2}\) nor \(a_{3}\) occurs in \(w\).
Then \(w\in\{a_{0},a_{1}\}^{N}\).
Conversely, for any \(u\in\{a_{0},a_{1}\}^{N}\), the sequence
\(\ldots a_{0}a_{0}a_{0}\,u\,a_{0}a_{0}a_{0}\ldots\)
belongs to \(\{a_{0},a_{1}\}^{\Z}\subseteq\Sigma\) and contains \(u\).
Thus there are exactly \(2^N\) words in \(\cL(\Sigma,N)\) in this case.

\item The letter \(a_{2}\) occurs in \(w\), but \(a_{3}\) does not.
Let \(p=\max\{i\in\{0,\ldots,N-1\}:w_i=a_{2}\}\).
Then \(w_i=a_{2}\) for \(0\leq i\leq p\), and
\(w_i\in\{a_{0},a_{1}\}\) for \(p<i<N\). Hence
\[
 w=a_{2}^{p+1}u,
 \quad u=w_{[p+1,N)}\in\{a_{0},a_{1}\}^{N-p-1}.
\]
Conversely, for any integer \(0\leq p<N\) and any
\(u\in\{a_{0},a_{1}\}^{N-p-1}\), the sequence
\(\ldots a_{2}a_{2}a_{2}\,u\,a_{0}a_{0}a_{0}\ldots\)
avoids every word in \(\mathcal{F}\), so it belongs to \(\Sigma\)
and contains \(a_{2}^{p+1}u\).
Since \(p\) is uniquely determined by the word \(a_{2}^{p+1}u\),
there are exactly \(\sum_{p=0}^{N-1}2^{N-p-1}\) words in \(\cL(\Sigma,N)\)
in this case.

\item The letter \(a_{3}\) occurs in \(w\), but \(a_{2}\) does not.
Let \(q=\min\{i\in\{0,\ldots,N-1\}:w_i=a_{3}\}\).
Then \(w_i\in\{a_{0},a_{1}\}\) for \(0\leq i<q\), and
\(w_i=a_{3}\) for \(q\leq i<N\). Hence
\[
 w=u a_{3}^{N-q},
 \quad u=w_{[0,q)}\in\{a_{0},a_{1}\}^{q}.
\]
Conversely, for any integer \(0\leq q<N\) and any
\(u\in\{a_{0},a_{1}\}^{q}\), the sequence
\(\ldots a_{0}a_{0}a_{0}\,u\,a_{3}a_{3}a_{3}\ldots\)
avoids every word in \(\mathcal{F}\), so it belongs to \(\Sigma\)
and contains \(u a_{3}^{N-q}\).
Since \(q\) is uniquely determined by the word \(u a_{3}^{N-q}\),
there are exactly \(\sum_{q=0}^{N-1}2^q\) words in \(\cL(\Sigma,N)\)
in this case.

\item Both \(a_{2}\) and \(a_{3}\) occur in \(w\).
Let \(p=\max\{i\in\{0,\ldots,N-1\}:w_i=a_{2}\}\) and
\(q=\min\{i\in\{0,\ldots,N-1\}:w_i=a_{3}\}\).
If \(q\leq p\), then \(w_q=a_{2}\) by the implication proved above,
contradicting \(w_q=a_{3}\). Thus \(p<q\), and
\[
 w_{[p,q+1)}=a_{2}u a_{3},
 \quad u=w_{[p+1,q)}\in\{a_{0},a_{1}\}^{q-p-1}.
\]
If \(q=p+1\), this contradicts \(a_{2}a_{3}\in\mathcal{F}\).
If \(q\geq p+2\), then \(|u|\geq1\), so \(a_{2}u a_{3}\) is forbidden
by the preceding argument, again contradicting \(w\in\cL(\Sigma,N)\).
Therefore this case is impossible.
\end{enumerate}
Adding the numbers of words in cases~(i)--(iii) gives
\[
 \#\cL(\Sigma,N)
 =2^N+\sum_{p=0}^{N-1}2^{N-p-1}+\sum_{q=0}^{N-1}2^q
 =3\cdot2^N-2.
\]
Thus \(\htop(\Sigma)=\log2\).

For any \(n\in\N\) and any \(u\in\{a_{0},a_{1}\}^{n}\),
case~(iv), applied with \(N=n+2\), shows that
\(a_{2}u a_{3}\notin\cL(\Sigma,n+2)\).
Cases~(ii) and~(iii), applied with \(N=n+1\) and with
\(p=0\) and \(q=n\), respectively, show that
\(a_{2}u,u a_{3}\in\cL(\Sigma,n+1)\).
By~\mbox{\eqref{eq:minimal-forbidden-characterization}},
\(a_{2}u a_{3}\in\cM(\Sigma,n+2)\). Consequently,
\(\#\cM(\Sigma,n+2)\geq2^n\) for all \(n\in\N\).
By~\mbox{\eqref{eq:entropy-definitions}},
\[
 \hM(\Sigma)
 =\limsup_{n\to\infty}\frac{1}{n+2}\log\#\cM(\Sigma,n+2)
 \geq\lim_{n\to\infty}\frac{n\log2}{n+2}
 =\log2.
\]
This, together with
\mbox{\eqref{eq:entropy-comparison}}, implies
\(\hM(\Sigma)=\htop(\Sigma)=\log2\).

By the definition of \(\mathcal{F}\), for any integer
\(N\geq3\),
\[
 \#\mathcal{F}_N=
 \begin{cases}
 2^n,&N=n^2+n+2\text{ for some }n\in\N,\\
 0,&\text{otherwise}.
 \end{cases}
\]
It follows from~\mbox{\eqref{eq:generating-forbidden-entropy}} that
\[
 \begin{aligned}
 \hG(\Sigma)
 &\leq\limsup_{N\to\infty}\frac{1}{N}\log\#\mathcal{F}_N
 =\lim_{n\to\infty}\frac{n\log2}{n^2+n+2}\\
 &=0<\htop(\Sigma)=\hM(\Sigma).
 \end{aligned}
\]
Therefore
\(\Sigma\in\mathscr{C}_{\mathcal{G}}^{(\mathcal{A})}
\setminus\mathscr{C}_{\mathcal{M}}^{(\mathcal{A})}\).
\end{example}

Let \(\Sigma\) be a subshift over a finite alphabet
\(\mathcal{A}\). Put
\[
 \mathcal{X}_{\Sigma}=\{((x_i,y_i))_{i\in\Z}:x,y\in\Sigma\}
 \subseteq(\mathcal{A}\times\mathcal{A})^{\Z},
\]
and define \(\Lambda:\Sigma^2\to \mathcal{X}_{\Sigma}\) by
\(\Lambda(x,y)=((x_i,y_i))_{i\in\Z}\). Then \(\mathcal{X}_{\Sigma}\)
is a subshift over $\mathcal{A}\times \mathcal{A}$, and
\(\Lambda\) is a homeomorphism satisfying
\(\Lambda\circ(\sigma\times\sigma)|_{\Sigma^2}
 =\sigma|_{\mathcal{X}_{\Sigma}}\circ\Lambda\).
In particular, \(\htop(\mathcal{X}_{\Sigma})=2\htop(\Sigma)\).

The following proposition shows that the condition
\(\hG(\Sigma)<\htop(\Sigma)\) passes to
\(\mathcal{X}_{\Sigma}\). In the proof of the uniform separation
estimate in Lemma~\mbox{\ref{lem:positive-entropy-separation}}, we use this fact
to apply Lemma~\mbox{\ref{lem:finite-replacement-invariance}} to the product
subshift and show that block exchanges preserve products of measures of maximal entropy.

\begin{proposition}
\label{prop:product-forbidden-entropy}
Let \(\Sigma\) be a subshift over a finite alphabet \(\mathcal{A}\).
Then
\(\hG(\mathcal{X}_\Sigma)\leq\hG(\Sigma)+\htop(\Sigma).\)
In particular, if \(\hG(\Sigma)<\htop(\Sigma)\), then
\(\hG(\mathcal{X}_\Sigma)<\htop(\mathcal{X}_\Sigma)\).
\end{proposition}

\begin{proof}
Fix any \(\varepsilon>0\) and any set of forbidden words \(\mathcal{F}\)
generating \(\Sigma\). By \mbox{\eqref{eq:language-limit}}, there exists
\(n_0\in\N\) such that
\(\frac{1}{n_0}\log\#\cL(\Sigma,n_0)<\htop(\Sigma)+\varepsilon\).
Let
\[
 \Sigma_{n_0}=\{x\in\mathcal{A}^\Z:
 x_{[i,i+n_0)}\in\cL(\Sigma,n_0)\text{ for all }i\in\Z\}.
\]
Then \(\Sigma_{n_0}\) is a subshift, and
\(\Sigma\subseteq\Sigma_{n_0}\) by the \(\sigma\)-invariance of \(\Sigma\).
This implies that \(\cL(\Sigma,n_0)\subseteq\cL(\Sigma_{n_0},n_0)\).
By the definition of \(\Sigma_{n_0}\), we also have
\(\cL(\Sigma_{n_0},n_0)\subseteq\cL(\Sigma,n_0)\).
Hence \(\cL(\Sigma_{n_0},n_0)=\cL(\Sigma,n_0)\).
By \mbox{\eqref{eq:language-limit}},
\[
 \htop(\Sigma_{n_0})\leq\frac{1}{n_0}\log\#\cL(\Sigma_{n_0},n_0)
 =\frac{1}{n_0}\log\#\cL(\Sigma,n_0)<\htop(\Sigma)+\varepsilon.
\]

We identify \(\mathcal{A}^n\times\mathcal{A}^n\)
with \((\mathcal{A}\times\mathcal{A})^n\) coordinatewise. Let
\[
 H
 =\left[(\mathcal{A}^{n_0}\times\mathcal{A}^{n_0})
 \setminus(\cL(\Sigma,n_0)\times\cL(\Sigma,n_0))\right]
 \cup\bigcup_{n\in\N}\left[
 (\mathcal{F}_n\times\cL(\Sigma_{n_0},n))
 \cup(\cL(\Sigma_{n_0},n)\times\mathcal{F}_n)
 \right].
\]
We claim that \(\Sigma_{H}=\mathcal{X}_\Sigma\).
For any \(z=((x_i,y_i))_{i\in\Z}\in\mathcal{X}_\Sigma\), we have
\(x=(x_i)_{i\in\Z}\in\Sigma\) and
\(y=(y_i)_{i\in\Z}\in\Sigma\) by the definition of \(\mathcal{X}_\Sigma\).
By the \(\sigma\)-invariance of \(\Sigma\), for any \(i\in\Z\),
\[
 z_{[i,i+n_0)}=(x_{[i,i+n_0)},y_{[i,i+n_0)})
 \in\cL(\Sigma,n_0)\times\cL(\Sigma,n_0).
\]
Since \(\Sigma=\Sigma_{\mathcal{F}}\), for any \(i\in\Z\)
and any \(n\in\N\), we have
\(x_{[i,i+n)},y_{[i,i+n)}\notin\mathcal{F}_n\), and hence
\(z_{[i,i+n)}=(x_{[i,i+n)},y_{[i,i+n)})
\notin(\mathcal{F}_n\times\cL(\Sigma_{n_0},n))
 \cup(\cL(\Sigma_{n_0},n)\times\mathcal{F}_n).\)
By the definition of \(H\), these two relations give
\(z_{[i,i+n)}\notin H\) for all \(i\in\Z\) and
\(n\in\N\),
and thus \(z\in\Sigma_{H}\).
Conversely, for any
\(z=((x_i,y_i))_{i\in\Z}\in\Sigma_{H}\), let
\(x=(x_i)_{i\in\Z}\) and \(y=(y_i)_{i\in\Z}\).
By the definition of \(H\), we have
\(x_{[i,i+n_0)}\), \(y_{[i,i+n_0)}\in\cL(\Sigma,n_0)\) for all \(i\in\Z\),
and hence \(x,y\in\Sigma_{n_0}\).
If \(x\notin\Sigma=\Sigma_{\mathcal{F}}\), then there exist
\(i\in\Z\) and \(n\in\N\) such that
\(x_{[i,i+n)}\in\mathcal{F}_n\).
Since \(y\in\Sigma_{n_0}\), this gives
\[
 z_{[i,i+n)}
 =(x_{[i,i+n)},y_{[i,i+n)})
 \in\mathcal{F}_n\times\cL(\Sigma_{n_0},n)\subseteq H,
\]
contradicting \(z\in\Sigma_{H}\).
Thus \(x\in\Sigma\). Similarly, it can be verified that
\(y\in\Sigma\). Hence \(z\in\mathcal{X}_\Sigma\), proving the claim.

For any \(n>n_0\), the definition of \(H\) gives
\[
 \#(H\cap(\mathcal{A}\times\mathcal{A})^n)
 \leq2\cdot\#\mathcal{F}_n\cdot\#\cL(\Sigma_{n_0},n).
\]
This, together with the claim, the definition of \(\hG\), and \mbox{\eqref{eq:language-limit}}, implies that
\[
 \begin{aligned}
 \hG(\mathcal{X}_\Sigma)
 &\leq\limsup_{n\to\infty}\frac{1}{n}
       \log\#(H\cap(\mathcal{A}\times\mathcal{A})^n)\\
 &\leq\limsup_{n\to\infty}\frac{1}{n}\log\#\mathcal{F}_n+\htop(\Sigma_{n_0})\\
 &\leq\limsup_{n\to\infty}\frac{1}{n}\log\#\mathcal{F}_n
       +\htop(\Sigma)+\varepsilon.
 \end{aligned}
\]
Taking the infimum over all sets \(\mathcal{F}\) generating \(\Sigma\)
and then letting \(\varepsilon\searrow0\) gives
\(\hG(\mathcal{X}_\Sigma)\leq\hG(\Sigma)+\htop(\Sigma)\).
If \(\hG(\Sigma)<\htop(\Sigma)\), this yields
\(\hG(\mathcal{X}_\Sigma)<2\htop(\Sigma)=\htop(\mathcal{X}_\Sigma)\).
\end{proof}

Given a finite or countable state space \(S\) and a probability vector
\(\mathbf{q}=(q_a)_{a\in S}\), the \emph{Bernoulli scheme} with probability
vector \(\mathbf{q}\) is \((S^\Z,\nu_{\mathbf{q}},\sigma)\), where
\(\sigma\) is the left-shift and \(\nu_{\mathbf{q}}\) is the product
probability measure determined by
\(\nu_{\mathbf{q}}([w]_{S^\Z})=\prod_{j=1}^{n}q_{w_j}\)
for any \(w=w_1\cdots w_n\in S^n\), \(n\in\N\)
(\cite[Section~2.2.3]{Climenhaga2018}).

When \(\Sigma\) is a two-sided subshift over a finite alphabet \(\mathcal{A}\),
an invariant probability measure \(\mu\) on \(\Sigma\)
is \emph{Bernoulli up to a finite period}
(\cite[Section~2.2.3]{Climenhaga2018})
if there exist \(p\in\N\) and
disjoint measurable sets \(\Sigma_0,\ldots,\Sigma_{p-1}\subseteq\Sigma\)
such that \(\mu(\bigcup_{i=0}^{p-1}\Sigma_i)=1\),
\(\sigma(\Sigma_i)=\Sigma_{i+1\bmod p}\), and
\((\Sigma_i,\text{\mbox{\(p\,\mu\)}}|_{\Sigma_i},\sigma^p|_{\Sigma_i})\) is
measure-theoretically isomorphic to a Bernoulli scheme for any
\(i\in\{0,\ldots,p-1\}\).

For any invariant probability measure \(\mu\) on \(\Sigma\) and any
\(n\in\N\), let \(\mathbb{P}^{(1)}_{\mu,n}\) and
\(\mathbb{P}^{(2)}_{\mu,n}\) be the probability vectors on
\(\cL(\Sigma,n)^2\) given by
\[
\mathbb{P}^{(1)}_{\mu,n}(u,v)=\mu([uv]_\Sigma) \quad\text{and}\quad
\mathbb{P}^{(2)}_{\mu,n}(u,v)=\mu([u]_\Sigma)\mu([v]_\Sigma),
\]
for any \(u,v\in\cL(\Sigma,n)\).
For any \(T\subseteq\cL(\Sigma,n)^2\), let
\(\mathbb{P}^{(1)}_{\mu,n}(T)
=\sum_{(u,v)\in T}\mathbb{P}^{(1)}_{\mu,n}(u,v)\) and
\(\mathbb{P}^{(2)}_{\mu,n}(T)
=\sum_{(u,v)\in T}\mathbb{P}^{(2)}_{\mu,n}(u,v)\).
By monotonicity and the \(\sigma\)-invariance of \(\mu\), we have
\begin{equation}\label{eq:joint-block-pointwise}
 0\leq\mathbb{P}^{(1)}_{\mu,n}(u,v)\leq\min\{\mu([u]_\Sigma),\mu([v]_\Sigma)\}
 \quad(u,v\in\cL(\Sigma,n)).
\end{equation}

For probability measures \(\mathbb{P}_1\) and \(\mathbb{P}_2\) on a common measurable
space \((X,\mathscr{B})\), their \emph{total variation distance} is
\(\TV(\mathbb{P}_1,\mathbb{P}_2)=\sup_{B\in\mathscr{B}}|\mathbb{P}_1(B)-\mathbb{P}_2(B)|\)
(\cite[p.~424, (1)]{GibbsSu2002}).
For any two probability measures \(\mathbb{P}_1,
\mathbb{P}_2\) on the same finite set \(X\), the total
variation distance satisfies
(\cite[p.~48, (4.2), and p.~51, (4.13)]{LevinPeres2017})
\begin{equation}
\label{eq:total-variation}
 \TV(\mathbb{P}_1,\mathbb{P}_2)=\frac{1}{2}\sum_{\alpha\in X}|\mathbb{P}_1(\alpha)-\mathbb{P}_2(\alpha)|\quad\text{and}\quad
 \sum_{\alpha\in X}\min\{\mathbb{P}_1(\alpha),\mathbb{P}_2(\alpha)\}=1-\TV(\mathbb{P}_1,\mathbb{P}_2).
\end{equation}

The following lemma is a direct consequence of
\cite[Theorem~2.2 and the subsequent note, p.~292]{Berbee1986}.
\begin{citedlemma}
\label{lem:berbee-symbolic}
Let \(\mu\) be an ergodic invariant probability measure on a subshift
\(\Sigma\) over a finite alphabet \(\mathcal{A}\).
Define \(\pi_-:\Sigma\to\mathcal{A}^{\Z\setminus\N_0}\) and
\(\pi_+:\Sigma\to\mathcal{A}^{\N_0}\) by
\[
 \pi_-(x)=(x_j)_{j<0}\quad\text{and}\quad
 \pi_+(x)=(x_j)_{j\geq0}
 \quad(x\in\Sigma).
\]
Let \(\mathbb{P}\) and \(\mathbb{Q}\) be the probability measures on
\(\pi_-(\Sigma)\times\pi_+(\Sigma)\) given by
\[
 \begin{aligned}
 \mathbb{P}(A\times B)
 &=\mu(\pi_-^{-1}(A)\cap\pi_+^{-1}(B)),\\
 \mathbb{Q}(A\times B)
 &=\mu(\pi_-^{-1}(A))\,\mu(\pi_+^{-1}(B)),
 \end{aligned}
\]
for Borel sets \(A\subseteq\pi_-(\Sigma)\) and
\(B\subseteq\pi_+(\Sigma)\).
If \(\TV(\mathbb{P},\mathbb{Q})<1\), then there exist
\(p\in\N\) and pairwise disjoint Borel sets
\(\Sigma_0,\ldots,\Sigma_{p-1}\subseteq\Sigma\)
with the following properties:
\begin{enumerate}
\item[(i)] For any \(i\in\{0,\ldots,p-1\}\),
\[
 \mu\left(\Sigma\setminus\bigcup_{j=0}^{p-1}\Sigma_j\right)=0,
 \quad \mu(\Sigma_i)=\frac{1}{p},\quad\text{and}\quad
 \mu(\sigma(\Sigma_i)\mathbin{\triangle}
                 \Sigma_{i+1\bmod p})=0.
\]
\item[(ii)] For any \(i\in\{0,\ldots,p-1\}\), the probability measure
\(\mu_i:=\text{\mbox{\(p\,\mu\)}}|_{\Sigma_i}\) satisfies
\[
 \mu_i(\sigma^{-p}(E))=\mu_i(E)
 \quad\text{for any Borel set }E\subseteq\Sigma.
\]
\item[(iii)] For any \(i\in\{0,\ldots,p-1\}\),
\[
 \lim_{n\to\infty}\sup_{r,s\in\N}\frac{1}{2}
 \sum_{\substack{u\in\cL(\Sigma,r)\\v\in\cL(\Sigma,s)}}
 \left|\mu_i(\sigma^r([u]_\Sigma)\cap\sigma^{-n}([v]_\Sigma))
 -\mu_i(\sigma^r([u]_\Sigma))\,
 \mu_i(\sigma^{-n}([v]_\Sigma))\right|=0.
\]
\end{enumerate}
\end{citedlemma}

Let \(T\) be an invertible measure-preserving transformation on
a standard Borel probability space \((X,\mathscr{B},\mu)\).
For a finite measurable partition \(\mathscr{P}=\{P_1,\ldots,P_N\}\)
and integers \(m<n\), write
\(\mathscr{P}_m^n=\bigvee_{j=m}^{n}T^{-j}(\mathscr{P})\).
Following \cite[Definition~3.5]{Sarig2011}, \(\mathscr{P}\) is called
\emph{weak Bernoulli} if, for any \(\varepsilon>0\), there exists
\(k\in \N\) such that
\[
 \sum_{A\in\mathscr{P}_{-n}^{0}}\sum_{B\in\mathscr{P}_k^{k+n}}
 |\mu(A\cap B)-\mu(A)\mu(B)|<\varepsilon
 \quad(n\in\N).
\]
For a family \(\mathscr{F}\) of subsets of the set \(X\), let
\(\mathscr{A}(\mathscr{F})\) denote the smallest \(\sigma\)-algebra on \(X\)
containing \(\mathscr{F}\).

The following criterion follows from the results of
Friedman and Ornstein~\cite{FriedmanOrnstein1970}
and Ornstein~\cite[Theorem~2]{Ornstein1970Infinite}
(see also \cite{Ornstein1970Factors} and \cite[Section~3.3]{Sarig2011}).

\begin{citedlemma}[{\cite[Section~3.3]{Sarig2011}}]
\label{lem:friedman-ornstein}
Let \((X,\mathscr{B},\mu)\) be a standard Borel probability space, and let
\(T\) be an invertible measure-preserving transformation on \(X\).
If there exists a sequence \((\mathscr{P}_n)_{n\in\N_0}\) of finite
measurable partitions satisfying the following conditions\footnote{Sarig~\cite[Section~3.3]{Sarig2011}
uses the term ``generating increasing sequence of weak Bernoulli
partitions'' for such a sequence.}:
\begin{enumerate}
\item[(i)] For any \(n\in\N_0\),
\(\mathscr{P}_{n+1}\) is a refinement of \(\mathscr{P}_n\);
\item[(ii)] \(\mathscr{A}(\bigcup_{n\in\N_0}\mathscr{P}_n)=\mathscr{B}\);
\item[(iii)] For any \(n\in\N_0\), \(\mathscr{P}_n\) is weak Bernoulli for \(T\);
\end{enumerate}
then \((X,\mu,T)\) is measure-theoretically isomorphic to a Bernoulli scheme.
\end{citedlemma}

\section{The zero-entropy case}
\label{sec:zero-entropy}

In this section, we prove Theorem~\ref{thm:finiteness}
in the zero-entropy case. We first show that all invariant probability
measures are supported on a common finite union of periodic orbits.

\begin{lemma}
\label{lem:zero-entropy}
Let \(\Sigma\) be a subshift over a finite alphabet \(\mathcal{A}\) satisfying
\(\hG(\Sigma)<\htop(\Sigma)=0.\)
Then the following statements hold:
\begin{enumerate}
\item[(i)] There exist periodic orbits
\(\mathcal{O}_{1},\ldots,\mathcal{O}_{s_{1}}\) such that
\(\mu\left(\bigcup_{a=1}^{s_{1}}\mathcal{O}_{a}\right)=1\)
for any invariant probability measure \(\mu\) on \(\Sigma\).
\item[(ii)] \(\Sigma\) has only finitely many ergodic measures of maximal entropy.
\end{enumerate}
\end{lemma}

\begin{proof}
(i)
Since \(\hG(\Sigma)<0\), there exists a set of forbidden words
\(\mathcal{F}\) generating \(\Sigma\) such that
\(\limsup_{n\to\infty}\frac{1}{n}\log\#\mathcal{F}_n<0\).
If \(\mathcal{F}_n\neq\varnothing\) for infinitely many \(n\in\N\),
then \(\#\mathcal{F}_n\geq1\) for those \(n\), contradicting this inequality.
Thus \(\mathcal{F}_n=\varnothing\) for all sufficiently large \(n\),
so \(\mathcal{F}\) is finite.
Put \(r_0=\max(\{2\}\cup\{n\in\N:\mathcal{F}_n\neq\varnothing\})\).

\begin{claim}\label{clm:local-characterization}
\(\Sigma=\big\{x\in\mathcal{A}^{\Z}: x_{[i,i+r_{0})}\in\cL(\Sigma,r_{0})
\text{ for all } i\in\Z\big\}.\)
\end{claim}

\begin{proof}[Proof of Claim~\ref{clm:local-characterization}]
By the \(\sigma\)-invariance of \(\Sigma\), we have
\(\Sigma\subseteq\allowbreak\big\{x\in\mathcal{A}^{\Z}:\allowbreak
 x_{[i,i+r_{0})}\in\cL(\Sigma,r_{0})\allowbreak
 \text{ for all }{i\in\Z}\big\}.\)
For the reverse inclusion, suppose, to the contrary, that
\(\big\{x\in\mathcal{A}^{\Z}:\allowbreak
 x_{[i,i+r_{0})}\in\cL(\Sigma,r_{0})\allowbreak
 \text{ for all }{i\in\Z}\big\}\nsubseteq\Sigma.\)
Then there exists
\(x\in\mathcal{A}^{\Z}\setminus\Sigma\) such that
\(x_{[i,i+r_{0})}\in\cL(\Sigma,r_{0})\) for all \(i\in\Z\).
Since \(x\notin\Sigma=\Sigma_{\mathcal{F}}\), there exist
\(i\in\Z\) and \(n_0\in\N\) such that
\(x_{[i,i+n_0)}\mathrel{\in}\mathcal{F}_{n_0}\subseteq
\mathcal{A}^{n_0}\mathbin{\setminus}\cL(\Sigma,n_0)\).
By the definition of \(r_0\), we have \(n_0\leq r_0\).
However, \(x_{[i,i+r_{0})}\in\cL(\Sigma,r_{0})\) by assumption,
so its prefix satisfies
\(x_{[i,i+n_{0})}\in\cL(\Sigma,n_{0}),\)
contradicting \(x_{[i,i+n_{0})}\notin\cL(\Sigma,n_{0})\).
Therefore, \(\big\{x\in\mathcal{A}^{\Z}:\allowbreak
 x_{[i,i+r_{0})}\in\cL(\Sigma,r_{0})\allowbreak
 \text{ for all }{i\in\Z}\big\}\subseteq\Sigma.\)
\end{proof}

Fix any recurrent point \(x\in\Sigma\). Then there exists a
sequence \(k_{\ell}\in\N\) with \(k_{\ell}\nearrow\infty\)
such that \(\sigma^{k_{\ell}}(x)\rightarrow x\).

\begin{claim}\label{clm:unique-extension}
For any \(i,j\in\Z\),
\(x_{[i,i+r_{0}-1)}=x_{[j,j+r_{0}-1)}\) implies
\(x_{i+r_{0}-1}=x_{j+r_{0}-1}\).
\end{claim}

\begin{proof}[Proof of Claim~\ref{clm:unique-extension}]
Suppose, to the contrary, that there exist \(i_{0},j_{0}\in\Z\) such that
\(x_{[i_{0},i_{0}+r_{0}-1)}=x_{[j_{0},j_{0}+r_{0}-1)}\) and
\(x_{i_{0}+r_{0}-1}\neq x_{j_{0}+r_{0}-1}.\)
Since \(\sigma^{k_{\ell}}(x)\rightarrow x\) and \(k_{\ell}\nearrow\infty\),
there exists \(\ell_{0}\in\N\) with \(k_{\ell_{0}}\geq r_{0}\) such that
\(x_{[i_{0}+k_{\ell_{0}},i_{0}+k_{\ell_{0}}+r_{0}-1)}=x_{[i_{0},i_{0}+r_{0}-1)}\) and
\(x_{[j_{0}+k_{\ell_{0}},j_{0}+k_{\ell_{0}}+r_{0}-1)}=x_{[j_{0},j_{0}+r_{0}-1)}\).
Fix any \(m\in\N\). For any
\((a_{0},\ldots,a_{m-1})\in\{i_{0},j_{0}\}^{m}\), define
\(y^{(a_{0},\ldots,a_{m-1})}\in\mathcal{A}^{\Z}\) by
\begin{equation}
\label{eq:periodic-concatenation}
 y^{(a_{0},\ldots,a_{m-1})}_{[0,mk_{\ell_{0}})}
 =x_{[a_{0},a_{0}+k_{\ell_{0}})}
  \cdots x_{[a_{m-1},a_{m-1}+k_{\ell_{0}})},
\end{equation}
and extend periodically so that
\(y^{(a_{0},\ldots,a_{m-1})}_{i+m k_{\ell_{0}}}
=y^{(a_{0},\ldots,a_{m-1})}_{i}\) for all \(i\in\Z\).
Then each length-\(r_{0}\) subword of \(y^{(a_{0},\ldots,a_{m-1})}\)
is a factor of
\(x_{[i_{0},i_{0}+k_{\ell_{0}}+r_{0}-1)}\) or
\(x_{[j_{0},j_{0}+k_{\ell_{0}}+r_{0}-1)}\).
Indeed, for any \(n\in\Z\), there exists a unique \(s\in\Z\)
such that \(s k_{\ell_{0}}\leq n<(s+1)k_{\ell_{0}}\).
By \eqref{eq:periodic-concatenation} and
the periodicity of \(y^{(a_{0},\ldots,a_{m-1})}\),
there exist \(a,b\in\{i_{0},j_{0}\}\)
such that
\(y^{(a_{0},\ldots,a_{m-1})}_{[s k_{\ell_{0}},(s+1)k_{\ell_{0}})}
=x_{[a,a+k_{\ell_{0}})}\) and
\(y^{(a_{0},\ldots,a_{m-1})}_{[(s+1)k_{\ell_{0}},(s+2)k_{\ell_{0}})}
   =x_{[b,b+k_{\ell_{0}})}.\) We consider the following two cases:

\smallskip
\textbf{Case \ref{clm:unique-extension}.1.} If \(n+r_{0}\leq(s+1)k_{\ell_{0}}\), then
\([n,n+r_{0})\subseteq[s k_{\ell_{0}},(s+1)k_{\ell_{0}})\).
Hence \(y^{(a_{0},\ldots,a_{m-1})}_{[n,n+r_{0})}\) is a factor of
\(x_{[a,a+k_{\ell_{0}})}\), and therefore of
\(x_{[a,a+k_{\ell_{0}}+r_{0}-1)}\).

\smallskip
\textbf{Case \ref{clm:unique-extension}.2.} If \(n+r_{0}>(s+1)k_{\ell_{0}}\), then
\((s+1)k_{\ell_{0}}<n+r_{0}
 \leq(s+1)k_{\ell_{0}}+r_{0}-1
 <(s+2)k_{\ell_{0}},\)
since \(n\leq(s+1)k_{\ell_{0}}-1\) and \(k_{\ell_{0}}\geq r_{0}\).
The choice of \(\ell_{0}\) and
\(x_{[i_{0},i_{0}+r_{0}-1)}=x_{[j_{0},j_{0}+r_{0}-1)}\) give
\(x_{[b,b+r_{0}-1)}
 =x_{[a,a+r_{0}-1)}
 =x_{[a+k_{\ell_{0}},a+k_{\ell_{0}}+r_{0}-1)},\)
implying that
\(y^{(a_{0},\ldots,a_{m-1})}_{[s k_{\ell_{0}},(s+1)k_{\ell_{0}}+r_{0}-1)}
   =x_{[a,a+k_{\ell_{0}})}x_{[b,b+r_{0}-1)}
   =x_{[a,a+k_{\ell_{0}}+r_{0}-1)}.\)
Since \([n,n+r_{0})\subseteq
[s k_{\ell_{0}},(s+1)k_{\ell_{0}}+r_{0}-1)\),
the word \(y^{(a_{0},\ldots,a_{m-1})}_{[n,n+r_{0})}\) is a factor of
\(x_{[a,a+k_{\ell_{0}}+r_{0}-1)}\) in this case as well.

\smallskip
Since \(n\in\Z\) is arbitrary and \(x\in\Sigma\),
the two cases imply that
\(y^{(a_{0},\ldots,a_{m-1})}_{[n,n+r_{0})}\in\cL(\Sigma,r_{0})\)
for all \(n\in\Z\).
By Claim~\ref{clm:local-characterization},
\(y^{(a_{0},\ldots,a_{m-1})}\in\Sigma\).
Since \(x_{i_{0}+r_{0}-1}\neq x_{j_{0}+r_{0}-1}\), the words
\(y^{(a_{0},\ldots,a_{m-1})}_{[0,mk_{\ell_{0}})}\) are pairwise distinct
as \((a_{0},\ldots,a_{m-1})\) ranges over \(\{i_{0},j_{0}\}^{m}\).
Hence \(\#\cL(\Sigma,mk_{\ell_{0}})\geq 2^{m}\) for all \(m\in\N\), and thus
\[
\htop(\Sigma)=\lim_{m\rightarrow\infty}
\frac{\log\#\cL(\Sigma,mk_{\ell_{0}})}{mk_{\ell_{0}}}
\geq\frac{\log 2}{k_{\ell_{0}}}>0,
\]
contrary to \(\htop(\Sigma)=0\).
\end{proof}

\begin{claim}\label{clm:bounded-period}
The point \(x\) is periodic with least period at most
\(\#\cL(\Sigma,r_{0}-1)\).
Consequently, every recurrent point of \(\Sigma\)
is periodic with least period at most \(\#\cL(\Sigma,r_{0}-1)\).
\end{claim}

\begin{proof}[Proof of Claim~\ref{clm:bounded-period}]
Since \(x\in\Sigma\), each word
\(x_{[j,j+r_{0}-1)}\) with \(0\leq j\leq\#\cL(\Sigma,r_{0}-1)\)
belongs to \(\cL(\Sigma,r_{0}-1)\).
There exist \(\#\cL(\Sigma,r_{0}-1)+1\) such starting positions,
whereas \(\cL(\Sigma,r_{0}-1)\) has only \(\#\cL(\Sigma,r_{0}-1)\) elements.
By the pigeonhole principle, there exist \(i_{1},j_{1}\in\Z\) with
\(0\leq i_{1}<j_{1}\leq\#\cL(\Sigma,r_{0}-1)\) such that
\(x_{[i_{1},i_{1}+r_{0}-1)}=x_{[j_{1},j_{1}+r_{0}-1)}.\)
By Claim~\ref{clm:unique-extension}, equal length-\((r_{0}-1)\) words in \(x\)
have equal following symbols.  Induction therefore
gives
\(x_{[i_{1}+n,i_{1}+n+r_{0}-1)}=x_{[j_{1}+n,j_{1}+n+r_{0}-1)}\)
for all \(n\in\N_{0}\).
In particular, \(x_{t+j_{1}-i_{1}}=x_{t}\) for all integers \(t\geq i_{1}\).
Since \(\sigma^{k_{\ell}}(x)\rightarrow x\) and \(k_{\ell}\nearrow\infty\),
for any \(t\in\Z\), there exists \(\ell_{1}(t)\in\N\) such that
\(t+k_{\ell_{1}(t)}\geq i_{1}\) and
\(x_{[t+k_{\ell_{1}(t)},\,t+j_{1}-i_{1}+k_{\ell_{1}(t)}+1)}
=x_{[t,\,t+j_{1}-i_{1}+1)}\).  Then
\(x_{t}=x_{t+k_{\ell_{1}(t)}}=x_{t+j_{1}-i_{1}+k_{\ell_{1}(t)}}
=x_{t+j_{1}-i_{1}}\).
Thus \(\sigma^{j_{1}-i_{1}}(x)=x\), where
\(1\leq j_{1}-i_{1}\leq\#\cL(\Sigma,r_{0}-1)\).
Since \(x\) is an arbitrary recurrent point of \(\Sigma\),
the second assertion follows.
\end{proof}

Set
\(P=\bigcup_{p=1}^{\#\cL(\Sigma,r_{0}-1)}
\{x\in\Sigma:\sigma^{p}(x)=x\}\).
Clearly, \(P\) is finite and \(\sigma\)-invariant.
Let \(\mu\) be any invariant probability measure on \(\Sigma\).
For any \(\ell\in\N\), the cylinder sets
\(C_{\ell,w}:=\{x\in\Sigma:x_{[-\ell,\ell)}=w\}\), \(w\in\cL(\Sigma,2\ell)\),
form a finite measurable partition of \(\Sigma\).
By Poincar\'e's Recurrence Theorem (\cite[Theorem~1.4]{Walters1982}),
for any \(\ell\in\N\) and any \(w\in\cL(\Sigma,2\ell)\)
with \(\mu(C_{\ell,w})>0\), there exists a measurable subset
\(G_{\ell,w}\subseteq C_{\ell,w}\) with
\(\mu(G_{\ell,w})=\mu(C_{\ell,w})\) such that every \(x\in G_{\ell,w}\)
satisfies \(\sigma^k(x)\in C_{\ell,w}\) for infinitely many \(k\in\N\).
For any \(\ell\in\N\), these subsets are pairwise disjoint and their
measures sum to \(1\). Thus, the set
\[
 G=\bigcap_{\ell\in\N}
 \bigcup_{\substack{w\in\cL(\Sigma,2\ell)\\\mu(C_{\ell,w})>0}}G_{\ell,w}
\]
satisfies \(\mu(G)=1\).
Clearly, every point of \(G\) is recurrent
and hence belongs to \(P\) by Claim~\ref{clm:bounded-period}.
Therefore \(\mu(P)=1\).
Writing \(P=\bigsqcup_{a=1}^{s_{1}}\mathcal{O}_{a}\)
as the disjoint union of its finitely many periodic orbits, we obtain
\(\mu(\bigcup_{a=1}^{s_{1}}\mathcal{O}_{a})=1\).

\smallskip

(ii)
Since \(\sigma\) is invertible and \(\mu\) is \(\sigma\)-invariant,
we have \(\mu(\{\sigma(x)\})
 =\mu(\sigma^{-1}(\{\sigma(x)\}))
 =\mu(\{x\})\) for any \(x\in\Sigma\).
Thus, the masses of the points
on each orbit \(\mathcal{O}_{a}\) are equal, so
\(\mu(\{x\})=\frac{\mu(\mathcal{O}_{a})}{\#\mathcal{O}_{a}}\)
for any \(a\in\{1,\ldots,s_{1}\}\) and any \(x\in\mathcal{O}_{a}\).
Since \(P=\bigsqcup_{a=1}^{s_{1}}\mathcal{O}_{a}\)
is finite and \(\mu(P)=1\),
\[
 \mu=\sum_{x\in P}\mu(\{x\})\delta_{x}
 =\sum_{a=1}^{s_{1}}\sum_{x\in\mathcal{O}_{a}}
 \mu(\{x\})\delta_{x}\\
 =\sum_{a=1}^{s_{1}}\mu(\mathcal{O}_{a})
 \frac{1}{\#\mathcal{O}_{a}}\sum_{x\in\mathcal{O}_{a}}\delta_{x}.
\]
If \(\mu\) is ergodic, then \(\mu(\mathcal{O}_{a})\in\{0,1\}\)
for any \(a\in\{1,\ldots,s_{1}\}\), since
\(\sigma^{-1}(\mathcal{O}_{a})=\mathcal{O}_{a}\). This, together with
\(\sum_{a=1}^{s_{1}}\mu(\mathcal{O}_{a})=\mu(P)=1\), implies that
exactly one of the coefficients
\(\mu(\mathcal{O}_{a})\) is \(1\).
Hence
\(\mu=\frac{1}{\#\mathcal{O}_{a}}
\sum_{x\in\mathcal{O}_{a}}\delta_{x}\)
for some \(a\in\{1,\ldots,s_{1}\}\).
Since \(\mu\) is an arbitrary invariant probability measure
and the periodic orbits
\(\mathcal{O}_{1},\ldots,\mathcal{O}_{s_{1}}\) do not depend on \(\mu\),
\(\Sigma\) has at most \(s_{1}\) ergodic invariant probability measures.
Hence,
\(\Sigma\) has only finitely many ergodic measures of maximal entropy.
\end{proof}

\begin{theorem}
\label{thm:zero-entropy}
Let \(\Sigma\) be a subshift over a finite alphabet
\(\mathcal{A}\) satisfying \(\hG(\Sigma)<\htop(\Sigma)=0\).
Then \(\Sigma\) has only finitely many ergodic measures of maximal entropy,
each of which is Bernoulli up to a finite period.
\end{theorem}

\begin{proof}
Lemma~\ref{lem:zero-entropy} gives only finitely many ergodic measures of maximal entropy,
each of which is the uniform probability measure on a periodic orbit.
Fix any ergodic measure of maximal entropy \(\mu\). Choose \(x\in\Sigma\) of least period
\(p\in\N\) such that \(\mu=\frac{1}{p}\sum_{i=0}^{p-1}\delta_{\sigma^i(x)}\).
Let \(\Sigma_i=\{\sigma^i(x)\}\) for \mbox{\(i=0,\ldots,p-1\)}.
These sets are measurable and pairwise disjoint, and
\(\mu(\Sigma\setminus\bigsqcup_{i=0}^{p-1}\Sigma_i)=0\).
Fix any \(i\in\{0,\ldots,p-1\}\). We have
\(\mu(\Sigma_i)=\frac{1}{p}\), \(\sigma(\Sigma_i)=\Sigma_{i+1\bmod p}\),
\(p\,\mu|_{\Sigma_i}=\delta_{\sigma^i(x)}\), and
\(\sigma^p|_{\Sigma_i}=\operatorname{id}_{\Sigma_i}\).
The map \(\sigma^i(x)\mapsto(0)_{n\in\Z}\) therefore gives a
measure-theoretic isomorphism from
\((\Sigma_i,p\,\mu|_{\Sigma_i},\sigma^p|_{\Sigma_i})\)
onto the Bernoulli scheme on \(\{0\}^{\Z}\) with the Dirac measure
\(\delta_{(0)_{n\in\Z}}\).
Thus, \(\mu\) is Bernoulli up to a finite period.
\end{proof}

\begin{remark}
\label{rem:zero-entropy-finite-type}
Under the hypotheses of \mbox{Theorem~\ref{thm:zero-entropy}},
\(\hG(\Sigma)<0\) gives a finite set of forbidden words
generating \(\Sigma\), so
\(\Sigma\) is of finite type.
Thus the finiteness of ergodic measures of maximal entropy also follows
from the classical result recorded by Cyr and Kra
\mbox{\cite[Lemma~2.4]{CyrKra2023}}.
Our proof uses word concatenations and recurrence to show directly
that all invariant probability measures are supported on a common finite
union of periodic orbits, without invoking that result.
\end{remark}

\section{The positive-entropy case}
\label{sec:positive-entropy}

In this section, we prove Theorem~\ref{thm:finiteness}
in the positive-entropy case. We establish finiteness using the uniform
separation estimate in Lemma~\ref{lem:positive-entropy-separation}
and obtain the Bernoulli conclusion by applying
Lemmas~\ref{lem:berbee-symbolic} and~\ref{lem:friedman-ornstein}.

For any subshift \(\Sigma\) over a finite alphabet
\(\mathcal{A}\),
any \(K\in\N\) and any family \((\mathcal{B}_n)_{n\geq K}\) with
\(\mathcal{B}_n\subseteq\cL(\Sigma,n)\) for all \(n\geq K\), define
\[
 U((\mathcal{B}_{n})_{n\geq K})=
 \{x\in\Sigma:\ \exists i,j\in\Z
       \text{ with } i\leq 0<j
       \text{ and } j-i\geq K
       \text{ such that } x_{[i,j)}\in\mathcal{B}_{j-i}\}.
\]

We first show that, for fixed \(0<t<\htop(\Sigma)\),
\(\mu(U((\mathcal{B}_n)_{n\geq K}))\) tends to zero as \(K\to\infty\),
uniformly over all ergodic measures of maximal entropy \(\mu\) and families satisfying
\(\#\mathcal{B}_n\leq e^{tn}\) for all \(n\geq K\).

\begin{lemma}\label{lem:uniform-coverage}
Let \(\Sigma\) be a subshift over a finite alphabet
\(\mathcal{A}\) with \(\htop(\Sigma)>0\).
For any \(0<t<\htop(\Sigma)\) and any
\(\delta>0\), there exists a \(K_0\in\N\) such that
\(\mu(U((\mathcal{B}_n)_{n\geq K}))\leq\delta\)
for any integer \(K\geq K_0\), any family
\((\mathcal{B}_n)_{n\geq K}\) satisfying
\(\mathcal{B}_n\subseteq\cL(\Sigma,n)\) and
\(\#\mathcal{B}_n\leq e^{tn}\) for all integers \(n\geq K\), and
any ergodic measure of maximal entropy \(\mu\) on \(\Sigma\).
\end{lemma}

\begin{proof}
Fix any \(0<t<\htop(\Sigma)\) and any \(\delta>0\).
Choose \(\eta>0\) such that
\(\frac{3\eta}{\htop(\Sigma)+\eta-t}<\frac{\delta}{2}\).
By \eqref{eq:language-limit}, there exists a constant
\(C_\eta\geq1\) such that
\begin{equation}\label{eq:language-growth-bound}
 \#\cL(\Sigma,n)\leq C_\eta\cdot e^{(\htop(\Sigma)+\eta)n}
 \quad(n\in\N_0).
\end{equation}
Since
\(\lim_{K\to\infty}\left[\left(1+\frac{2}{K}\right)H_2\left(\frac{2}{K+2}\right)
+\frac{\log C_\eta}{K}\right]=0\),
there exists \(K_0\in\N\) such that
\[
 \frac{3}{\htop(\Sigma)+\eta-t}
 \left(\eta+\left(1+\frac{2}{K}\right)H_2\left(\frac{2}{K+2}\right)
 +\frac{\log C_\eta}{K}\right)<\delta\quad (K\geq K_0),
\]
where \(H_2(u)=-u\log u-(1-u)\log(1-u)\) for \(0\leq u\leq1\).

For any integer \(K\geq K_0\), any family
\((\mathcal{B}_n)_{n\geq K}\) satisfying
\(\mathcal{B}_n\subseteq\cL(\Sigma,n)\) and
\(\#\mathcal{B}_n\leq e^{tn}\) for all integers \(n\geq K\),
and any ergodic measure of maximal entropy \(\mu\) on \(\Sigma\), it suffices to prove
\begin{equation}
\label{eq:coverage-bound}
 \mu(U((\mathcal{B}_{n})_{n\geq K}))
 \leq
 \frac{3}{\htop(\Sigma)+\eta-t}
 \left(\eta+
     \left(1+\frac{2}{K}\right)H_{2}\left(\frac{2}{K+2}\right)
                    +\frac{\log C_{\eta}}{K}\right).
\end{equation}

By the definition of \(U((\mathcal{B}_n)_{n\geq K})\), we have
\(U((\mathcal{B}_n)_{n\geq K})
=\bigcup_{m=K}^{+\infty}\bigcup_{i=1-m}^{0}\bigcup_{w\in\mathcal{B}_m}
\{x\in\Sigma:x_{[i,i+m)}=w\}\),
which, being a countable union of cylinders, is open and hence Borel.
For any integer \(M\geq K\), define
\(U_M((\mathcal{B}_n)_{n\geq K})
=\bigcup_{m=K}^{M}\bigcup_{i=1-m}^{0}\bigcup_{w\in\mathcal{B}_m}
\{x\in\Sigma:x_{[i,i+m)}=w\}\).
Then
\(U_M((\mathcal{B}_n)_{n\geq K})\nearrow U((\mathcal{B}_n)_{n\geq K})\)
as \(M\to\infty\).

For any integer \(N>2M\) and any \(x\in\Sigma\), put
\[
 \mathscr{I}^{(N)}(x)=
 \left\{[i,j):i,j\in\Z,\ 0\leq i<j\leq N,\
 K\leq j-i\leq M,\ x_{[i,j)}\in\mathcal{B}_{j-i}\right\},
\]
and
\(c_N(x)=\#((\bigcup\mathscr{I}^{(N)}(x))\cap\Z)\).
Then
\(\{j\in[M,N-M)\cap\Z:
\sigma^j(x)\in U_M((\mathcal{B}_n)_{n\geq K})\}
\allowbreak\subseteq\bigcup_{[a,b)\in\mathscr{I}^{(N)}(x)}([a,b)\cap\Z)\).
Indeed, for any \(j\in[M,N-M)\cap\Z\) with
\(\sigma^j(x)\in U_M((\mathcal{B}_n)_{n\geq K})\), there exist integers
\(m,i\) with \(K\leq m\leq M\) and \(i\leq0<i+m\) (i.e.,
\(1-m\leq i\leq0\)) such that
\((\sigma^j(x))_{[i,i+m)}\in\mathcal{B}_m\).
Since \((\sigma^j(x))_k=x_{j+k}\) for any \(k\in\Z\),
we have \(x_{[j+i,j+i+m)}=(\sigma^j(x))_{[i,i+m)}\in\mathcal{B}_m\).
Set \(a=j+i\) and \(b=j+i+m\).  Then \(a\leq j<b\), because
\(i\leq0\) and \(i+m\geq1\).  Moreover, \(b-a=m\) gives
\(K\leq b-a\leq M\) and \(x_{[a,b)}\in\mathcal{B}_{b-a}\).
Finally, \(M\leq j<N-M\) implies
\(a=j+i\geq j+1-m\geq0\) and \(b=j+i+m\leq j+m\leq N\).
Thus \([a,b)\in\mathscr{I}^{(N)}(x)\) and
\(j\in[a,b)\cap\Z\), which proves the inclusion.
Consequently,
\begin{equation}
\label{eq:coverage-count-lower-bound}
 \begin{aligned}
 c_N(x)&=\#\left(
 \bigcup_{[a,b)\in\mathscr{I}^{(N)}(x)}
       ([a,b)\cap\Z)\right)
 \geq\sum_{j=M}^{N-M-1}
       \mathbf{1}_{U_M((\mathcal{B}_n)_{n\geq K})}(\sigma^j(x))\\
 &=\sum_{j=0}^{N-1}
       \mathbf{1}_{U_M((\mathcal{B}_n)_{n\geq K})}(\sigma^j(x))
       -\sum_{j=0}^{M-1}
       \mathbf{1}_{U_M((\mathcal{B}_n)_{n\geq K})}(\sigma^j(x))\\
 &\quad
       -\sum_{j=N-M}^{N-1}
       \mathbf{1}_{U_M((\mathcal{B}_n)_{n\geq K})}(\sigma^j(x))\\
 &\geq\sum_{j=0}^{N-1}
       \mathbf{1}_{U_M((\mathcal{B}_n)_{n\geq K})}(\sigma^j(x))-2M.
 \end{aligned}
\end{equation}

Since \(\mu\) is ergodic, the Birkhoff Ergodic Theorem
(\cite[Theorem~1.14]{Walters1982})
implies that, for \(\mu\)-a.e. \(x\in\Sigma\),
\[
 \lim_{N\to\infty}\frac{1}{N}\sum_{j=0}^{N-1}
 \mathbf{1}_{U_M((\mathcal{B}_n)_{n\geq K})}(\sigma^j(x))
 =\mu(U_M((\mathcal{B}_n)_{n\geq K})).
\]
This, together with \eqref{eq:coverage-count-lower-bound} and
\(\frac{2M}{N}\to0\) for fixed \(M\), implies
that, for \(\mu\)-a.e. \(x\in\Sigma\),
\[
 \liminf_{N\to\infty}\frac{c_N(x)}{N}
 \geq\lim_{N\to\infty}\left[
 \frac{1}{N}\sum_{j=0}^{N-1}
       \mathbf{1}_{U_M((\mathcal{B}_n)_{n\geq K})}(\sigma^j(x))
 -\frac{2M}{N}\right]
 =\mu(U_M((\mathcal{B}_n)_{n\geq K})).
\]
For any \(\zeta>0\), the Dominated Convergence Theorem
(\cite[Corollary~0.9.1]{Walters1982}) therefore gives
\begin{equation}\label{eq:coverage-probability-limit}
 \lim_{N\to\infty}\mu\left(\left\{x\in\Sigma:\frac{c_N(x)}{N}\geq
       \mu(U_M((\mathcal{B}_n)_{n\geq K}))-\zeta\right\}\right)
 =1.
\end{equation}

Fix any integer \(M\geq K\) and any \(\zeta>0\). We prove that
\begin{equation}\label{eq:truncated-coverage-bound}
 \mu(U_M((\mathcal{B}_n)_{n\geq K}))
 \leq\zeta+\frac{3}{\htop(\Sigma)+\eta-t}
 \left(\eta+\left(1+\frac{2}{K}\right)H_2\left(\frac{2}{K+2}\right)
 +\frac{\log C_\eta}{K}\right).
\end{equation}

If \(\mu(U_M((\mathcal{B}_n)_{n\geq K}))\leq\zeta\),
\eqref{eq:truncated-coverage-bound} is immediate.
We may therefore assume that
\(\mu(U_M((\mathcal{B}_n)_{n\geq K}))>\zeta\).

For any integer \(N>2M\), set
\[
 \mathcal{C}_N(\zeta)=\left\{x_{[0,N)}:x\in\Sigma\text{ and }
 \frac{c_N(x)}{N}\geq\mu(U_M((\mathcal{B}_n)_{n\geq K}))-\zeta\right\},
\]
and for any \(r\in\{1,\ldots,\lfloor \frac{N}{K}\rfloor\}\), define
\begin{equation}\label{eq:coverage-tuple-definition}
 \begin{aligned}
 \mathcal{T}_{N,r}(\zeta)=\Big\{&
 (i_1,j_1,\ldots,i_r,j_r)\in\Z^{2r}:\\
 &0\leq i_k<j_k\leq N,\quad K\leq j_k-i_k\leq M
       \quad(1\leq k\leq r),\\
 &j_k\leq i_{k+1}\quad(1\leq k<r),\\
 &\sum_{k=1}^{r}(j_k-i_k)\geq
 \frac{N}{3}(\mu(U_M((\mathcal{B}_n)_{n\geq K}))-\zeta)
 \Big\}.
 \end{aligned}
\end{equation}
Fix any integer \(N>2M\).
For any \(x\in\Sigma\) with
\(\frac{c_N(x)}{N}\geq\mu(U_M((\mathcal{B}_n)_{n\geq K}))-\zeta\),
choose \(\mathscr{J}^{(N)}(x)\subseteq\mathscr{I}^{(N)}(x)\)
of minimum cardinality subject to
\(
 \bigcup\mathscr{J}^{(N)}(x)
 =\bigcup\mathscr{I}^{(N)}(x).
\)
By minimality, no member of \(\mathscr{J}^{(N)}(x)\) is contained in another.
Otherwise, removing such a member preserves the union and
reduces the cardinality, contradicting minimality. Write
\(\mathscr{J}^{(N)}(x)=\{[a_k,b_k):1\leq k\leq\#\mathscr{J}^{(N)}(x)\}\),
with \(a_1<\cdots<a_{\#\!\mathscr{J}^{(N)}(x)}\). Then
\(b_{k_1}<b_{k_2}\) for any integers
\(1\leq k_1<k_2\leq\#\!\mathscr{J}^{(N)}(x)\).
Otherwise, \(a_{k_1}<a_{k_2}\) and \(b_{k_2}\leq b_{k_1}\) give
\([a_{k_2},b_{k_2})\subseteq[a_{k_1},b_{k_1})\), contradicting the fact that
no member of \(\mathscr{J}^{(N)}(x)\) is contained in another.
Moreover, if \(1\leq k\leq\#\!\mathscr{J}^{(N)}(x)-2\) and \(a_{k+2}<b_k\), then
\([a_{k+1},b_{k+1})\subseteq[a_k,b_{k+2})=[a_k,b_k)\cup[a_{k+2},b_{k+2})\),
so removing \([a_{k+1},b_{k+1})\) again contradicts minimality.
Thus \(b_k\leq a_{k+2}\) for any integer \(1\leq k\leq\#\!\mathscr{J}^{(N)}(x)-2\), and hence
\([a_{k_1},b_{k_1})\cap[a_{k_2},b_{k_2})=\varnothing\)
whenever \(1\leq k_1<k_2\leq\#\!\mathscr{J}^{(N)}(x)\)
and \(k_1\equiv k_2\pmod2\).
Moreover,
\[
\begin{aligned}
&\max_{\varepsilon\in\{0,1\}}
 \left\{
 \sum_{\substack{1\leq k\leq\#\!\mathscr{J}^{(N)}(x)\\k\equiv\varepsilon\pmod2}}(b_k-a_k)
 \right\}
 \geq\frac{1}{2}\sum_{k=1}^{\#\!\mathscr{J}^{(N)}(x)}(b_k-a_k)
  \\
 &\quad\geq\frac{1}{2}\#\!\left(\bigcup_{k=1}^{\#\!\mathscr{J}^{(N)}(x)}([a_k,b_k)\cap\Z)\right)
 =\frac{c_N(x)}{2}>0.
\end{aligned}
\]
Choose \(\varepsilon\in\{0,1\}\) attaining this maximum and list the
corresponding intervals from left to right as \([i_1,j_1),\ldots,[i_r,j_r)\).
By the constructions of \(\mathscr{I}^{(N)}(x)\) and
\(\mathscr{J}^{(N)}(x)\),
\[
 \begin{cases}
 0\leq i_k<j_k\leq N,\quad K\leq j_k-i_k\leq M,\quad
 x_{[i_k,j_k)}\in\mathcal{B}_{j_k-i_k}\ (1\leq k\leq r),\\
 j_k\leq i_{k+1}\ (1\leq k<r),\quad
 rK\leq\sum_{k=1}^{r}(j_k-i_k)\leq N,\\
 \sum_{k=1}^{r}(j_k-i_k)
 \geq\frac{c_N(x)}{2}>\frac{c_N(x)}{3}
 \geq\frac{N}{3}(\mu(U_M((\mathcal{B}_n)_{n\geq K}))-\zeta),
 \end{cases}
\]
implying that \(1\leq r\leq\lfloor \frac{N}{K}\rfloor\) and
\((i_1,j_1,\ldots,i_r,j_r)\in\mathcal{T}_{N,r}(\zeta)\).
Since \(x\) is arbitrary subject to the defining condition of
\(\mathcal{C}_N(\zeta)\), the preceding construction gives
\[
 \mathcal{C}_N(\zeta)
 \subseteq\bigcup_{r=1}^{\lfloor \frac{N}{K}\rfloor}
 \ \bigcup_{\tau=(i_1,j_1,\ldots,i_r,j_r)\in\mathcal{T}_{N,r}(\zeta)}
 \{x_{[0,N)}:x\in\Sigma,\
 x_{[i_k,j_k)}\in\mathcal{B}_{j_k-i_k}\ (1\leq k\leq r)\}.
\]
For any \(r\in\{1,\ldots,\lfloor\frac{N}{K}\rfloor\}\)
and any \mbox{\(\tau=(i_1,j_1,\ldots,i_r,j_r)\in\mathcal{T}_{N,r}(\zeta)\)},
put \(j_0=0\) and \(i_{r+1}=N\).
For any \(x\in\Sigma\) with
\(x_{[i_k,j_k)}\in\mathcal{B}_{j_k-i_k}\) for \(1\leq k\leq r\), the
word \(x_{[0,N)}\) is uniquely
determined by its subwords \(x_{[i_k,j_k)}\), \(1\leq k\leq r\), and
\(x_{[j_k,i_{k+1})}\), \(0\leq k\leq r\).
Since \(x_{[i_k,j_k)}\in\mathcal{B}_{j_k-i_k}\)
for any integer \(1\leq k\leq r\) and
\(x_{[j_k,i_{k+1})}\in\cL(\Sigma,i_{k+1}-j_k)\) for any integer \(0\leq k\leq r\),
we obtain
\begin{equation}\label{eq:coverage-word-count}
 \#\mathcal{C}_N(\zeta)
 \leq\sum_{r=1}^{\lfloor \frac{N}{K}\rfloor}
 \ \sum_{\tau=(i_1,j_1,\ldots,i_r,j_r)\in\mathcal{T}_{N,r}(\zeta)}
 \left(\prod_{k=1}^{r}\#\mathcal{B}_{j_k-i_k}\right)
 \left(\prod_{k=0}^{r}\#\cL(\Sigma,i_{k+1}-j_k)\right).
\end{equation}
For any \(r\in\{1,\ldots,\lfloor \frac{N}{K}\rfloor\}\), the elements of
\(\mathcal{T}_{N,r}(\zeta)\) are nondecreasing \(2r\)-tuples
in \(\{0,\ldots,N\}\). Hence
\begin{equation}\label{eq:coverage-tuple-count}
 \#\mathcal{T}_{N,r}(\zeta)\leq\binom{N+2r}{2r}
 \leq e^{(N+2r)H_2\left(\frac{2r}{N+2r}\right)}
 \leq e^{N\left(1+\frac{2}{K}\right)H_2\left(\frac{2}{K+2}\right)},
\end{equation}
where the second inequality follows from
\cite[p.~353, Eq.~(11.40)]{CoverThomas2006}, expressed using natural logarithms;
the last follows from \(\frac{2r}{N}\leq\frac{2}{K}\) and the monotonicity of
\(z\mapsto(1+z)H_2\left(\frac{z}{1+z}\right)\) on \([0,\infty)\).
Moreover, for any \(\tau=(i_1,j_1,\ldots,i_r,j_r)\in\mathcal{T}_{N,r}(\zeta)\), we have
\(\sum_{k=0}^{r}(i_{k+1}-j_k)=N-\sum_{k=1}^{r}(j_k-i_k)\).
This, together with \(\#\mathcal{B}_n\leq e^{tn}\),
\eqref{eq:language-growth-bound}, \eqref{eq:coverage-tuple-definition},
and the inequalities \(r\leq\frac{N}{K}\),
\(C_\eta\geq1\), and
\(\htop(\Sigma)+\eta-t>0\), implies
\begin{equation}\label{eq:coverage-product-bound}
 \begin{aligned}
 &\left(\prod_{k=1}^{r}\#\mathcal{B}_{j_k-i_k}\right)
  \left(\prod_{k=0}^{r}\#\cL(\Sigma,i_{k+1}-j_k)\right)\\
 &\quad\leq
 e^{t\sum_{k=1}^{r}(j_k-i_k)}\cdot
 C_\eta^{r+1}\cdot e^{(\htop(\Sigma)+\eta)(N-\sum_{k=1}^{r}(j_k-i_k))}\\
 &\quad=
 C_\eta^{r+1}\cdot e^{(\htop(\Sigma)+\eta)N
 -(\htop(\Sigma)+\eta-t)\sum_{k=1}^{r}(j_k-i_k)}\\
 &\quad\leq C_\eta\cdot
 e^{N\left(\htop(\Sigma)+\eta+\frac{\log C_\eta}{K}\right)}\cdot
 e^{-(\htop(\Sigma)+\eta-t)\sum_{k=1}^{r}(j_k-i_k)}\\
 &\quad\leq C_\eta\cdot
 e^{N\left(\htop(\Sigma)+\eta+\frac{\log C_\eta}{K}\right)}\cdot
 e^{-\frac{N(\htop(\Sigma)+\eta-t)
 (\mu(U_M((\mathcal{B}_n)_{n\geq K}))-\zeta)}{3}}.
 \end{aligned}
\end{equation}
Substituting \eqref{eq:coverage-tuple-count} and
\eqref{eq:coverage-product-bound} into \eqref{eq:coverage-word-count} and using
\(\lfloor \frac{N}{K}\rfloor\leq N\) yields
\begin{equation}\label{eq:coverage-cardinality-bound}
 \#\mathcal{C}_N(\zeta)\leq N\cdot C_\eta\cdot e^{N\left[\htop(\Sigma)+\eta
 +\left(1+\frac{2}{K}\right)H_2\left(\frac{2}{K+2}\right)
 +\frac{\log C_\eta}{K}\right]}\cdot
 e^{-\frac{N(\htop(\Sigma)+\eta-t)
 (\mu(U_M((\mathcal{B}_n)_{n\geq K}))-\zeta)}{3}}.
\end{equation}
Since \(N\) is arbitrary, \eqref{eq:coverage-cardinality-bound}
holds for all integers \(N>2M\).

The Shannon--McMillan--Breiman theorem
(\cite[p.~93, Remark]{Walters1982}),
applied to the finite generating partition
\(\{[w]_\Sigma:w\in\cL(\Sigma,1)\}\), gives
\(\lim_{N\to\infty}
\left(-\frac{1}{N}\log\mu([x_{[0,N)}]_{\Sigma})\right)
=h_{\mu}(\sigma|_\Sigma)=\htop(\Sigma)\) for \(\mu\)-a.e. \(x\in\Sigma\).
Fix any \(\gamma>0\).
This implies that, for \(\mu\)-a.e. \(x\in\Sigma\), there exists
\(N(x)\in\N\) such that,
for all integers \(N\geq N(x)\),
\(-\frac{1}{N}\log\mu([x_{[0,N)}]_{\Sigma})\geq\htop(\Sigma)-\gamma\),
or equivalently
\(\mu([x_{[0,N)}]_{\Sigma})\leq e^{-N(\htop(\Sigma)-\gamma)}\).
Hence, for \(\mu\)-a.e. \(x\in\Sigma\),
\[
 \lim_{N\to\infty}\mathbf{1}_{\left\{y\in\Sigma:\,
 \mu([y_{[0,N)}]_{\Sigma})\leq e^{-N(\htop(\Sigma)-\gamma)}\right\}}(x)=1.
\]
The Dominated Convergence Theorem
(\cite[Corollary~0.9.1]{Walters1982}) therefore gives
\[
 \begin{aligned}
 &\lim_{N\to\infty}\mu\left(\left\{x\in\Sigma:
 \mu([x_{[0,N)}]_{\Sigma})\leq e^{-N(\htop(\Sigma)-\gamma)}
 \right\}\right)\\
 &\quad=\lim_{N\to\infty}\int_{\Sigma}
 \mathbf{1}_{\left\{y\in\Sigma:\,
 \mu([y_{[0,N)}]_{\Sigma})\leq e^{-N(\htop(\Sigma)-\gamma)}\right\}}(x)
 \,d\mu(x)\\
 &\quad=\int_{\Sigma}1\,d\mu=1.
 \end{aligned}
\]
This, together with \eqref{eq:coverage-probability-limit}
and the definition of \(\mathcal{C}_N(\zeta)\), implies that there exists
an integer \(N_0>2M\) such that, for all integers \(N\geq N_0\),
\(\mu(\{x\in\Sigma:x_{[0,N)}\in\mathcal{C}_N(\zeta)\})
\geq\frac{3}{4}\)
and \(\mu(\{x\in\Sigma:
\mu([x_{[0,N)}]_{\Sigma})\leq e^{-N(\htop(\Sigma)-\gamma)}\})
\geq\frac{3}{4}\), and thus
\[
 \begin{aligned}
 &\mu\left(\bigcup\left\{[w]_{\Sigma}: w\in\mathcal{C}_N(\zeta),\;
       \mu([w]_{\Sigma})\leq e^{-N(\htop(\Sigma)-\gamma)}\right\}\right)\\
 &\quad=\mu(\{x\in\Sigma:x_{[0,N)}\in\mathcal{C}_N(\zeta)\}
       \cap\{x\in\Sigma:
       \mu([x_{[0,N)}]_{\Sigma})\leq e^{-N(\htop(\Sigma)-\gamma)}\})\\
 &\quad\geq\mu(\{x\in\Sigma:x_{[0,N)}\in\mathcal{C}_N(\zeta)\})
 +\mu(\{x\in\Sigma:
       \mu([x_{[0,N)}]_{\Sigma})\leq e^{-N(\htop(\Sigma)-\gamma)}\})-1\\
 &\quad\geq\frac{3}{4}+\frac{3}{4}-1=\frac{1}{2}.
 \end{aligned}
\]
By subadditivity of \(\mu\),
\begin{equation}\label{eq:coverage-small-cylinder-count}
 \begin{aligned}
 \frac{1}{2}
 &\leq\mu\left(\bigcup\left\{[w]_{\Sigma}: w\in\mathcal{C}_N(\zeta),\;
       \mu([w]_{\Sigma})\leq e^{-N(\htop(\Sigma)-\gamma)}\right\}\right)\\
 &\leq\sum_{\substack{w\in\mathcal{C}_N(\zeta)\\
       \mu([w]_{\Sigma})\leq e^{-N(\htop(\Sigma)-\gamma)}}}
       \mu([w]_{\Sigma})
       \leq\#\mathcal{C}_N(\zeta)\cdot e^{-N(\htop(\Sigma)-\gamma)}.
 \end{aligned}
\end{equation}
Combining this with \eqref{eq:coverage-cardinality-bound}, we obtain
for all integers \(N\geq N_0\)
\[
 \begin{aligned}
 \htop(\Sigma)-\gamma-\frac{\log 2}{N}
 &\leq\frac{1}{N}\log\#\mathcal{C}_N(\zeta)\\
 &\leq\htop(\Sigma)+\eta
  +\left(1+\frac{2}{K}\right)H_2\left(\frac{2}{K+2}\right)
  +\frac{\log C_\eta}{K}\\
 &\quad+\frac{\log(NC_\eta)}{N}
  -\frac{\htop(\Sigma)+\eta-t}{3}
  (\mu(U_M((\mathcal{B}_n)_{n\geq K}))-\zeta),
 \end{aligned}
\]
and hence
\[
 \begin{aligned}
 &\frac{\htop(\Sigma)+\eta-t}{3}
  (\mu(U_M((\mathcal{B}_n)_{n\geq K}))-\zeta)\\
 &\qquad\leq\eta+\gamma
  +\left(1+\frac{2}{K}\right)H_2\left(\frac{2}{K+2}\right)
  +\frac{\log C_\eta}{K}+\frac{\log(2NC_\eta)}{N}.
 \end{aligned}
\]
Letting \(N\to\infty\) and then \(\gamma\searrow0\)
gives \eqref{eq:truncated-coverage-bound}.

Since \(U_M((\mathcal{B}_n)_{n\geq K})\nearrow U((\mathcal{B}_n)_{n\geq K})\),
letting \(\zeta\searrow0\) and then \(M\to\infty\) in
\eqref{eq:truncated-coverage-bound} yields \eqref{eq:coverage-bound}.
\end{proof}

Motivated by the replacement arguments of Garc\'ia-Ramos and
Pavlov~\mbox{\cite{GarciaRamosPavlov2019}}, we introduce the following
maps that replace one word by another of the same length.

For any subshift \(\Sigma\) over a finite alphabet
\(\mathcal{A}\),
any \(r\in\N\), and any pair of distinct
words \(u,v\in\cL(\Sigma,r)\), define \(\Phi^{\Sigma}_{u,v}: [u]_\Sigma\rightarrow
\mathcal{A}^{\Z}\) by
\begin{equation}\label{eq:replacement-map}
 \Phi^{\Sigma}_{u,v}(x)_{i}=
 \begin{cases}
 v_{i},&i\in \Z\cap [0, r),\\
 x_{i},&i\in\Z\setminus [0,r),
 \end{cases}
\end{equation}
and put \(D_{\Sigma}(u,v)=\{x\in[u]_\Sigma:
\Phi^{\Sigma}_{u,v}(x)\in\Sigma\}\).
It can be verified that the restriction
\(\Phi^{\Sigma}_{u,v}|_{D_{\Sigma}(u,v)}:D_{\Sigma}(u,v)\to D_{\Sigma}(v,u)\) is a homeomorphism
between compact sets, with inverse \(\Phi^{\Sigma}_{v,u}|_{D_{\Sigma}(v,u)}\).

We next combine Lemma~\ref{lem:uniform-coverage} with the
hypothesis \(\hG(\Sigma)<\htop(\Sigma)\) to establish
measure preservation on every Borel subset of
\(D_{\Sigma}(u,v)\), without assuming either inclusion
between the extender sets of \(u\) and \(v\).
For related results under inclusion assumptions, see
\mbox{\cite[Theorem~3.1]{Meyerovitch2013}} and
\mbox{\cite[Theorems~1.1 and~1.3]{Hedges2025}}.

\begin{lemma}\label{lem:finite-replacement-invariance}
Let \(\Sigma\) be a subshift over a finite alphabet
\(\mathcal{A}\) satisfying \(\htop(\Sigma)>0\) and
\(\hG(\Sigma)<\htop(\Sigma)\).
For any \(r\in\N\), any pair of distinct words
\(u,v\in\cL(\Sigma,r)\), and any invariant probability measure
\(\mu\) on \(\Sigma\) with
\(h_\mu(\sigma|_\Sigma)=\htop(\Sigma)\), we have
\(\mu(B)=\mu(\Phi^{\Sigma}_{u,v}(B))\) for any Borel set \(B\subseteq D_{\Sigma}(u,v)\).
\end{lemma}

\begin{proof}
Since \(\htop(\Sigma)>0\), we have
\(\#\mathcal{A}\geq2\). Put
\(s^*=\frac{1}{2}(\htop(\Sigma)+\max\{0,\hG(\Sigma)\})\), so that
\(\max\{0,\hG(\Sigma)\}<s^*<\htop(\Sigma)\).
Fix any \(r\in\N\), any pair of distinct words
\(u,v\in\cL(\Sigma,r)\), and any invariant probability measure
\(\mu\) on \(\Sigma\) with \(h_\mu(\sigma|_\Sigma)=\htop(\Sigma)\).

\smallskip

We first assume that \(\mu\) is ergodic and prove
\(\mu(D_{\Sigma}(u,v))\leq\mu([v]_\Sigma)\).

\smallskip

Fix \(0<\epsilon_0<\frac{1}{2}\) with
\(s^*+H_2(\epsilon_0)+\epsilon_0\log(\#\mathcal{A}-1)<\htop(\Sigma)\),
and put \(t_0=\frac{1}{2}(s^*+H_2(\epsilon_0)+\epsilon_0\log(\#\mathcal{A}-1)
+\htop(\Sigma))\), so that \(s^*+H_2(\epsilon_0)+\epsilon_0\log(\#\mathcal{A}-1)
<t_0<\htop(\Sigma)\).
Since \(\hG(\Sigma)<s^*\), there exists a set of forbidden words
\(\mathcal{F}\) generating \(\Sigma\) such that
\(\limsup_{n\to\infty}\frac{1}{n}\log\#\mathcal{F}_n<s^*\).
This and
\(\lim_{n\to\infty}\frac{\log(n+1)}{n}=0\) give \(n_0\in\N\)
such that, for all integers \(n\geq n_0\),
\[
 \#\mathcal{F}_n\leq e^{s^*n}
 \quad\text{and}\quad
 \frac{\log(n+1)}{n}
 \leq t_0-s^*-H_2(\epsilon_0)-\epsilon_0\log(\#\mathcal{A}-1).
\]
For any \(n\in\N\), let
\(\mathcal{B}_{n}=\bigcup_{f\in\mathcal{F}_n}
 \left\{w\in\cL(\Sigma,n):
 \#\{0\leq i<n:w_i\neq f_i\}\leq\epsilon_0 n\right\}.\)
For any integer \(n\geq n_0\) and any integer \(q\) with
\(0\leq q\leq\lfloor\epsilon_0 n\rfloor\), we have
\(0\leq\frac{q}{n}\leq\epsilon_0<\frac{1}{2}\).
Hence, by \cite[p.~353, Eq.~(11.40)]{CoverThomas2006}, expressed using natural logarithms,
and the monotonicity of \(H_2\) on \([0,\frac{1}{2}]\), we obtain
\[
 \binom{n}{q}\leq e^{nH_2(\frac{q}{n})}\leq e^{nH_2(\epsilon_0)}.
\]
Since \(\#\mathcal{A}\geq2\) and \(q\leq\epsilon_0 n\), we also have
\((\#\mathcal{A}-1)^q\leq e^{\epsilon_0 n\log(\#\mathcal{A}-1)}\).
Consequently, for any integer \(n\geq n_0\), the definition of \(\mathcal{B}_n\),
the preceding estimates, and \(\lfloor\epsilon_0 n\rfloor+1\leq n+1\) give
\[
 \begin{aligned}
 \#\mathcal{B}_{n}
 &\leq\sum_{f\in\mathcal{F}_n}
 \sum_{q=0}^{\lfloor\epsilon_0 n\rfloor}
 \binom{n}{q}(\#\mathcal{A}-1)^q\\
 &=\#\mathcal{F}_n\sum_{q=0}^{\lfloor\epsilon_0n\rfloor}
 \binom{n}{q}(\#\mathcal{A}-1)^{q}\\
 &\leq(n+1)\cdot
 e^{n(s^*+H_{2}(\epsilon_0)+\epsilon_0\log(\#\mathcal{A}-1))}
 \leq e^{t_{0}n}.
 \end{aligned}
\]
For any integer \(K\geq n_0\), put
\(V_{K}=\bigcup_{j=0}^{r-1}\sigma^{-j}
(U((\mathcal{B}_{n})_{n\geq K}))\).
By Lemma~\ref{lem:uniform-coverage}, subadditivity, and the \(\sigma\)-invariance of \(\mu\),
\[
 \mu(V_K)\leq\sum_{j=0}^{r-1}
 \mu\left(\sigma^{-j}(U((\mathcal{B}_n)_{n\geq K}))\right)
 =r\cdot\mu(U((\mathcal{B}_n)_{n\geq K}))\rightarrow 0
 \quad\text{as }K\to\infty.
\]
Thus \(\lim_{K\to\infty}\mu(V_K)=0\).

Suppose, to the contrary, that \(\mu(D_{\Sigma}(u,v))>\mu([v]_\Sigma)\).
By \(\lim_{K\to\infty}\mu(V_K)=0\), there exists an integer
\(K_0\geq\max\{n_0,\frac{4r}{\epsilon_0}\}\) such that
\(\mu(V_{K_0})<\mu(D_{\Sigma}(u,v))-\mu([v]_\Sigma)\).
Then \(\mu(D_{\Sigma}(u,v)\setminus V_{K_0})\geq\mu(D_{\Sigma}(u,v))-\mu(V_{K_0})>\mu([v]_\Sigma)\).

\begin{claim}\label{clm:simultaneous-replacement}
For any \(x\in\Sigma\) and any finite set \(J\subseteq\Z\) satisfying
\(\sigma^j(x)\in D_{\Sigma}(u,v)\setminus V_{K_0}\) for all \(j\in J\) and
\(|j-j'|\geq K_0+r\) for any distinct \(j,j'\in J\),
we have \(y(x,J)\in\Sigma\), where
\(y(x,J)\in\mathcal{A}^{\Z}\) is defined by
\[
 y(x,J)_{i}=
 \begin{cases}
 v_{i-j},&i\in[j,j+r)\text{ for some }j\in J,\\
 x_{i},&i\in\Z\setminus\bigcup_{j\in J}[j,j+r).
 \end{cases}
\]
\end{claim}

\begin{proof}[Proof of Claim~\ref{clm:simultaneous-replacement}]
Suppose, to the contrary, that there exist \(x\in\Sigma\)
and a finite set \(J\subseteq\Z\) satisfying the hypotheses of the claim
such that \(y(x,J)\notin\Sigma\).
Since \(\Sigma=\Sigma_{\mathcal{F}}\), there exist \(i_0\in\Z\)
and \(n\in\N\) such that
\(y(x,J)_{[i_0,i_0+n)}\in\mathcal{F}_n\).
This implies that the interval \([i_0,i_0+n)\) meets at least one
\([j,j+r)\), \(j\in J\).
Otherwise, by the definition of \(y(x,J)\), we have
\(y(x,J)_{[i_0,i_0+n)}=x_{[i_0,i_0+n)}\in\cL(\Sigma,n)\),
contradicting \(y(x,J)_{[i_0,i_0+n)}\in\mathcal{F}_n\).
We consider the following two cases:

\smallskip

\textbf{Case \ref{lem:finite-replacement-invariance}.1.} If \(n<K_0\), then for any \(j,j'\in J\) with \(j<j'\),
we have \(j'-(j+r)\geq K_0>n\), implying that
\([i_0,i_0+n)\) meets at most one interval \([j,j+r)\), \(j\in J\),
and therefore exactly one. For this \(j\), \eqref{eq:replacement-map}
shows that \(\sigma^{-j}(\Phi^{\Sigma}_{u,v}(\sigma^j(x)))\) is obtained from \(x\)
by replacing the word on \([j,j+r)\) with \(v\) and leaving all other
coordinates unchanged. Since \([i_0,i_0+n)\) meets no other interval
\([j',j'+r)\), \(j'\in J\setminus\{j\}\), the definition of \(y(x,J)\) gives
\(y(x,J)_{[i_0,i_0+n)}
 =(\sigma^{-j}(\Phi^{\Sigma}_{u,v}(\sigma^j(x))))_{[i_0,i_0+n)}\).
By the assumption on \(J\), we have
\(\sigma^j(x)\in D_{\Sigma}(u,v)\setminus V_{K_0}\subseteq D_{\Sigma}(u,v)\),
so the definition of \(D_{\Sigma}(u,v)\) gives \(\Phi^{\Sigma}_{u,v}(\sigma^j(x))\in\Sigma\).
By the \(\sigma\)-invariance of \(\Sigma\), we have
\(\sigma^{-j}(\Phi^{\Sigma}_{u,v}(\sigma^j(x)))\in\Sigma\), implying that
\(y(x,J)_{[i_0,i_0+n)}
 =(\sigma^{-j}(\Phi^{\Sigma}_{u,v}(\sigma^j(x))))_{[i_0,i_0+n)}
 \in\cL(\Sigma,n),\)
contradicting
\(y(x,J)_{[i_0,i_0+n)}\in\mathcal{F}_n\).

\smallskip

\textbf{Case \ref{lem:finite-replacement-invariance}.2.} If \(n\geq K_0\), let \(j_1<\cdots<j_m\) be all the
elements of \(\{j\in J:i_0-r<j<i_0+n\}\).
We have proved that \([i_0,i_0+n)\) meets at least one interval
\([j,j+r)\) with \(j\in J\).
Since \([j,j+r)\) meets \([i_0,i_0+n)\) if and only if
\(i_0-r<j<i_0+n\), it follows that
\(\{j\in J:i_0-r<j<i_0+n\}\neq\varnothing\), and hence \(m\geq1\).
By \(|j-j'|\geq K_0+r\) for any distinct \(j,j'\in J\), we have
\[
(m-1)(K_0+r)
 \leq\sum_{\ell=1}^{m-1}(j_{\ell+1}-j_\ell)
 =j_m-j_1<n+r,
\]
implying that
\(m<1+\frac{n+r}{K_0+r}<\frac{n}{K_0+r}+2.\)
By the definition of \(y(x,J)\), any coordinate
\(i\in[i_0,i_0+n)\) with \(x_i\neq y(x,J)_i\) belongs to one of the
intervals \([j_\ell,j_\ell+r)\), \(1\leq\ell\leq m\).
Each of these intervals contains \(r\) integer coordinates, so
\[
\#\{i\in[i_0,i_0+n):x_i\neq y(x,J)_i\}
 \leq rm\leq r\left(\frac{n}{K_0+r}+2\right).
\]
This, together with \(n\geq K_0\) and \(K_0\geq\frac{4r}{\epsilon_0}\), implies
\begin{equation}\label{eq:coordinate-difference-bound}
 \#\{i\in[i_0,i_0+n):x_i\neq y(x,J)_i\}
 \leq r\left(\frac{n}{K_0+r}+2\right)
 \leq r\left(\frac{n}{K_0}+\frac{2n}{K_0}\right)
 =\frac{3r}{K_0}n\leq\frac{3}{4}\epsilon_0n\leq\epsilon_0n.
\end{equation}
Since \(x_{[i_0,i_0+n)}\in\cL(\Sigma,n)\) and
\(y(x,J)_{[i_0,i_0+n)}\in\mathcal{F}_n\),
\eqref{eq:coordinate-difference-bound}
and the definition of \(\mathcal{B}_n\) give
\(x_{[i_0,i_0+n)}\in\mathcal{B}_n\).
Since \(\{j\in J:i_0-r<j<i_0+n\}\neq\varnothing\),
there exists \(j_0\in J\) such that \(i_0-r<j_0<i_0+n\).
Put \(k=\max\{j_0,i_0\}\in\Z\).
The inequalities \(i_0<j_0+r\), \(j_0<i_0+n\), and \(r,n\geq1\) give
\(k=\max\{j_0,i_0\}<\min\{j_0+r,i_0+n\}\), implying that
\(k\in[j_0,j_0+r)\cap[i_0,i_0+n)\).
Then \(0\leq k-j_0\leq r-1\), \(i_0-k\leq0<i_0+n-k\), and
\((\sigma^k(x))_{[i_0-k,i_0+n-k)}
 =x_{[i_0,i_0+n)}\in\mathcal{B}_n.\)
By \(x\in\Sigma\) and the \(\sigma\)-invariance of \(\Sigma\), we have
\(\sigma^k(x)\in\Sigma\).
Since \(n\geq K_0\), it follows that
\(\sigma^k(x)\in U((\mathcal{B}_{n})_{n\geq K_0}).\)
This, together with \(0\leq k-j_0\leq r-1\) and the definition of \(V_{K_0}\),
implies that
\(\sigma^{j_0}(x)=\sigma^{-(k-j_0)}(\sigma^{k}(x))\in\sigma^{-(k-j_0)}
 (U((\mathcal{B}_{n})_{n\geq K_0}))
 \subseteq V_{K_0}.\)
Since \(j_0\in J\), the assumption on \(J\) gives
\(\sigma^{j_0}(x)\in D_{\Sigma}(u,v)
\setminus V_{K_0}\), a contradiction.

Thus \(y(x,J)\in\Sigma\).
\end{proof}

Choose \(a,b \in\mathbb{R}\) such that
\(\mu(D_{\Sigma}(u,v)\setminus V_{K_0})>a>b>\mu([v]_\Sigma)\),
and put \(d=\frac{a-b}{4K_0+8r}\).
For any integer \(N\geq r\), define
\[
 \mathcal{W}_{N}=
 \left\{x_{[0,N)}:
 x\in\Sigma,\
 \frac{1}{N}\sum_{j=0}^{N-r}\mathbf{1}_{D_{\Sigma}(u,v)\setminus V_{K_0}}(\sigma^{j}(x))
 \geq a\quad\text{and}\quad
 \frac{1}{N}\sum_{j=0}^{N-r}\mathbf{1}_{[v]_{\Sigma}}(\sigma^{j}(x))\leq b\right\}.
\]
The Birkhoff Ergodic Theorem (\cite[Theorem~1.14]{Walters1982}) gives
\(\lim_{N\to\infty}\frac{1}{N}\sum_{j=0}^{N-r}\mathbf{1}_{D_{\Sigma}(u,v)\setminus V_{K_0}}(\sigma^j(x))
=\mu(D_{\Sigma}(u,v)\setminus V_{K_0})>a\) and
\(\lim_{N\to\infty}\frac{1}{N}\sum_{j=0}^{N-r}\mathbf{1}_{[v]_\Sigma}(\sigma^j(x))
=\mu([v]_\Sigma)<b\)
for \(\mu\)-a.e. \(x\in\Sigma\). The Dominated Convergence Theorem
(\cite[Corollary~0.9.1]{Walters1982}) therefore gives
\(\lim_{N\to\infty}\mu\left(\bigcup_{w\in\mathcal{W}_N}[w]_\Sigma\right)=1.\)

Fix any integer \(N\geq\max\{r,\lceil\frac{1}{d}\rceil\}\).
For each \(w\in\mathcal{W}_N\), fix some \(x^{(w)}\in\Sigma\) satisfying
\(x^{(w)}_{[0,N)}=w\) and the two inequalities defining \(\mathcal{W}_N\).
For any \(q\in\{1,\ldots,\lfloor dN\rfloor\}\), let \(\mathcal{J}_q(w)\)
be the family of sets \(J\subseteq\{0,\ldots,N-r\}\) such that
\(\#J=q\), \(\sigma^j(x^{(w)})\in D_{\Sigma}(u,v)\setminus V_{K_0}\) for all \(j\in J\),
and \(|j-j'|\geq K_0+r\) for any distinct \(j,j'\in J\).
Since \(\frac{1}{N}\sum_{j=0}^{N-r}\mathbf{1}_{D_{\Sigma}(u,v)
\setminus V_{K_0}}(\sigma^j(x^{(w)}))\geq a\),
we have
\begin{equation}\label{eq:replacement-singleton-count}
 \begin{aligned}
 \#\mathcal{J}_1(w)
 &=\#\{j\in\{0,\ldots,N-r\}:\sigma^j(x^{(w)})\in D_{\Sigma}(u,v)\setminus V_{K_0}\}\\
 &=\sum_{j=0}^{N-r}\mathbf{1}_{D_{\Sigma}(u,v)\setminus V_{K_0}}(\sigma^j(x^{(w)}))
 \geq aN.
 \end{aligned}
\end{equation}
For any integer \(q\) with \(1\leq q<\lfloor dN\rfloor\) and any
\(J\in\mathcal{J}_q(w)\), we have
\begin{equation}\label{eq:replacement-exclusion-bound}
 \begin{aligned}
 &\#\left(\bigcup_{j\in J}((j-K_0-r,j+K_0+r)\cap\Z)\right)\\
 &\quad\leq\sum_{j\in J}\#((j-K_0-r,j+K_0+r)\cap\Z)\\
 &\quad=\sum_{j\in J}(2(K_0+r)-1)
 =(2(K_0+r)-1)q.
 \end{aligned}
\end{equation}
For any \(J\in\mathcal{J}_q(w)\) and any \(j\in\{0,\ldots,N-r\}\),
the definition of \(\mathcal{J}_{q+1}(w)\) gives
\begin{equation}\label{eq:admissible-extension}
 \begin{aligned}
 J\cup\{j\}\in\mathcal{J}_{q+1}(w)
 &\iff
 \begin{cases}
 \sigma^j(x^{(w)})\in D_{\Sigma}(u,v)\setminus V_{K_0},\\
 |j-j'|\geq K_0+r\quad\text{for all }j'\in J
 \end{cases}\\
 &\iff
 \begin{cases}
 \sigma^j(x^{(w)})\in D_{\Sigma}(u,v)\setminus V_{K_0},\\
 j\notin\displaystyle\bigcup_{j'\in J}
 ((j'-K_0-r,j'+K_0+r)\cap\Z).
 \end{cases}
 \end{aligned}
\end{equation}
For any \(J'\in\mathcal{J}_{q+1}(w)\), the definition of
\(\mathcal{J}_q(w)\) gives \(J'\setminus\{j\}\in\mathcal{J}_q(w)\)
for all \(j\in J'\), and thus
\(\{(J,j):J\in\mathcal{J}_q(w),\ j\in\{0,\ldots,N-r\},\ J\cup\{j\}=J'\}
=\{(J'\setminus\{j\},j):j\in J'\}.\)
Summing over \(J'\in\mathcal{J}_{q+1}(w)\) and using
\eqref{eq:replacement-singleton-count},
\eqref{eq:replacement-exclusion-bound}, and \eqref{eq:admissible-extension}
yields
\begin{equation}\label{eq:replacement-recursion}
 \begin{aligned}
 &(q+1)\cdot\#\mathcal{J}_{q+1}(w)=\sum_{J'\in\mathcal{J}_{q+1}(w)}\#J'\\
 &\quad=\sum_{J\in\mathcal{J}_q(w)}
 \#\{j\in\{0,\ldots,N-r\}:J\cup\{j\}\in\mathcal{J}_{q+1}(w)\}\\
 &\quad=\sum_{J\in\mathcal{J}_q(w)}\#\left(
 \begin{aligned}
 &\{j\in\{0,\ldots,N-r\}:\sigma^j(x^{(w)})\in D_{\Sigma}(u,v)\setminus V_{K_0}\}\\
 &\qquad\setminus\bigcup_{j'\in J}((j'-K_0-r,j'+K_0+r)\cap\Z)
 \end{aligned}\right)\\
 &\quad\geq\sum_{J\in\mathcal{J}_q(w)}\left(
 \begin{aligned}
 &\#\{j\in\{0,\ldots,N-r\}:\sigma^j(x^{(w)})\in D_{\Sigma}(u,v)\setminus V_{K_0}\}\\
 &\qquad-\#\Bigl(\bigcup_{j'\in J}((j'-K_0-r,j'+K_0+r)\cap\Z)\Bigr)
 \end{aligned}\right)\\
 &\quad\geq(aN-(2(K_0+r)-1)q)\cdot \#\mathcal{J}_q(w)\\
 &\quad\geq(aN-2(K_0+r)q)\cdot \#\mathcal{J}_q(w).
 \end{aligned}
\end{equation}
The choice of \(d\) gives \(a-2(K_0+r)d>b+2rd>0\), so
\(aN-2(K_0+r)q>0\) for all integers \(0\leq q\leq\lfloor dN\rfloor\).
Thus, if \(\lfloor dN\rfloor\geq2\), starting from \(\#\mathcal{J}_1(w)\geq aN\)
and applying \eqref{eq:replacement-recursion} successively
for \(q=1,\ldots,\lfloor dN\rfloor-1\) gives
\begin{equation}\label{eq:replacement-family-lower-bound}
 \begin{aligned}
 \#\mathcal{J}_{\lfloor dN\rfloor}(w)
 &\geq\frac{aN-2(K_0+r)(\lfloor dN\rfloor-1)}{\lfloor dN\rfloor}
          \cdot\#\mathcal{J}_{\lfloor dN\rfloor-1}(w)\\
 &\geq\cdots\geq
 \left(\prod_{q=1}^{\lfloor dN\rfloor-1}
          \frac{aN-2(K_0+r)q}{q+1}\right)\cdot\#\mathcal{J}_1(w)\\
 &\geq aN\prod_{q=1}^{\lfloor dN\rfloor-1}
          \frac{aN-2(K_0+r)q}{q+1}\\
 &=\prod_{q=1}^{\lfloor dN\rfloor}
          \frac{aN-2(K_0+r)(q-1)}{q}\\
 &\geq
 \frac{(aN-2(K_0+r)\lfloor dN\rfloor)^{\lfloor dN\rfloor}}
 {\lfloor dN\rfloor!}.
 \end{aligned}
\end{equation}
If \(\lfloor dN\rfloor=1\), then
\(\#\mathcal{J}_1(w)\geq aN\geq aN-2(K_0+r)\), so the lower bound in
\eqref{eq:replacement-family-lower-bound} also holds in this case.

For any \(w\in\mathcal{W}_N\) and any
\(J\in\mathcal{J}_{\lfloor dN\rfloor}(w)\),
Claim~\ref{clm:simultaneous-replacement} gives \(y(x^{(w)},J)\in\Sigma\).
If \(0\leq j\leq N-r\) and \(y(x^{(w)},J)_{[j,j+r)}\neq x^{(w)}_{[j,j+r)}\), then
there exists \(k\in[j,j+r)\cap\Z\) such that
\(y(x^{(w)},J)_k\neq x^{(w)}_k\).
By the definition of \(y(x^{(w)},J)\), there exists \(j'\in J\) such that
\(k\in[j',j'+r)\). Thus
\(j-j'\leq k-j'\leq r-1\) and
\(j'-j\leq k-j\leq r-1,\)
implying that \(|j-j'|\leq r-1<r\).
Therefore,
\[
 \begin{aligned}
 &\{0\leq j\leq N-r:y(x^{(w)},J)_{[j,j+r)}=v\}\\
 &\quad=\{0\leq j\leq N-r:y(x^{(w)},J)_{[j,j+r)}=v,\ x^{(w)}_{[j,j+r)}=v\}\\
 &\qquad\cup\{0\leq j\leq N-r:y(x^{(w)},J)_{[j,j+r)}=v,\ x^{(w)}_{[j,j+r)}\neq v\}\\
 &\quad\subseteq\{0\leq j\leq N-r:x^{(w)}_{[j,j+r)}=v\}
 \cup\bigcup_{j'\in J}((j'-r,j'+r)\cap\Z).
 \end{aligned}
\]
This, together with \(\frac{1}{N}\sum_{j=0}^{N-r}\mathbf{1}_{[v]_\Sigma}(\sigma^j(x^{(w)}))\leq b\)
and \(\#J=\lfloor dN\rfloor\), implies
\begin{equation}\label{eq:replacement-occurrence-bound}
 \begin{split}
 &\#\{0\leq j\leq N-r:y(x^{(w)},J)_{[j,j+r)}=v\}\\
 &\quad\leq\#\{0\leq j\leq N-r:x^{(w)}_{[j,j+r)}=v\}
 +\#(\bigcup_{j'\in J}((j'-r,j'+r)\cap\Z))\\
 &\quad\leq\sum_{j=0}^{N-r}\mathbf{1}_{[v]_\Sigma}(\sigma^j(x^{(w)}))
 +\sum_{j'\in J}\#((j'-r,j'+r)\cap\Z)\\
 &\quad=\sum_{j=0}^{N-r}\mathbf{1}_{[v]_\Sigma}(\sigma^j(x^{(w)}))
 +\sum_{j'\in J}(2r-1)\leq bN+(2r-1)\cdot \#J\\
 &\quad=bN+(2r-1)\cdot \lfloor dN\rfloor
 \leq bN+2r\cdot \lfloor dN\rfloor.
 \end{split}
\end{equation}

Define \(\Psi_N:\{(w,J):w\in\mathcal{W}_N,\ J\in\mathcal{J}_{\lfloor dN\rfloor}(w)\}
\rightarrow\cL(\Sigma,N)\) by
\(\Psi_N(w,J)=y(x^{(w)},J)_{[0,N)}\).
The map \(\Psi_N\) is well defined, since
Claim~\ref{clm:simultaneous-replacement} gives \(y(x^{(w)},J)\in\Sigma\)
for every pair in its domain.
Fix any \(z\in\operatorname{im}\Psi_N\).
For any \((w,J)\in\Psi_N^{-1}(\{z\})\), we have
\(\#J=\lfloor dN\rfloor\) by the definition of
\(\mathcal{J}_{\lfloor dN\rfloor}(w)\), and \(z=\Psi_N(w,J)=y(x^{(w)},J)_{[0,N)}\).
For \mbox{any integer \(j\) with} \(0\leq j\leq N-r\), the inclusion \([j,j+r)\subseteq[0,N)\)
therefore gives \(z_{[j,j+r)}=y(x^{(w)},J)_{[j,j+r)}\).
Using \eqref{eq:replacement-occurrence-bound}, we obtain
\begin{equation}\label{eq:replacement-image-occurrence-bound}
 \#\{0\leq j\leq N-r:z_{[j,j+r)}=v\}
 =\#\{0\leq j\leq N-r:y(x^{(w)},J)_{[j,j+r)}=v\}
 \leq bN+2r\lfloor dN\rfloor.
\end{equation}
The definition of \(y(x^{(w)},J)\) also gives
\(z_{[j,j+r)}=v\) for all \(j\in J\). Hence
\(J\subseteq\{0\leq j\leq N-r:z_{[j,j+r)}=v\}\).
Moreover, by the definitions of \(\mathcal{J}_{\lfloor dN\rfloor}(w)\) and
\(D_{\Sigma}(u,v)\), we have
\(\sigma^j(x^{(w)})\in D_{\Sigma}(u,v)\setminus V_{K_0}\subseteq D_{\Sigma}(u,v)\subseteq[u]_\Sigma\)
for all \(j\in J\).
Since \(x^{(w)}_{[0,N)}=w\) and \(J\subseteq\{0,\ldots,N-r\}\), it follows that
\(w_{[j,j+r)}=x^{(w)}_{[j,j+r)}=(\sigma^j(x^{(w)}))_{[0,r)}=u\)
for all \(j\in J\).
For any \(i\in\{0,\ldots,N-1\}\setminus\bigcup_{j\in J}[j,j+r)\),
the definition of \(y(x^{(w)},J)\) gives
\(y(x^{(w)},J)_i=x^{(w)}_i\), and hence \(w_i=z_i\). Therefore,
\[
 w_i=
 \begin{cases}
 u_{i-j},&i\in[j,j+r)\text{ for some }j\in J,\\
 z_i,&i\notin\bigcup_{j\in J}[j,j+r),
 \end{cases}
 \quad(0\leq i<N),
\]
implying that
\((w,J),(w',J)\in\Psi_N^{-1}(\{z\})\Longrightarrow w=w'\).
Combining this with \eqref{eq:replacement-image-occurrence-bound} and
the facts that \(\#J=\lfloor dN\rfloor\) and
\(J\subseteq\{0\leq j\leq N-r:z_{[j,j+r)}=v\}\)
for any \((w,J)\in\Psi_N^{-1}(\{z\})\), we obtain
\[
 \#\Psi_N^{-1}(\{z\})
 \leq\binom{\#\{0\leq j\leq N-r:z_{[j,j+r)}=v\}}{\lfloor dN\rfloor}
 \leq\binom{\lceil bN+2r\lfloor dN\rfloor\rceil}{\lfloor dN\rfloor},
\]
implying that
\[
 \begin{aligned}
 &\#\{(w,J):w\in\mathcal{W}_N,\ J\in\mathcal{J}_{\lfloor dN\rfloor}(w)\}
 =\sum_{w\in\mathcal{W}_N}\#\mathcal{J}_{\lfloor dN\rfloor}(w)\\
 &\quad\leq\sum_{z\in\cL(\Sigma,N)}\#\Psi_N^{-1}(\{z\})
 \leq\#\cL(\Sigma,N)
 \binom{\lceil bN+2r\lfloor dN\rfloor\rceil}{\lfloor dN\rfloor}.
 \end{aligned}
\]
Using the preceding inequality and
\eqref{eq:replacement-family-lower-bound} yields
\begin{equation}\label{eq:replacement-language-bound}
 \begin{aligned}
 \#\cL(\Sigma,N)
 &\geq\frac{\sum_{w\in\mathcal{W}_N}\#\mathcal{J}_{\lfloor dN\rfloor}(w)}
 {\binom{\lceil bN+2r\lfloor dN\rfloor\rceil}{\lfloor dN\rfloor}}\\
 &\geq\#\mathcal{W}_N\,
 \frac{(aN-2(K_0+r)\lfloor dN\rfloor)^{\lfloor dN\rfloor}}
 {\lfloor dN\rfloor!\binom{\lceil bN+2r\lfloor dN\rfloor\rceil}{\lfloor dN\rfloor}}\\
 &\geq\#\mathcal{W}_N
 \left(\frac{aN-2(K_0+r)\lfloor dN\rfloor}
 {\lceil bN+2r\lfloor dN\rfloor\rceil}\right)^{\lfloor dN\rfloor}.
 \end{aligned}
\end{equation}

Fix any \(\gamma>0\).
By \(\lim_{N\to\infty}\mu(\bigcup_{w\in\mathcal{W}_N}[w]_\Sigma)=1\)
and an argument similar to that used to prove
\eqref{eq:coverage-small-cylinder-count},
there exists an integer
\(N_1\geq\max\{r,\lceil\frac{1}{d}\rceil\}\) such that,
for all integers \(N\geq N_1\),
\[
 \begin{aligned}
 \frac{1}{2}
 &\leq\mu\left(\bigcup\left\{[w]_{\Sigma}: w\in\mathcal{W}_N,\;
       \mu([w]_{\Sigma})\leq e^{-N(\htop(\Sigma)-\gamma)}\right\}\right)\\
 &\leq\sum_{\substack{w\in\mathcal{W}_N,\\
       \mu([w]_{\Sigma})\leq e^{-N(\htop(\Sigma)-\gamma)}}}
       \mu([w]_{\Sigma})
       \leq\#\mathcal{W}_N\cdot e^{-N(\htop(\Sigma)-\gamma)}.
 \end{aligned}
\]
Taking logarithms and using \(\mathcal{W}_N\subseteq\cL(\Sigma,N)\)
yields, for all integers \(N\geq N_1\),
\[\htop(\Sigma)-\gamma-\frac{\log 2}{N}
\leq\frac{1}{N}\log\#\mathcal{W}_N
\leq\frac{1}{N}\log\#\cL(\Sigma,N).
\]
By \eqref{eq:language-limit}, letting \(N\to\infty\)
and then \(\gamma\searrow0\) gives
\[
\lim_{N\to\infty}\frac{1}{N}\log\#\mathcal{W}_N=\htop(\Sigma).
\]
This, together with \eqref{eq:language-limit},
\eqref{eq:replacement-language-bound}, and
\(\lim_{N\to\infty}\frac{\lfloor dN\rfloor}{N}=d\), implies
\[
 \begin{aligned}
 0&=\lim_{N\to\infty}\left(
       \frac{1}{N}\log\#\cL(\Sigma,N)
       -\frac{1}{N}\log\#\mathcal{W}_N\right)\\
  &=\lim_{N\to\infty}\frac{1}{N}
       \log\frac{\#\cL(\Sigma,N)}{\#\mathcal{W}_N}
       \\
  &\geq\lim_{N\to\infty}\frac{\lfloor dN\rfloor}{N}
       \log\frac{aN-2(K_0+r)\lfloor dN\rfloor}
       {\lceil bN+2r\lfloor dN\rfloor\rceil}\\
  &=d\log\frac{a-2(K_0+r)d}{b+2rd}>0,
 \end{aligned}
\]
a contradiction. Thus \(\mu(D_{\Sigma}(u,v))\leq\mu([v]_\Sigma)\).

\smallskip

Since \(r,u,v\) are arbitrary, we have
\(\mu(D_{\Sigma}(\hat{u},\hat{v}))\leq\mu([\hat{v}]_\Sigma)\) for any pair of distinct words
\(\hat{u},\hat{v}\) of the same length occurring in \(\Sigma\).

\smallskip

We next prove that
\(\mu(E)=\mu(\Phi^{\Sigma}_{u,v}(E))\) for any Borel set \(E\subseteq D_{\Sigma}(u,v)\).

Let
\(\pi_r:\Sigma\to\mathcal{A}^{\Z\setminus[0,r)}\) be the projection
\(\pi_r(x)=x_{\Z\setminus[0,r)}\).
For any \(w\in\{u,v\}\) and any \(x,y\in[w]_\Sigma\),
\(\pi_r(x)=\pi_r(y)\) gives \(x_i=y_i\) for all \(i\in\Z\setminus[0,r)\),
while \(x_{[0,r)}=y_{[0,r)}=w\). Thus \(x=y\), so
\(\pi_r|_{[w]_\Sigma}\) is injective. Since \([w]_\Sigma\) is compact and
\(\mathcal{A}^{\Z\setminus[0,r)}\) is Hausdorff, the continuous injection
\(\pi_r|_{[w]_\Sigma}\) is a homeomorphism onto the compact set
\(\pi_r([w]_\Sigma)\). In particular,
\(\pi_r([u]_\Sigma)\cap\pi_r([v]_\Sigma)\) is compact and hence Borel.
For any \(w\in\{u,v\}\), define the finite Borel measure
\(\mu_w\) by \(\mu_w(B)=\mu([w]_\Sigma\cap\pi_r^{-1}(B))\)
for any Borel set \(B\subseteq\mathcal{A}^{\Z\setminus[0,r)}\).

For any finite words \(a',b'\), possibly empty, let
\[
 Q(a',b')=\{z\in\mathcal{A}^{\Z\setminus[0,r)}:
 z_{[-|a'|,0)}=a',\ z_{[r,r+|b'|)}=b'\}.
\]
To prove that
\(\mu_u(\pi_r([u]_\Sigma)\cap\pi_r([v]_\Sigma)\cap Q(a',b'))
\leq\mu_v(Q(a',b'))\), we consider the following two cases:

\smallskip

\textbf{Case \ref{lem:finite-replacement-invariance}.3.} \(\{a'ub',a'vb'\}\subseteq\cL(\Sigma,|a'|+r+|b'|)\).
We claim that
\begin{equation}\label{eq:replacement-domain-identity}
[u]_\Sigma\cap\pi_r^{-1}(\pi_r([u]_\Sigma)\cap\pi_r([v]_\Sigma)\cap Q(a',b'))
=\sigma^{|a'|}(D_{\Sigma}(a'ub',a'vb')).
\end{equation}
Indeed, for any
\(x\in[u]_\Sigma\cap\pi_r^{-1}(\pi_r([u]_\Sigma)\cap\pi_r([v]_\Sigma)\cap Q(a',b'))\),
the definitions of \([u]_\Sigma\) and \(Q(a',b')\) give
\((\sigma^{-|a'|}(x))_{[0,|a'|+r+|b'|)}
=x_{[-|a'|,r+|b'|)}=a'ub'\).
Since \(\sigma^{-|a'|}(x)\in\Sigma\) by the \(\sigma\)-invariance of \(\Sigma\),
we have \(\sigma^{-|a'|}(x)\in[a'ub']_\Sigma\), and thus
\(\Phi^{\Sigma}_{a'ub',a'vb'}(\sigma^{-|a'|}(x))
=\sigma^{-|a'|}(\Phi^{\Sigma}_{u,v}(x))\)
by \eqref{eq:replacement-map}.
Moreover, \(\pi_r(x)\in\pi_r([v]_\Sigma)\) implies that there exists
\(z\in[v]_\Sigma\) such that $\pi_{r}(z)=\pi_{r}(x)$, i.e.,
\(z_{\Z\setminus[0,r)}=x_{\Z\setminus[0,r)}\).
Since \(x\in[u]_\Sigma\) and \(z_{[0,r)}=v\), \eqref{eq:replacement-map} gives
\(\Phi^{\Sigma}_{u,v}(x)=z\in[v]_\Sigma\subseteq \Sigma\), implying that
\(\sigma^{-|a'|}(\Phi^\Sigma_{u,v}(x))\in\Sigma\), and hence
\(\Phi^{\Sigma}_{a'ub',a'vb'}(\sigma^{-|a'|}(x))\in\Sigma\).
This, together with \(\sigma^{-|a'|}(x)\in[a'ub']_\Sigma\)
and the definition of \(D_{\Sigma}(a'ub',a'vb')\), gives
\(\sigma^{-|a'|}(x)\in D_{\Sigma}(a'ub',a'vb')\) and hence
\(x=\sigma^{|a'|}(\sigma^{-|a'|}(x))\in\sigma^{|a'|}(D_{\Sigma}(a'ub',a'vb'))\).
Conversely, for any \(t\in D_{\Sigma}(a'ub',a'vb')\), put \(x=\sigma^{|a'|}(t)\).
By the definition of \(D_{\Sigma}(a'ub',a'vb')\) and \eqref{eq:replacement-map},
we have \(t\in[a'ub']_\Sigma\) and \(\Phi^{\Sigma}_{a'ub',a'vb'}(t)\in[a'vb']_\Sigma\), implying that
\(x=\sigma^{|a'|}(t)\in\Sigma\) and thus
\(x_{[-|a'|,r+|b'|)}=t_{[0,|a'|+r+|b'|)}=a'ub'\).
Hence \(x\in[u]_\Sigma\) and \(\pi_r(x)\in Q(a',b')\).
By \eqref{eq:replacement-map}, we also have
\(\Phi^{\Sigma}_{u,v}(x)=\sigma^{|a'|}(\Phi^{\Sigma}_{a'ub',a'vb'}(t))\).
This, together with \(\Phi^{\Sigma}_{a'ub',a'vb'}(t)\in[a'vb']_\Sigma\), implies that
\(\Phi^{\Sigma}_{u,v}(x)_{[0,r)}=v\) and, by the \(\sigma\)-invariance of \(\Sigma\),
\(\Phi^{\Sigma}_{u,v}(x)\in\Sigma\). Thus \(\Phi^{\Sigma}_{u,v}(x)\in[v]_\Sigma\).
Moreover, \(\Phi^\Sigma_{u,v}\) leaves all coordinates outside \([0,r)\) unchanged,
so \(\pi_r(\Phi^\Sigma_{u,v}(x))=\pi_r(x)\in\pi_r([v]_\Sigma)\).
Since \(x\in[u]_\Sigma\), \(\pi_r(x)\in\pi_r([v]_\Sigma)\),
and \(\pi_r(x)\in Q(a',b')\), we have
\(x\in[u]_\Sigma\cap\pi_r^{-1}(\pi_r([u]_\Sigma)\cap\pi_r([v]_\Sigma)\cap Q(a',b'))\),
proving \eqref{eq:replacement-domain-identity}.

Also, for any \(x\in\Sigma\),
the definitions of \([v]_\Sigma\), \(\pi_r\), and \(Q(a',b')\) show that
\(x\in[v]_\Sigma\cap\pi_r^{-1}(Q(a',b'))\) if and only if
\(x_{[-|a'|,r+|b'|)}=a'vb'\).
By the \(\sigma\)-invariance of \(\Sigma\),
this is equivalent to \(\sigma^{-|a'|}(x)\in[a'vb']_\Sigma\),
and hence equivalent to \(x\in\sigma^{|a'|}([a'vb']_\Sigma)\). Therefore,
\begin{equation}\label{eq:replacement-cylinder-identity}
[v]_\Sigma\cap\pi_r^{-1}(Q(a',b'))=\sigma^{|a'|}([a'vb']_\Sigma).
\end{equation}

Since \(u\neq v\), we have \(a'ub'\neq a'vb'\).
As \(a'ub',a'vb'\in\cL(\Sigma,|a'|+r+|b'|)\),
applying the inequality \(\mu(D_{\Sigma}(\hat{u},\hat{v}))\leq\mu([\hat{v}]_\Sigma)\) proved above with
\(\hat{u}=a'ub'\) and \(\hat{v}=a'vb'\) gives
\(\mu(D_{\Sigma}(a'ub',a'vb'))\leq\mu([a'vb']_\Sigma)\).
By the definition of \(\mu_u\), \eqref{eq:replacement-domain-identity},
and the \(\sigma\)-invariance of \(\mu\), we have
\(\mu_u(\pi_r([u]_\Sigma)\cap\pi_r([v]_\Sigma)\cap Q(a',b'))
=\mu([u]_\Sigma\cap\pi_r^{-1}(\pi_r([u]_\Sigma)\cap\pi_r([v]_\Sigma)\cap Q(a',b')))
=\mu(\sigma^{|a'|}(D_{\Sigma}(a'ub',a'vb')))
=\mu(D_{\Sigma}(a'ub',a'vb'))\).
Likewise, the definition of \(\mu_v\), \eqref{eq:replacement-cylinder-identity},
and the \(\sigma\)-invariance of \(\mu\) give
\(\mu_v(Q(a',b'))\allowbreak=\mu([v]_\Sigma\cap\pi_r^{-1}(Q(a',b')))
=\mu(\sigma^{|a'|}([a'vb']_\Sigma))
=\mu([a'vb']_\Sigma)\).
Therefore
\(\mu_u(\pi_r([u]_\Sigma)\cap\pi_r([v]_\Sigma)\cap Q(a',b'))
\leq\mu_v(Q(a',b'))\).

\smallskip

\textbf{Case \ref{lem:finite-replacement-invariance}.4.} \(\{a'ub',a'vb'\}\nsubseteq\cL(\Sigma,|a'|+r+|b'|)\).
Then
\(\pi_r([u]_\Sigma)\cap\pi_r([v]_\Sigma)\cap Q(a',b')=\varnothing\).
Indeed, suppose, to the contrary, that \(\pi_r([u]_\Sigma)\cap\pi_r([v]_\Sigma)\cap Q(a',b')\) contains a point \(z\).
Then there exist \(x\in[u]_\Sigma\) and \(y\in[v]_\Sigma\)
with \(\pi_r(x)=\pi_r(y)=z\).
This, together with \(z\in Q(a',b')\), gives
\(x_{[-|a'|,r+|b'|)}=a'ub'\) and \(y_{[-|a'|,r+|b'|)}=a'vb'\),
so \(\{a'ub',a'vb'\}\subseteq\cL(\Sigma,|a'|+r+|b'|)\),
contradicting
\(\{a'ub',a'vb'\}\nsubseteq\cL(\Sigma,|a'|+r+|b'|)\).
Thus \(\mu_u(\pi_r([u]_\Sigma)\cap\pi_r([v]_\Sigma)\cap Q(a',b'))
=0\leq\mu_v(Q(a',b'))\).

\smallskip

To extend the inequalities proved in Cases~\ref{lem:finite-replacement-invariance}.3 and~\ref{lem:finite-replacement-invariance}.4
to all Borel subsets of \(\mathcal{A}^{\Z\setminus[0,r)}\), define
\(\mathscr{F}=\{\bigcup_{(a',b')\in F}Q(a',b'):
 m,n\in\N_0,\ F\subseteq\mathcal{A}^m\times\mathcal{A}^n\}.\)
It can be verified that the following three conditions hold:

(i) \(\varnothing\in\mathscr{F}\), by taking \(F=\varnothing\)
in the definition of \(\mathscr{F}\).

(ii) \(B_1\cap B_2\in\mathscr{F}\) for any \(B_1,B_2\in\mathscr{F}\).
For each \(i\in\{1,2\}\), choose
\(m_i,n_i\in\N_0\) and \(F_i\subseteq\mathcal{A}^{m_i}\times\mathcal{A}^{n_i}\)
such that \(B_i=\bigcup_{(a',b')\in F_i}Q(a',b')\).
Put \(m=\max\{m_1,m_2\}\) and \(n=\max\{n_1,n_2\}\), and define
\(\hat{B}_i=\left\{(\alpha a',b'\beta):
 (a',b')\in F_i,\ \alpha\in\mathcal{A}^{m-m_i},\
 \beta\in\mathcal{A}^{n-n_i}\right\}\)
for each \(i\in\{1,2\}\).
Then \(\hat{B}_1,\hat{B}_2\subseteq\mathcal{A}^m\times\mathcal{A}^n\).
By the definitions of \(Q\) and \(\hat{B}_i\) (\(i\in\{1,2\}\)), for any
\(z\in\mathcal{A}^{\Z\setminus[0,r)}\), we have \(z\in B_i\) if and only if
\((z_{[-m,0)},z_{[r,r+n)})\in\hat{B}_i\). Hence
\(B_{1}=\bigcup_{(a',b')\in\hat{B}_{1}}Q(a',b')\)
and \(B_2=\bigcup_{(a',b')\in\hat{B}_2}Q(a',b')\), and thus
\(B_1\cap B_2=\bigcup_{(a',b')\in\hat{B}_1\cap\hat{B}_2}Q(a',b')
\in\mathscr{F}\).

(iii) \(\mathcal{A}^{\Z\setminus[0,r)}\setminus B\in\mathscr{F}\)
for any \(B\in\mathscr{F}\).
Choose \(m,n\in\N_0\) and \(F\subseteq\mathcal{A}^m\times\mathcal{A}^n\)
such that \(B=\bigcup_{(a',b')\in F}Q(a',b')\).
For any \(z\in\mathcal{A}^{\Z\setminus[0,r)}\), the definition of \(Q\) gives
\(z\in B\) if and only if \((z_{[-m,0)},z_{[r,r+n)})\in F\), and thus
\(\mathcal{A}^{\Z\setminus[0,r)}\setminus B
=\bigcup_{(a',b')\in(\mathcal{A}^m\times\mathcal{A}^n)\setminus F}
Q(a',b')\in\mathscr{F}\).

\smallskip

Thus, \(\mathscr{F}\) is an algebra of subsets of \(\mathcal{A}^{\Z\setminus[0,r)}\).

\smallskip

We now show that \(\mathscr{A}(\mathscr{F})\) is the Borel
\(\sigma\)-algebra of \(\mathcal{A}^{\Z\setminus[0,r)}\).
Every member of \(\mathscr{F}\) is a finite union of clopen cylinders
in \mbox{\(\mathcal{A}^{\Z\setminus[0,r)}\)}.
The Borel \(\sigma\)-algebra therefore contains \(\mathscr{F}\) and hence
contains \(\mathscr{A}(\mathscr{F})\).
Conversely, every cylinder specified by finitely many coordinates in
\(\Z\setminus[0,r)\) belongs to \(\mathscr{F}\): choosing \(m,n\in\N_0\)
so that all prescribed coordinates lie in \([-m,0)\cup[r,r+n)\) expresses
it as a finite union of cylinders \(Q(a',b')\) with
\(a'\in\mathcal{A}^m\) and \(b'\in\mathcal{A}^n\).
Since the finite-coordinate cylinders form a countable basis for the product
topology, every open set is a countable union of members of \(\mathscr{F}\)
and hence belongs to \(\mathscr{A}(\mathscr{F})\).
As \(\mathscr{A}(\mathscr{F})\) is a \(\sigma\)-algebra,
it contains all Borel sets.

Define
\(\mathscr{C}=\{B\in\mathscr{A}(\mathscr{F}):
\mu_u(\pi_r([u]_\Sigma)\cap\pi_r([v]_\Sigma)\cap B)\leq\mu_v(B)\}\).
We first show that \(\mathscr{F}\subseteq\mathscr{C}\).
For any \(B\in\mathscr{F}\), choose \(m,n\in\N_0\) and \(F\subseteq\mathcal{A}^m\times\mathcal{A}^n\)
such that \(B=\bigcup_{(a',b')\in F}Q(a',b')\).
Since the cylinders \(Q(a',b')\), \((a',b')\in F\), are pairwise disjoint,
the finite additivity of \(\mu_u\) and \(\mu_v\), together with the estimates
proved in Cases~\ref{lem:finite-replacement-invariance}.3 and~\ref{lem:finite-replacement-invariance}.4, gives
\(\mu_u(\pi_r([u]_\Sigma)\cap\pi_r([v]_\Sigma)\cap B)
=\sum_{(a',b')\in F}
\mu_u(\pi_r([u]_\Sigma)\cap\pi_r([v]_\Sigma)\cap Q(a',b'))
\leq\sum_{(a',b')\in F}\mu_v(Q(a',b'))=\mu_v(B)\),
implying that \(B\in\mathscr{C}\).

We next show that \(\mathscr{C}\) is a monotone class.
For any increasing or decreasing sequence \((B_k)_{k\in\N}\) in
\(\mathscr{C}\), put \(B=\bigcup_{k\in\N}B_k\) or
\(B=\bigcap_{k\in\N}B_k\), respectively.
In either case, \(B\in\mathscr{A}(\mathscr{F})\) by
\cite[Theorem~0.6]{Walters1982}.
Since \(\mu_u\) and \(\mu_v\) are finite measures, continuity from below
or from above, respectively, gives
\(\mu_u(\pi_r([u]_\Sigma)\cap\pi_r([v]_\Sigma)\cap B)
=\lim_{k\to\infty}
\mu_u(\pi_r([u]_\Sigma)\cap\pi_r([v]_\Sigma)\cap B_k)
\leq\lim_{k\to\infty}\mu_v(B_k)=\mu_v(B)\).
Hence \(B\in\mathscr{C}\), so \(\mathscr{C}\) is a monotone class
containing \(\mathscr{F}\).

By \cite[Theorem~0.6]{Walters1982}, \(\mathscr{A}(\mathscr{F})\) is the
smallest monotone class containing \(\mathscr{F}\), and hence
\(\mathscr{A}(\mathscr{F})\subseteq\mathscr{C}\).
The reverse inclusion holds by the definition of \(\mathscr{C}\).
Thus \(\mathscr{C}=\mathscr{A}(\mathscr{F})\), which proves
\(\mu_u(\pi_r([u]_\Sigma)\cap\pi_r([v]_\Sigma)\cap B)\leq\mu_v(B)\)
for any Borel set \(B\subseteq\mathcal{A}^{\Z\setminus[0,r)}\).

Similarly, with \(u\) and \(v\) interchanged, it can be verified that
\(\mu_v(\pi_r([v]_\Sigma)\cap\pi_r([u]_\Sigma)\cap B)\leq\mu_u(B)\)
for any Borel set \(B\subseteq\mathcal{A}^{\Z\setminus[0,r)}\).
Thus
\(\mu_u|_{\pi_r([u]_\Sigma)\cap\pi_r([v]_\Sigma)}
=\mu_v|_{\pi_r([u]_\Sigma)\cap\pi_r([v]_\Sigma)}\).

For any Borel set \(E\subseteq D_{\Sigma}(u,v)\), \(\pi_r(E)\) is a Borel subset of
\(\mathcal{A}^{\Z\setminus[0,r)}\), since \(\pi_r|_{[u]_\Sigma}\) is a
homeomorphism onto the compact set \(\pi_r([u]_\Sigma)\).
Since \(E\subseteq D_{\Sigma}(u,v)\), the definition of \(D_{\Sigma}(u,v)\) and \eqref{eq:replacement-map}
give \(E\subseteq[u]_\Sigma\) and \(\Phi^{\Sigma}_{u,v}(E)\subseteq[v]_\Sigma\).
Using \(\pi_r\circ(\Phi^{\Sigma}_{u,v}|_{D_{\Sigma}(u,v)})=\pi_r|_{D_{\Sigma}(u,v)}\) (by \eqref{eq:replacement-map} and the definition of \(\pi_r\)) yields
\begin{equation}\label{eq:replacement-projection-images}
 \pi_r(E)\subseteq\pi_r([u]_\Sigma)\text{ and }
 \pi_r(E)=\pi_r(\Phi^{\Sigma}_{u,v}(E))\subseteq\pi_r([v]_\Sigma),
\end{equation}
implying that \(\pi_r(E)\subseteq\pi_r([u]_\Sigma)\cap\pi_r([v]_\Sigma)\),
and therefore
\begin{equation}\label{eq:replacement-projected-measure-equality}
 \mu_u(\pi_r(E))=\mu_v(\pi_r(E))
\end{equation}
by
\(\mu_u|_{\pi_r([u]_\Sigma)\cap\pi_r([v]_\Sigma)}
=\mu_v|_{\pi_r([u]_\Sigma)\cap\pi_r([v]_\Sigma)}\).
By \eqref{eq:replacement-projection-images} and the injectivity of
\(\pi_r|_{[u]_\Sigma}\) and \(\pi_r|_{[v]_\Sigma}\), we have
\([u]_\Sigma\cap\pi_r^{-1}(\pi_r(E))=E\) and
\([v]_\Sigma\cap\pi_r^{-1}(\pi_r(E))
=[v]_\Sigma\cap\pi_r^{-1}(\pi_r(\Phi^{\Sigma}_{u,v}(E)))=\Phi^{\Sigma}_{u,v}(E).\)
These, together with the definitions of \(\mu_u\) and \(\mu_v\)
and \eqref{eq:replacement-projected-measure-equality}, imply
\(\mu(E)=\mu([u]_\Sigma\cap\pi_r^{-1}(\pi_r(E)))
=\mu_u(\pi_r(E))=\mu_v(\pi_r(E))
=\mu([v]_\Sigma\cap\pi_r^{-1}(\pi_r(E)))
=\mu(\Phi^{\Sigma}_{u,v}(E)).\)

\smallskip

Finally, we prove that \(\mu(B)=\mu(\Phi^{\Sigma}_{u,v}(B))\)
for any invariant probability measure \(\mu\) on \(\Sigma\)
with \(h_\mu(\sigma|_\Sigma)=\htop(\Sigma)\) and any Borel set \(B\subseteq D_{\Sigma}(u,v)\).

\smallskip

Write the ergodic decomposition of \(\mu\) as
\(\mu=\int_{\mathscr{E}(\Sigma,\sigma)}\rho\,d\lambda(\rho)\)
(\cite[p.~153, Remark~(2)]{Walters1982}),
where \(\lambda\) is a probability measure on the Borel subsets of the space
of invariant probability measures on \(\Sigma\), endowed with the weak* topology,
and \(\lambda(\mathscr{E}(\Sigma,\sigma))=1\).
Jacobs' theorem (\cite[Theorem~8.4~(ii)]{Walters1982}) gives
\[
\int_{\mathscr{E}(\Sigma,\sigma)}
(\htop(\Sigma)-h_\rho(\sigma|_\Sigma))\,d\lambda(\rho)
=\htop(\Sigma)-h_\mu(\sigma|_\Sigma)=0.
\]
This, together with
\(h_\rho(\sigma|_\Sigma)\leq\htop(\Sigma)\) for any
\(\rho\in\mathscr{E}(\Sigma,\sigma)\), implies
\(\lambda(\{\rho\in\mathscr{E}(\Sigma,\sigma):
h_\rho(\sigma|_\Sigma)=\htop(\Sigma)\})=1\).
For any \(\rho\) in this set, the ergodic case gives
\(\rho(B)=\rho(\Phi^{\Sigma}_{u,v}(B))\)
for any Borel set \(B\subseteq D_{\Sigma}(u,v)\),
implying that
\[
\mu(B)=\int_{\mathscr{E}(\Sigma,\sigma)}\rho(B)\,d\lambda(\rho)
=\int_{\mathscr{E}(\Sigma,\sigma)}\rho(\Phi^{\Sigma}_{u,v}(B))\,d\lambda(\rho)
=\mu(\Phi^{\Sigma}_{u,v}(B)).\qedhere
\]
\end{proof}

We now combine Lemma~\ref{lem:uniform-coverage} with the
measure-preservation property in Lemma~\ref{lem:finite-replacement-invariance},
applied to the product subshift \(\mathcal{X}_\Sigma\), to obtain a uniform
separation estimate for distinct ergodic measures of maximal entropy.

\begin{lemma}
\label{lem:positive-entropy-separation}
Let \(\Sigma\) be a subshift over a finite alphabet \(\mathcal{A}\) satisfying
\(\htop(\Sigma)>0\) and \(\hG(\Sigma)<\htop(\Sigma).\)
Then there exists \(\ell\in\N\) such that
\begin{equation}
\label{eq:finite-word-separation}
 \frac{1}{2}\sum_{w\in\cL(\Sigma,2\ell)}
               |\mu([w]_{\Sigma})-\nu([w]_{\Sigma})|\leq\frac{1}{16}
 \ \Longrightarrow\ \mu=\nu
\end{equation}
for any pair of ergodic measures of maximal entropy \(\mu,\nu\) on \(\Sigma\).
\end{lemma}

\begin{proof}
Put
\(s^*=\frac{1}{2}(\htop(\Sigma)+\max\{0,\hG(\Sigma)\})\) and
\(t^*=\frac{1}{4}(3\htop(\Sigma)+s^*)\). Then
\(\max\{0,\hG(\Sigma)\}<s^*<\frac{1}{2}(\htop(\Sigma)+s^*)
<t^*<\htop(\Sigma)\) by \(\htop(\Sigma)>0\) and \(\hG(\Sigma)<\htop(\Sigma)\).
\begin{claim}\label{clm:uniform-cylinder-estimates}
There exists a constant \(C_0\geq1\) such that, for any ergodic measure of maximal entropy
\(\mu\) on \(\Sigma\), there exists a closed set \(G_\mu\subseteq\Sigma\) with
\(\mu(G_\mu)\geq\frac{15}{16}\) satisfying
the following two properties:
\begin{enumerate}
\item[(i)] For any \(i,j\in\Z\) with \(i\leq0<j\)
and any \(w\in\mathcal{A}^{j-i}\),
\begin{equation}
\label{eq:good-cylinder}
 \mu(G_{\mu}\cap\{x\in\Sigma:x_{[i,j)}=w\})\leq C_{0}
 e^{-t^*(j-i)}.
\end{equation}
\item[(ii)] For any \(n\in\N\) and any \(w\in\mathcal{A}^n\),
\begin{equation}
\label{eq:half-cylinder}
 \begin{cases}
 \mu(G_{\mu}\cap\{x\in\Sigma:x_{[-n,0)}=w\})\leq C_{0} e^{-t^*n},\\
 \mu(G_{\mu}\cap\{x\in\Sigma:x_{[0,n)}=w\})\leq C_{0} e^{-t^*n}.
 \end{cases}
\end{equation}
\end{enumerate}
\end{claim}

\begin{proof}[Proof of Claim~\ref{clm:uniform-cylinder-estimates}]
Choose \(K_0\in\N\) as in Lemma~\ref{lem:uniform-coverage}
with \(t=t^*\) and \(\delta=\frac{1}{16}\), and put
\(C_0=(\#\mathcal{A})e^{t^*K_0}\geq2>1\).
Fix any ergodic measure of maximal entropy \(\mu\). For any \(n\in\N\), set
\(\mathcal{B}_n^{(\mu)}
=\{w\in\cL(\Sigma,n):\mu([w]_\Sigma)>e^{-t^*n}\}\).
Since the cylinders \([w]_\Sigma\), \(w\in\mathcal{B}_n^{(\mu)}\),
are pairwise disjoint, the definition of \(\mathcal{B}_n^{(\mu)}\) gives
\(e^{-t^*n}\#\mathcal{B}_n^{(\mu)}
\leq\sum_{w\in\mathcal{B}_n^{(\mu)}}\mu([w]_\Sigma)\leq1\), implying that
\(\#\mathcal{B}_n^{(\mu)}\leq e^{t^*n}\). Let
\(G_\mu=\Sigma\setminus U
((\mathcal{B}_n^{(\mu)})_{n\geq K_0})\).
Then \(G_\mu\) is closed and, by the choice of \(K_0\)
and Lemma~\ref{lem:uniform-coverage},
satisfies \(\mu(G_\mu)\geq\frac{15}{16}\).
For any \(i,j\in \Z\) with \(i\leq0<j\) and any word
\(w\in\mathcal{A}^{j-i}\), we consider the following two cases:

\smallskip
\textbf{Case \ref{clm:uniform-cylinder-estimates}.1.} If \(j-i\geq K_0\),
the \(\sigma\)-invariance of \(\mu\) and the definitions of \(G_{\mu}\),
\(U((\mathcal{B}_n^{(\mu)})_{n\geq K_0})\), and
\(\mathcal{B}_{j-i}^{(\mu)}\) give
\[
 \begin{cases}
 G_{\mu}\cap\{x\in\Sigma:x_{[i,j)}=w\}\subseteq G_{\mu}\cap U
((\mathcal{B}_n^{(\mu)})_{n\geq K_0})=\varnothing,
                                      &w\in\mathcal{B}_{j-i}^{(\mu)},\\
 \mu(G_{\mu}\cap\{x\in\Sigma:x_{[i,j)}=w\})\leq
 \mu(\{x\in\Sigma:x_{[i,j)}=w\})=
 \mu([w]_{\Sigma})\leq e^{-t^*(j-i)},
                                      &w\notin\mathcal{B}_{j-i}^{(\mu)}.
 \end{cases}
\]

\smallskip
\textbf{Case \ref{clm:uniform-cylinder-estimates}.2.} If \(j-i<K_0\), then
\(\mu(G_{\mu}\cap\{x\in\Sigma:x_{[i,j)}=w\})\leq 1\leq e^{t^*K_0-t^*(j-i)}\).

In both cases,
\(\mu(G_\mu\cap\{x\in\Sigma:x_{[i,j)}=w\})
\leq e^{t^*K_0-t^*(j-i)}\leq C_0 e^{-t^*(j-i)}\),
which proves \eqref{eq:good-cylinder}.
For any \(n\in\N\) and any \(w\in\mathcal{A}^n\), taking
\(i=0\) and \(j=n\) in \eqref{eq:good-cylinder} yields the second inequality
in \eqref{eq:half-cylinder}.
For the first inequality in \eqref{eq:half-cylinder}, taking
\(i=-n\) and \(j=1\) in the common estimate from Cases~\ref{clm:uniform-cylinder-estimates}.1
and~\ref{clm:uniform-cylinder-estimates}.2
and partitioning the cylinder \(\{x\in\Sigma:x_{[-n,0)}=w\}\)
according to the value of \(x_0\), we obtain
\begin{align*}
 \mu(G_{\mu}\cap\{x\in\Sigma:x_{[-n,0)}=w\})
 &=\sum_{\alpha\in\mathcal{A}}
                 \mu(G_{\mu}\cap\{x\in\Sigma:x_{[-n,1)}=w\alpha\})\\
 &\leq(\#\mathcal{A}) e^{t^*K_0
 -t^*(n+1)}\leq C_0 e^{-t^*n}.\qedhere
\end{align*}
\end{proof}

Fix \(C_0\geq1\) as in Claim~\ref{clm:uniform-cylinder-estimates}.
For each ergodic measure of maximal entropy \(\mu\) on \(\Sigma\), choose
a closed set \(G_\mu\subseteq\Sigma\) satisfying
the conclusions of Claim~\ref{clm:uniform-cylinder-estimates}
with this fixed \(C_0\).
Since \(\hG(\Sigma)<s^*\), there exists a set of forbidden words
\(\mathcal{F}\) generating \(\Sigma\) such that
\(\limsup_{n\to\infty}\frac{1}{n}\log\#\mathcal{F}_n<s^*\),
implying that there exists a constant \(C_{\mathcal{F}}\geq1\) such that
\(\#\mathcal{F}_n\leq C_{\mathcal{F}}\cdot e^{s^*n}\) for all \(n\in\N\).

Define \(\Gamma:\Sigma^2\to\mathcal{A}^{\Z}\) by
\begin{equation}\label{eq:splice-map}
 \Gamma(x,y)_i=
 \begin{cases}
 x_i,&i<0,\\
 y_i,&i\geq0.
 \end{cases}
\end{equation}
It can be verified that the map \(\Gamma\) is continuous.
For any pair of ergodic measures of maximal entropy \(\mu,\nu\) on \(\Sigma\), put
\(E_{\mu,\nu}=\{(x,y)\in G_{\mu}\times G_{\nu}:
\Gamma(x,y)\in\Sigma\text{ and }\Gamma(y,x)\in\Sigma\}\).
Since \(G_\mu,G_\nu,\Sigma\) are closed,
the continuity of \(\Gamma\) implies that
the set \(E_{\mu,\nu}\) is closed, hence Borel.

For any \(\ell\in\N\), put
\(\xi_{\ell}=e^{\frac{64}{13}\left[2(\log\#\cL(\Sigma,\ell)-\ell\htop(\Sigma))+\frac{1}{e}\right]}\)
and \(\kappa_{\ell}=\left(\frac{13}{64}\right)^2\left(\frac{13}{32\xi_{\ell}}\right)^4\). By \eqref{eq:language-limit}, we have
\begin{equation}
\lim_{\ell\to\infty}\frac{\log \xi_{\ell}}{\ell}=\frac{64}{13}\lim_{\ell\to\infty}
 \left[2\left(\frac{\log\#\cL(\Sigma,\ell)}{\ell}-\htop(\Sigma)\right)
       +\frac{1}{e\ell}\right]=0,
\end{equation}
and
\begin{equation}\label{eq:splice-normalization-limits}
 \lim_{\ell\to\infty}\frac{\log((\ell+1)/\kappa_{\ell})}{\ell}=\lim_{\ell\to\infty}\left[
 \frac{\log(\ell+1)+4\log \xi_{\ell}}{\ell}
 +\frac{2\log(64/13)+4\log(32/13)}{\ell}\right]=0.
\end{equation}

\begin{claim}\label{clm:two-splices}
There exists \(\ell^{*}\in\N\) such that
\((\mu\times\nu)(E_{\mu,\nu})\geq\frac{\kappa_{\ell^{*}}}{2}>0\) for any pair of ergodic measures of maximal entropy
\(\mu,\nu\) on \(\Sigma\) with \(\TV(\mathbb{P}^{(1)}_{\mu,\ell^{*}},\mathbb{P}^{(1)}_{\nu,\ell^{*}})\leq\frac{1}{16}\).
\end{claim}
\begin{proof}[Proof of Claim~\ref{clm:two-splices}]
Fix any \(\ell\in\N\) and any pair of ergodic measures of maximal entropy \(\mu,\nu\) on \(\Sigma\) with
\(
 \TV(\mathbb{P}^{(1)}_{\mu,\ell},\mathbb{P}^{(1)}_{\nu,\ell})\leq\frac{1}{16}.
\)
Define
\(\mathbf{A}_\mu,\mathbf{A}_\nu,\mathbf{R}_{\ell}:\cL(\Sigma,\ell)^2\to[0,1]\) by
\[
 \begin{cases}
 \mathbf{A}_{\mu}(u,v)=\mu(G_{\mu}\cap
             \{x\in\Sigma:x_{[-\ell,\ell)}=uv\}),\\
 \mathbf{A}_{\nu}(u,v)=\nu(G_{\nu}\cap
             \{y\in\Sigma:y_{[-\ell,\ell)}=uv\}),\\
 \mathbf{R}_{\ell}(u,v)=\min\{\mathbf{A}_{\mu}(u,v),\mathbf{A}_{\nu}(u,v)\}.
 \end{cases}
\]
For any \(u,v\in\cL(\Sigma,\ell)\), the \(\sigma\)-invariance of
\(\mu,\nu\) gives
\(\mathbf{R}_{\ell}(u,v)\leq \mathbf{A}_\mu(u,v)
\leq\mu(\{x\in\Sigma:x_{[-\ell,\ell)}=uv\})
=\mu(\sigma^{\ell}([uv]_\Sigma))=\mu([uv]_\Sigma)=\mathbb{P}^{(1)}_{\mu,\ell}(u,v)\)
and
\(\mathbf{R}_{\ell}(u,v)\leq \mathbf{A}_\nu(u,v)
\leq\nu(\{y\in\Sigma:y_{[-\ell,\ell)}=uv\})
=\nu(\sigma^{\ell}([uv]_\Sigma))=\nu([uv]_\Sigma)=\mathbb{P}^{(1)}_{\nu,\ell}(u,v)\),
implying that
\begin{equation}\label{eq:common-block-pointwise}
 \begin{aligned}
 \mathbf{R}_{\ell}(u,v)&=\min\{\mathbf{A}_\mu(u,v),\mathbf{A}_\nu(u,v)\}\\
 &\geq\min\{\mathbf{A}_\mu(u,v)-(\mathbb{P}^{(1)}_{\nu,\ell}(u,v)-\mathbf{A}_\nu(u,v)),\
                 \mathbf{A}_\nu(u,v)-(\mathbb{P}^{(1)}_{\mu,\ell}(u,v)-\mathbf{A}_\mu(u,v))\}\\
 &=\min\{\mathbb{P}^{(1)}_{\mu,\ell}(u,v),\mathbb{P}^{(1)}_{\nu,\ell}(u,v)\}
 -(\mathbb{P}^{(1)}_{\mu,\ell}(u,v)-\mathbf{A}_\mu(u,v))
        -(\mathbb{P}^{(1)}_{\nu,\ell}(u,v)-\mathbf{A}_\nu(u,v)).
 \end{aligned}
\end{equation}
The cylinder sets \(\{x\in\Sigma:x_{[-\ell,\ell)}=uv\}\),
\(u,v\in\cL(\Sigma,\ell)\), form a finite measurable partition of \(\Sigma\).
Consequently,
\begin{equation}\label{eq:restricted-block-masses}
 \sum_{u,v\in\cL(\Sigma,\ell)}\mathbf{A}_\mu(u,v)=\mu(G_\mu)
 \quad \text{and}\quad
 \sum_{u,v\in\cL(\Sigma,\ell)}\mathbf{A}_\nu(u,v)=\nu(G_\nu).
\end{equation}
Since the cylinder sets \([w]_\Sigma\), \(w\in\cL(\Sigma,2\ell)\),
form a finite measurable partition of \(\Sigma\), we have
\[
\sum_{u,v\in\cL(\Sigma,\ell)}\mathbb{P}^{(1)}_{\mu,\ell}(u,v)
\allowbreak=\nobreak\sum_{u,v\in\cL(\Sigma,\ell)}\mu([uv]_\Sigma)
\allowbreak=\nobreak\sum_{w\in\cL(\Sigma,2\ell)}\mu([w]_\Sigma)
\allowbreak=\nobreak\mu(\Sigma)\allowbreak=\nobreak 1,
\]
and
\[
\sum_{u,v\in\cL(\Sigma,\ell)}\mathbb{P}^{(1)}_{\nu,\ell}(u,v)
\allowbreak=\nobreak\sum_{u,v\in\cL(\Sigma,\ell)}\nu([uv]_\Sigma)
\allowbreak=\nobreak\sum_{w\in\cL(\Sigma,2\ell)}\nu([w]_\Sigma)
\allowbreak=\nobreak\nu(\Sigma)\allowbreak=\nobreak1.
\]
These, together with \eqref{eq:common-block-pointwise},
\eqref{eq:restricted-block-masses}, \eqref{eq:total-variation},
\(\TV(\mathbb{P}^{(1)}_{\mu,\ell},\mathbb{P}^{(1)}_{\nu,\ell})\leq\frac{1}{16}\), and
Claim~\ref{clm:uniform-cylinder-estimates}, yield
\begin{equation}\label{eq:common-block-mass}
 \begin{aligned}
 \sum_{u,v\in\cL(\Sigma,\ell)}\mathbf{R}_{\ell}(u,v)
 &\geq\sum_{u,v\in\cL(\Sigma,\ell)}\min\{\mathbb{P}^{(1)}_{\mu,\ell}(u,v),\mathbb{P}^{(1)}_{\nu,\ell}(u,v)\}\\
 &\quad-\sum_{u,v\in\cL(\Sigma,\ell)}(\mathbb{P}^{(1)}_{\mu,\ell}(u,v)-\mathbf{A}_\mu(u,v))\\
 &\quad-\sum_{u,v\in\cL(\Sigma,\ell)}(\mathbb{P}^{(1)}_{\nu,\ell}(u,v)-\mathbf{A}_\nu(u,v))\\
 &=1-\TV(\mathbb{P}^{(1)}_{\mu,\ell},\mathbb{P}^{(1)}_{\nu,\ell})
        -(1-\mu(G_\mu))-(1-\nu(G_\nu))\\
 &\geq1-\frac{1}{16}-\frac{1}{16}-\frac{1}{16}=\frac{13}{16}.
 \end{aligned}
\end{equation}

Let
\(\mathscr{T}_{\ell}=\{(u,v)\in\cL(\Sigma,\ell)^2:\mathbb{P}^{(2)}_{\mu,\ell}(u,v)>0,\
\frac{13}{64}\mathbb{P}^{(2)}_{\mu,\ell}(u,v)\leq \mathbf{R}_{\ell}(u,v)\leq \xi_{\ell}\mathbb{P}^{(2)}_{\mu,\ell}(u,v)\}\).
We first verify that \(\mathbb{P}^{(2)}_{\mu,\ell}(\mathscr{T}_{\ell})\geq\frac{13}{32\xi_{\ell}}\).
By the definition of \(\mathbb{P}^{(2)}_{\mu,\ell}\),
\[
\sum_{u,v\in\cL(\Sigma,\ell)}\mathbb{P}^{(2)}_{\mu,\ell}(u,v)
\allowbreak=\nobreak\sum_{u,v\in\cL(\Sigma,\ell)}\mu([u]_\Sigma)\mu([v]_\Sigma)
\allowbreak=\nobreak\left(\sum_{u\in\cL(\Sigma,\ell)}\mu([u]_\Sigma)\right)^2
\allowbreak=\nobreak1,
\]
implying that
\begin{equation}\label{eq:small-ratio-mass}
 \sum_{\substack{u,v\in\cL(\Sigma,\ell)\\\mathbf{R}_{\ell}(u,v)<\frac{13}{64}\mathbb{P}^{(2)}_{\mu,\ell}(u,v)}}\mathbf{R}_{\ell}(u,v)
 \leq\sum_{\substack{u,v\in\cL(\Sigma,\ell)\\\mathbf{R}_{\ell}(u,v)<\frac{13}{64}\mathbb{P}^{(2)}_{\mu,\ell}(u,v)}}
 \frac{13}{64}\mathbb{P}^{(2)}_{\mu,\ell}(u,v)
 \leq\frac{13}{64}.
\end{equation}
For any \(u\in\cL(\Sigma,\ell)\), the
\(\sigma\)-invariance of \(\mu\) gives
\[
 \begin{cases}
 \displaystyle\sum_{v\in\cL(\Sigma,\ell)}\mathbb{P}^{(1)}_{\mu,\ell}(u,v)
  =\sum_{v\in\cL(\Sigma,\ell)}
      \mu([uv]_\Sigma)=\mu([u]_\Sigma),\\
 \displaystyle\sum_{v\in\cL(\Sigma,\ell)}\mathbb{P}^{(1)}_{\mu,\ell}(v,u)
  =\sum_{v\in\cL(\Sigma,\ell)}
      \mu([vu]_\Sigma)
      =\mu(\sigma^{-\ell}([u]_\Sigma))=\mu([u]_\Sigma).
 \end{cases}
\]
By \eqref{eq:joint-block-pointwise},
\(\mathbb{P}^{(1)}_{\mu,\ell}(u,v)>0\) implies \(\mathbb{P}^{(2)}_{\mu,\ell}(u,v)>0\)
for all \(u,v\in\cL(\Sigma,\ell)\).
These, together with the definition of \(\mathbb{P}^{(1)}_{\mu,\ell}\) and
\eqref{eq:block-entropy}, imply
\begin{equation}
\label{eq:mutual-information}
 \begin{split}
 &\sum_{\substack{u,v\in\cL(\Sigma,\ell)\\\mathbb{P}^{(1)}_{\mu,\ell}(u,v)>0}}
       \mathbb{P}^{(1)}_{\mu,\ell}(u,v)\log\frac{\mathbb{P}^{(1)}_{\mu,\ell}(u,v)}{\mathbb{P}^{(2)}_{\mu,\ell}(u,v)}
 =\sum_{\substack{u,v\in\cL(\Sigma,\ell)\\\mathbb{P}^{(1)}_{\mu,\ell}(u,v)>0}}
       \mathbb{P}^{(1)}_{\mu,\ell}(u,v)\log\frac{\mathbb{P}^{(1)}_{\mu,\ell}(u,v)}
       {\mu([u]_\Sigma)\mu([v]_\Sigma)}\\
 &=\sum_{\substack{u,v\in\cL(\Sigma,\ell)\\\mathbb{P}^{(1)}_{\mu,\ell}(u,v)>0}}
       \mathbb{P}^{(1)}_{\mu,\ell}(u,v)(\log \mathbb{P}^{(1)}_{\mu,\ell}(u,v)
       -\log\mu([u]_\Sigma)-\log\mu([v]_\Sigma))\\
 &=\sum_{u,v\in\cL(\Sigma,\ell)}\mathbb{P}^{(1)}_{\mu,\ell}(u,v)\log \mathbb{P}^{(1)}_{\mu,\ell}(u,v)
 -\sum_{\substack{u\in\cL(\Sigma,\ell)\\\mu([u]_\Sigma)>0}}
      \left(\sum_{v\in\cL(\Sigma,\ell)}\mathbb{P}^{(1)}_{\mu,\ell}(u,v)\right)\log\mu([u]_\Sigma)\\
 &\quad-\sum_{\substack{v\in\cL(\Sigma,\ell)\\\mu([v]_\Sigma)>0}}
      \left(\sum_{u\in\cL(\Sigma,\ell)}\mathbb{P}^{(1)}_{\mu,\ell}(u,v)\right)\log\mu([v]_\Sigma)\\
 &=\sum_{u,v\in\cL(\Sigma,\ell)}\mu([uv]_\Sigma)\log\mu([uv]_\Sigma)
   -2\sum_{u\in\cL(\Sigma,\ell)}\mu([u]_\Sigma)\log\mu([u]_\Sigma)\\
 &=\sum_{w\in\cL(\Sigma,2\ell)}\mu([w]_\Sigma)\log\mu([w]_\Sigma)
   -2\sum_{u\in\cL(\Sigma,\ell)}\mu([u]_\Sigma)\log\mu([u]_\Sigma)\\
 &=-H_\mu(2\ell)+2H_\mu(\ell).
 \end{split}
\end{equation}
Moreover, using \(z\log\frac{1}{z}\leq \frac{1}{e}\) for \(0<z\leq1\) yields
\begin{equation}\label{eq:negative-information-bound}
 \begin{aligned}
 &\sum_{\substack{u,v\in\cL(\Sigma,\ell)\\0<\mathbb{P}^{(1)}_{\mu,\ell}(u,v)\leq\mathbb{P}^{(2)}_{\mu,\ell}(u,v)}}
   \mathbb{P}^{(1)}_{\mu,\ell}(u,v)\log\frac{\mathbb{P}^{(2)}_{\mu,\ell}(u,v)}{\mathbb{P}^{(1)}_{\mu,\ell}(u,v)}\\
   &\quad=\sum_{\substack{u,v\in\cL(\Sigma,\ell)\\0<\mathbb{P}^{(1)}_{\mu,\ell}(u,v)\leq\mathbb{P}^{(2)}_{\mu,\ell}(u,v)}}
   \mathbb{P}^{(2)}_{\mu,\ell}(u,v)\frac{\mathbb{P}^{(1)}_{\mu,\ell}(u,v)}{\mathbb{P}^{(2)}_{\mu,\ell}(u,v)}
   \log\frac{\mathbb{P}^{(2)}_{\mu,\ell}(u,v)}{\mathbb{P}^{(1)}_{\mu,\ell}(u,v)}\\
 &\quad\leq\frac{1}{e}
   \sum_{\substack{u,v\in\cL(\Sigma,\ell)\\0<\mathbb{P}^{(1)}_{\mu,\ell}(u,v)\leq\mathbb{P}^{(2)}_{\mu,\ell}(u,v)}}\mathbb{P}^{(2)}_{\mu,\ell}(u,v)
   \leq\frac{1}{e}.
 \end{aligned}
\end{equation}
By \eqref{eq:language-limit} and the definition of \(\xi_{\ell}\), we have
\(\log\#\cL(\Sigma,\ell)\geq \ell\htop(\Sigma)\) and
\(\xi_{\ell}\geq e^{64/(13e)}>1\).
Since \(\mu\) is a measure of maximal entropy, by \eqref{eq:mutual-information},
\eqref{eq:block-entropy-bounds}, \eqref{eq:negative-information-bound}, and
\(\log \xi_{\ell}=\frac{64}{13}[2(\log\#\cL(\Sigma,\ell)-\ell\htop(\Sigma))+\frac{1}{e}]\),
we obtain
\[
 \begin{aligned}
 &\sum_{\substack{u,v\in\cL(\Sigma,\ell)\\\mathbb{P}^{(1)}_{\mu,\ell}(u,v)>\xi_{\ell}\mathbb{P}^{(2)}_{\mu,\ell}(u,v)}}
   \mathbb{P}^{(1)}_{\mu,\ell}(u,v)\log \xi_{\ell}
   \leq\sum_{\substack{u,v\in\cL(\Sigma,\ell)\\\mathbb{P}^{(1)}_{\mu,\ell}(u,v)>\xi_{\ell}\mathbb{P}^{(2)}_{\mu,\ell}(u,v)}}
   \mathbb{P}^{(1)}_{\mu,\ell}(u,v)\log\frac{\mathbb{P}^{(1)}_{\mu,\ell}(u,v)}{\mathbb{P}^{(2)}_{\mu,\ell}(u,v)}\\
 &\leq\sum_{\substack{u,v\in\cL(\Sigma,\ell)\\\mathbb{P}^{(1)}_{\mu,\ell}(u,v)>\mathbb{P}^{(2)}_{\mu,\ell}(u,v)}}
    \mathbb{P}^{(1)}_{\mu,\ell}(u,v)\log\frac{\mathbb{P}^{(1)}_{\mu,\ell}(u,v)}{\mathbb{P}^{(2)}_{\mu,\ell}(u,v)}\\
 &=\sum_{\substack{u,v\in\cL(\Sigma,\ell)\\\mathbb{P}^{(1)}_{\mu,\ell}(u,v)>0}}
    \mathbb{P}^{(1)}_{\mu,\ell}(u,v)\log\frac{\mathbb{P}^{(1)}_{\mu,\ell}(u,v)}{\mathbb{P}^{(2)}_{\mu,\ell}(u,v)}
    -\sum_{\substack{u,v\in\cL(\Sigma,\ell)\\0<\mathbb{P}^{(1)}_{\mu,\ell}(u,v)\leq\mathbb{P}^{(2)}_{\mu,\ell}(u,v)}}
    \mathbb{P}^{(1)}_{\mu,\ell}(u,v)\log\frac{\mathbb{P}^{(1)}_{\mu,\ell}(u,v)}{\mathbb{P}^{(2)}_{\mu,\ell}(u,v)}\\
 &=\sum_{\substack{u,v\in\cL(\Sigma,\ell)\\\mathbb{P}^{(1)}_{\mu,\ell}(u,v)>0}}
    \mathbb{P}^{(1)}_{\mu,\ell}(u,v)\log\frac{\mathbb{P}^{(1)}_{\mu,\ell}(u,v)}{\mathbb{P}^{(2)}_{\mu,\ell}(u,v)}
    +\sum_{\substack{u,v\in\cL(\Sigma,\ell)\\0<\mathbb{P}^{(1)}_{\mu,\ell}(u,v)\leq\mathbb{P}^{(2)}_{\mu,\ell}(u,v)}}
    \mathbb{P}^{(1)}_{\mu,\ell}(u,v)\log\frac{\mathbb{P}^{(2)}_{\mu,\ell}(u,v)}{\mathbb{P}^{(1)}_{\mu,\ell}(u,v)}\\
 &\leq-H_\mu(2\ell)+2H_\mu(\ell)+\frac{1}{e}
 \leq2(\log\#\cL(\Sigma,\ell)-\ell\htop(\Sigma))+\frac{1}{e}
   =\frac{13}{64}\log \xi_{\ell},
 \end{aligned}
\]
and hence
\[
\sum_{\substack{u,v\in\cL(\Sigma,\ell)\\\mathbb{P}^{(1)}_{\mu,\ell}(u,v)>\xi_{\ell}\mathbb{P}^{(2)}_{\mu,\ell}(u,v)}}
\mathbb{P}^{(1)}_{\mu,\ell}(u,v)\allowbreak\leq\nobreak\frac{13}{64}.
\]
This, together with \(\mathbf{R}_{\ell}\leq \mathbb{P}^{(1)}_{\mu,\ell}\), implies
\begin{equation}\label{eq:large-ratio-mass}
 \sum_{\substack{u,v\in\cL(\Sigma,\ell)\\\mathbf{R}_{\ell}(u,v)>\xi_{\ell}\mathbb{P}^{(2)}_{\mu,\ell}(u,v)}}\mathbf{R}_{\ell}(u,v)
 \leq\sum_{\substack{u,v\in\cL(\Sigma,\ell)\\\mathbb{P}^{(1)}_{\mu,\ell}(u,v)>\xi_{\ell}\mathbb{P}^{(2)}_{\mu,\ell}(u,v)}}\mathbb{P}^{(1)}_{\mu,\ell}(u,v)
 \leq\frac{13}{64}.
\end{equation}
Moreover, by \eqref{eq:joint-block-pointwise} and \(\mathbf{R}_{\ell}\leq \mathbb{P}^{(1)}_{\mu,\ell}\),
\(\mathbf{R}_{\ell}=0\) wherever \(\mathbb{P}^{(2)}_{\mu,\ell}=0\).
Combining this with the definition of \(\mathscr{T}_{\ell}\),
\eqref{eq:common-block-mass}, \eqref{eq:small-ratio-mass}, and
\eqref{eq:large-ratio-mass} yields
\[
 \begin{aligned}
 \sum_{(u,v)\in \mathscr{T}_{\ell}}\mathbf{R}_{\ell}(u,v)
 &=\sum_{u,v\in\cL(\Sigma,\ell)}\mathbf{R}_{\ell}(u,v)
 -\sum_{\substack{u,v\in\cL(\Sigma,\ell)\\\mathbf{R}_{\ell}(u,v)<\frac{13}{64}\mathbb{P}^{(2)}_{\mu,\ell}(u,v)}}\mathbf{R}_{\ell}(u,v)
       -\sum_{\substack{u,v\in\cL(\Sigma,\ell)\\\mathbf{R}_{\ell}(u,v)>\xi_{\ell}\mathbb{P}^{(2)}_{\mu,\ell}(u,v)}}\mathbf{R}_{\ell}(u,v)\\
 &\geq\frac{13}{16}-\frac{13}{64}-\frac{13}{64}
       =\frac{13}{32}.
 \end{aligned}
\]
The definition of \(\mathscr{T}_{\ell}\) then gives
\begin{equation}\label{eq:pi-good-pair-mass}
\mathbb{P}^{(2)}_{\mu,\ell}(\mathscr{T}_{\ell})=\sum_{(u,v)\in \mathscr{T}_{\ell}}\mathbb{P}^{(2)}_{\mu,\ell}(u,v)
\geq\sum_{(u,v)\in \mathscr{T}_{\ell}}\frac{\mathbf{R}_{\ell}(u,v)}{\xi_{\ell}}
\geq\frac{13}{32\xi_{\ell}}.
\end{equation}

For any \(\ell\in\N\), let
\(E_{\ell}^{(\mu,\nu)}=\{(x,y)\in G_\mu\times G_\nu:
\Gamma(x,y)_{[-\ell,\ell)}\in\cL(\Sigma,2\ell)\text{ and }
\Gamma(y,x)_{[-\ell,\ell)}\in\cL(\Sigma,2\ell)\}\).
Clearly, \(E_{\mu,\nu}\subseteq E_{\ell}^{(\mu,\nu)}\) for any \(\ell\in\N\).
We next bound \((\mu\times\nu)(E_{\ell}^{(\mu,\nu)})\) from below.

Since
\(\sum_{u\in\cL(\Sigma,\ell)}\mu([u]_\Sigma)=1\) and
\(\sum_{v_1,v_2\in\cL(\Sigma,\ell)}\) \(\mu([v_1]_\Sigma)\mu([v_2]_\Sigma)
=(\sum_{v\in\cL(\Sigma,\ell)}\mu([v]_\Sigma))^2=1\),
by the Cauchy--Schwarz inequality,
\begin{equation}
\label{eq:four-cycle}
\begin{aligned}
 &\sum_{\substack{u_{1},u_{2}\in\cL(\Sigma,\ell)\\v_{1},v_{2}\in\cL(\Sigma,\ell)}}
 \mu([u_{1}]_{\Sigma})\mu([u_{2}]_{\Sigma})\mu([v_{1}]_{\Sigma})\mu([v_{2}]_{\Sigma})
          \prod_{i,j=1}^{2}\mathbf{1}_{\mathscr{T}_{\ell}}(u_{i},v_{j})\\
 &=\sum_{v_{1},v_{2}\in\cL(\Sigma,\ell)}\mu([v_{1}]_{\Sigma})\mu([v_{2}]_{\Sigma})
     \left(\sum_{u\in\cL(\Sigma,\ell)}\mu([u]_{\Sigma})\mathbf{1}_{\mathscr{T}_{\ell}}(u,v_{1})\mathbf{1}_{\mathscr{T}_{\ell}}(u,v_{2})\right)^{2}\\
 &\geq\left(\sum_{v_{1},v_{2}\in\cL(\Sigma,\ell)}\mu([v_{1}]_{\Sigma})\mu([v_{2}]_{\Sigma})
       \sum_{u\in\cL(\Sigma,\ell)}\mu([u]_{\Sigma})\mathbf{1}_{\mathscr{T}_{\ell}}(u,v_{1})\mathbf{1}_{\mathscr{T}_{\ell}}(u,v_{2})\right)^{2}\\
 &=\left(\sum_{u\in\cL(\Sigma,\ell)}\mu([u]_{\Sigma})
       \sum_{v_{1},v_{2}\in\cL(\Sigma,\ell)}\mu([v_{1}]_{\Sigma})\mu([v_{2}]_{\Sigma})
       \mathbf{1}_{\mathscr{T}_{\ell}}(u,v_{1})\mathbf{1}_{\mathscr{T}_{\ell}}(u,v_{2})\right)^{2}\\
 &=\left[\sum_{u\in\cL(\Sigma,\ell)}\mu([u]_{\Sigma})
       \left(\sum_{v\in\cL(\Sigma,\ell)}\mu([v]_{\Sigma})\mathbf{1}_{\mathscr{T}_{\ell}}(u,v)\right)^{2}\right]^{2}\\
 &\geq\left(\sum_{u\in\cL(\Sigma,\ell)}\mu([u]_{\Sigma})
       \sum_{v\in\cL(\Sigma,\ell)}\mu([v]_{\Sigma})\mathbf{1}_{\mathscr{T}_{\ell}}(u,v)\right)^{4}\\
 &=\left(\sum_{(u,v)\in \mathscr{T}_{\ell}}\mathbb{P}^{(2)}_{\mu,\ell}(u,v)\right)^{4}
 =(\mathbb{P}^{(2)}_{\mu,\ell}(\mathscr{T}_{\ell}))^{4}.
\end{aligned}
\end{equation}
By the definition of \(E_{\ell}^{(\mu,\nu)}\),
\[
 E_{\ell}^{(\mu,\nu)}
 =\bigsqcup_{\substack{u_1,u_2,v_1,v_2\in\cL(\Sigma,\ell)\\
                            u_1v_2,\,u_2v_1\in\cL(\Sigma,2\ell)}}
 (\{x\in G_\mu:x_{[-\ell,\ell)}=u_1v_1\}
       \times\{y\in G_\nu:y_{[-\ell,\ell)}=u_2v_2\}).
\]
For any \((u,v)\in \mathscr{T}_{\ell}\), the definition of \(\mathscr{T}_{\ell}\)
and \(\mathbf{R}_{\ell}\leq\mathbb{P}^{(1)}_{\mu,\ell}\) give
\(0<\frac{13}{64}\mathbb{P}^{(2)}_{\mu,\ell}(u,v)
\leq \mathbf{R}_{\ell}(u,v)\leq\mathbb{P}^{(1)}_{\mu,\ell}(u,v)=\mu([uv]_\Sigma)\),
so \([uv]_\Sigma\neq\varnothing\), and hence \(uv\in\cL(\Sigma,2\ell)\).
This, together with the preceding decomposition, the definition of \(\mathscr{T}_{\ell}\), \eqref{eq:four-cycle},
and \eqref{eq:pi-good-pair-mass}, implies
\begin{equation}\label{eq:finite-splice-mass}
 \begin{split}
 &(\mu\times\nu)(E_{\ell}^{(\mu,\nu)})
 =\sum_{\substack{u_1,u_2,v_1,v_2\in\cL(\Sigma,\ell)\\
                 u_1v_2,\,u_2v_1\in\cL(\Sigma,2\ell)}}
                     \mathbf{A}_\mu(u_1,v_1) \mathbf{A}_\nu(u_2,v_2)\\
 &\quad \geq
 \sum_{\substack{u_{1},u_{2},v_{1},v_{2}\in\cL(\Sigma,\ell)\\
                 (u_{i},v_{j})\in \mathscr{T}_{\ell}\ (1\leq i,j\leq 2)}}
                     \mathbf{A}_{\mu}(u_{1},v_{1}) \mathbf{A}_{\nu}(u_{2},v_{2})\\
 &\quad \geq
 \sum_{\substack{u_1,u_2,v_1,v_2\in\cL(\Sigma,\ell)\\
                 (u_i,v_j)\in \mathscr{T}_{\ell}\ (1\leq i,j\leq2)}}
                     \mathbf{R}_{\ell}(u_1,v_1) \mathbf{R}_{\ell}(u_2,v_2)\\
 &\quad\geq\left(\frac{13}{64}\right)^2
 \sum_{\substack{u_1,u_2,v_1,v_2\in\cL(\Sigma,\ell)\\
                 (u_i,v_j)\in \mathscr{T}_{\ell}\ (1\leq i,j\leq2)}}
 \mathbb{P}^{(2)}_{\mu,\ell}(u_1,v_1) \mathbb{P}^{(2)}_{\mu,\ell}(u_2,v_2)\\
 &\quad=\left(\frac{13}{64}\right)^2
 \sum_{u_1,u_2,v_1,v_2\in\cL(\Sigma,\ell)}
 \mu([u_1]_\Sigma)\mu([u_2]_\Sigma)\mu([v_1]_\Sigma)\mu([v_2]_\Sigma)
 \prod_{i,j=1}^2\mathbf{1}_{\mathscr{T}_{\ell}}(u_i,v_j)\\
 &\quad\geq
 \left(\frac{13}{64}\right)^{2}(\mathbb{P}^{(2)}_{\mu,\ell}(\mathscr{T}_{\ell}))^{4}
 \geq
 \left(\frac{13}{64}\right)^2
 \left(\frac{13}{32\xi_{\ell}}\right)^4=\kappa_{\ell}.
 \end{split}
\end{equation}

We finally bound
\((\mu\times\nu)(E_{\ell}^{(\mu,\nu)}\setminus E_{\mu,\nu})\) from above.

Fix any \((x,y)\in E_{\ell}^{(\mu,\nu)}\setminus E_{\mu,\nu}\).
By the definition of \(E_{\mu,\nu}\) and
\(\Sigma=\Sigma_{\mathcal{F}}\),
\[
 \begin{gathered}
 (x,y)\notin E_{\mu,\nu}
 \Longrightarrow \Gamma(x,y)\notin\Sigma
 \quad\text{or}\quad \Gamma(y,x)\notin\Sigma\\
 \Longrightarrow
 \exists n\in\N,\ \exists i\in\Z,\ \exists f\in\mathcal{F}_n:\\
 \Gamma(x,y)_{[i,i+n)}=f\quad\text{or}\quad \Gamma(y,x)_{[i,i+n)}=f.
 \end{gathered}
\]
Since \(x,y\in\Sigma\), \eqref{eq:splice-map} gives
\[
 \begin{aligned}
 i+n\leq0
 &\Longrightarrow
 \begin{cases}
 \Gamma(x,y)_{[i,i+n)}=x_{[i,i+n)}\in\cL(\Sigma,n),\\
 \Gamma(y,x)_{[i,i+n)}=y_{[i,i+n)}\in\cL(\Sigma,n),
 \end{cases}\\[2pt]
 i\geq0
 &\Longrightarrow
 \begin{cases}
 \Gamma(x,y)_{[i,i+n)}=y_{[i,i+n)}\in\cL(\Sigma,n),\\
 \Gamma(y,x)_{[i,i+n)}=x_{[i,i+n)}\in\cL(\Sigma,n).
 \end{cases}
 \end{aligned}
\]
In either case, \(f\in\cL(\Sigma,n)\), a contradiction.
Hence \(i<0<i+n\), equivalently \(1-n\leq i\leq-1\).
Moreover, \(n\geq \ell+2\). Indeed,
\((x,y)\in E_{\ell}^{(\mu,\nu)}\) implies
\(\Gamma(x,y)_{[-\ell,\ell)},\,\Gamma(y,x)_{[-\ell,\ell)}\in\cL(\Sigma,2\ell)\).
If \(n\leq \ell+1\), then
\(-\ell\leq1-n\leq i<0<i+n\leq n-1\leq \ell\),
and hence
\(\Gamma(x,y)_{[i,i+n)},\ \Gamma(y,x)_{[i,i+n)}\in\cL(\Sigma,n)\).
Since one of these blocks equals \(f\), we have \(f\in\cL(\Sigma,n)\), a contradiction.
Consequently,
\[
E_{\ell}^{(\mu,\nu)}\setminus E_{\mu,\nu}\subseteq
 \bigcup_{n=\ell+2}^{\infty}\bigcup_{f\in\mathcal{F}_n}\bigcup_{i=1-n}^{-1}
 \left[\begin{aligned}
 &\ \{(x,y)\in G_\mu\times G_\nu:\Gamma(x,y)_{[i,i+n)}=f\}\\
 &{}\cup\{(x,y)\in G_\mu\times G_\nu:\Gamma(y,x)_{[i,i+n)}=f\}
 \end{aligned}\right].
\]
For any \(n\geq \ell+2\), any \(1-n\leq i\leq-1\),
and any \(f\in\mathcal{F}_n\), \eqref{eq:splice-map} gives
\(\{(x,y)\in G_\mu\times G_\nu:\Gamma(x,y)_{[i,i+n)}=f\}
=(G_\mu\cap\{x\in\Sigma:x_{[i,0)}=f_{[0,-i)}\})
 \times(G_\nu\cap\{y\in\Sigma:y_{[0,i+n)}=f_{[-i,n)}\})\)
and
\( \{(x,y)\in G_\mu\times G_\nu:\Gamma(y,x)_{[i,i+n)}=f\}
=(G_\mu\cap\{x\in\Sigma:x_{[0,i+n)}=f_{[-i,n)}\})
 \times(G_\nu\cap\{y\in\Sigma:y_{[i,0)}=f_{[0,-i)}\}).\)
These, together with \eqref{eq:half-cylinder}, give
\begin{equation*}
 \begin{split}
 &(\mu\times\nu)(\{(x,y)\in G_\mu\times G_\nu:\Gamma(x,y)_{[i,i+n)}=f\})\\
 &\quad=
 \mu(G_{\mu}\cap\{x\in\Sigma:x_{[i,0)}=f_{[0,-i)}\})
 \nu(G_{\nu}\cap\{y\in\Sigma:y_{[0,i+n)}=f_{[-i,n)}\})\\
 &\quad\leq(C_0e^{-t^*(-i)})(C_0e^{-t^*(i+n)})
       =C_0^2 e^{-t^*n},
 \end{split}
\end{equation*}
and
\begin{equation*}
 \begin{split}
 &(\mu\times\nu)(\{(x,y)\in G_\mu\times G_\nu:\Gamma(y,x)_{[i,i+n)}=f\})\\
 &\quad=
 \mu(G_{\mu}\cap\{x\in\Sigma:x_{[0,i+n)}=f_{[-i,n)}\})
 \nu(G_{\nu}\cap\{y\in\Sigma:y_{[i,0)}=f_{[0,-i)}\})\\
 &\quad\leq(C_0e^{-t^*(i+n)})(C_0e^{-t^*(-i)})
       =C_0^2 e^{-t^*n},
 \end{split}
\end{equation*}
implying that
\begin{equation}\label{eq:finite-splice-error}
 \begin{aligned}
 &(\mu\times\nu)(E_{\ell}^{(\mu,\nu)}\setminus E_{\mu,\nu})\\
 &\quad\leq\sum_{n=\ell+2}^{\infty}\sum_{f\in\mathcal{F}_n}\sum_{i=1-n}^{-1}
 \left[\begin{aligned}
 &\ (\mu\times\nu)(\{(x,y)\in G_\mu\times G_\nu:\Gamma(x,y)_{[i,i+n)}=f\})\\
 &{}+(\mu\times\nu)(\{(x,y)\in G_\mu\times G_\nu:\Gamma(y,x)_{[i,i+n)}=f\})
 \end{aligned}\right]\\
 &\quad\leq\sum_{n=\ell+2}^{\infty}\sum_{f\in\mathcal{F}_n}\sum_{i=1-n}^{-1}
        2 C_0^2 e^{-t^*n}\\
 &\quad=2 C_0^2\sum_{n=\ell+2}^{\infty}(n-1) \#\mathcal{F}_n e^{-t^*n}.
 \end{aligned}
\end{equation}

By~\eqref{eq:splice-normalization-limits} and \(s^*<t^*\), there exists
\(\ell_0\in\N\) such that, for all integers \(\ell\geq \ell_0\),
\(\frac{\ell+1}{\kappa_{\ell}}\leq e^{(t^*-s^*)\ell/2}\).
Using \(\#\mathcal{F}_n\leq C_{\mathcal{F}} e^{s^*n}\)
yields, for all integers \(\ell\geq \ell_0\),
\[
 \begin{aligned}
 &0\leq\frac{2\cdot C_0^2}{\kappa_{\ell}}\sum_{n>\ell+1}(n-1)
 \#\mathcal{F}_n e^{-t^*n}
 \leq\frac{2 C_0^2 C_{\mathcal{F}}}{\kappa_{\ell}}
 \sum_{n>\ell+1}(n-1) e^{-(t^*-s^*)n}\\
 &\quad=\frac{2\cdot C_0^2\cdot C_{\mathcal{F}}}{\kappa_{\ell}}e^{-(t^*-s^*)(\ell+2)}
   \left(\frac{\ell+1}{1-e^{-(t^*-s^*)}}
   +\frac{e^{-(t^*-s^*)}}{(1-e^{-(t^*-s^*)})^2}\right)\\
 &\quad\leq\frac{2 C_0^2 C_{\mathcal{F}}}{(1-e^{-(t^*-s^*)})^2}
   \frac{\ell+1}{\kappa_{\ell}}e^{-(t^*-s^*)(\ell+2)}\\
 &\quad\leq\frac{2 C_0^2 C_{\mathcal{F}}}{(1-e^{-(t^*-s^*)})^2}
   e^{-(t^*-s^*)(\ell/2+2)}\rightarrow0\quad(\ell\to\infty),
 \end{aligned}
\]
implying that there exists an integer \(\ell^{*}\geq \ell_0\) such that
\begin{equation*}
 2 C_0^2\sum_{n>\ell^{*}+1}(n-1)
 \#\mathcal{F}_n e^{-t^*n}<\frac{\kappa_{\ell^{*}}}{2}.
\end{equation*}
For any pair of ergodic measures of maximal entropy \(\mu,\nu\) on \(\Sigma\) with
\(\TV(\mathbb{P}^{(1)}_{\mu,\ell^{*}},\mathbb{P}^{(1)}_{\nu,\ell^{*}})\leq\frac{1}{16}\), this,
together with \(E_{\mu,\nu}\subseteq E_{\ell^{*}}^{(\mu,\nu)}\),
\eqref{eq:finite-splice-mass}, and \eqref{eq:finite-splice-error}, implies
 \begin{align*}
 &(\mu\times\nu)(E_{\mu,\nu})=(\mu\times\nu)(E_{\ell^{*}}^{(\mu,\nu)})
 -(\mu\times\nu)(E_{\ell^{*}}^{(\mu,\nu)}\setminus E_{\mu,\nu})\\
 &\quad\geq \kappa_{\ell^{*}}-2
 C_0^2\sum_{n=\ell^{*}+2}^{\infty}(n-1) \#\mathcal{F}_n e^{-t^*n}\\
 &\quad \geq \kappa_{\ell^{*}}-\frac{\kappa_{\ell^{*}}}{2}=\frac{\kappa_{\ell^{*}}}{2}>0.\qedhere
 \end{align*}
\end{proof}

Fix \(\ell^{*}\in\N\) as in Claim~\ref{clm:two-splices}.
For any \(n\in\N\), define
\(\Upsilon_n:(\mathcal{A}^{\Z})^2\to(\mathcal{A}^{\Z})^2\)
by \(\Upsilon_n(x,y)=(x',y')\), where
\begin{equation}\label{eq:block-exchange-map}
 (x'_{i},y'_{i})=
 \begin{cases}
 (y_{i},x_{i}),&i\in \Z\cap [0, n),\\
 (x_{i},y_{i}),&i\in\Z\setminus [0,n).
 \end{cases}
\end{equation}
For any \(n\in\N\), it can be verified that \(\Upsilon_n\) is continuous and
satisfies \(\Upsilon_n\circ\Upsilon_n=\mathrm{id}\).
Put \(D_n=\{(x,y)\in\Sigma^2:\Upsilon_n(x,y)\in\Sigma^2\}
=\Sigma^2\cap\Upsilon_n^{-1}(\Sigma^2)\).
Then \(D_n\) is a closed subset of the compact set \(\Sigma^2\), and hence
\(D_n\) is compact.
Moreover, \(\Upsilon_n\circ\Upsilon_n=\mathrm{id}\) implies
\(\Upsilon_n(D_n)=D_n\), so
\(\Upsilon_n|_{D_n}:D_n\to D_n\) is a homeomorphism
with itself as inverse.
\begin{claim}\label{clm:positive-exchange}
For any pair of ergodic measures of maximal entropy \(\mu,\nu\) on \(\Sigma\) with
\(\TV(\mathbb{P}^{(1)}_{\mu,\ell^{*}},\mathbb{P}^{(1)}_{\nu,\ell^{*}})\leq\frac{1}{16}\),
we have
\(\liminf_{N\rightarrow\infty}\frac{1}{N}
\sum_{n=1}^{N}(\mu\times\nu)(D_{n})\geq\frac{\kappa_{\ell^{*}}^2}{4}>0\).
\end{claim}

\begin{proof}[Proof of Claim~\ref{clm:positive-exchange}]
Fix any pair of ergodic measures of maximal entropy \(\mu,\nu\) on \(\Sigma\) with
\(\TV(\mathbb{P}^{(1)}_{\mu,\ell^{*}},\mathbb{P}^{(1)}_{\nu,\ell^{*}})\leq\frac{1}{16}\).
For all \(N\in\N\), we have
\begin{equation}\label{eq:exchange-average-comparison}
 \begin{aligned}
 \frac{1}{N}\sum_{n=1}^N(\mu\times\nu)(D_n)
 &\geq\frac{1}{N}\sum_{n=1}^N(\mu\times\nu)
       ((E_{\mu,\nu}\cap(\sigma\times\sigma)^{-n}(E_{\mu,\nu}))
       \cap D_n)\\
 &=\frac{1}{N}\sum_{n=1}^N(\mu\times\nu)
       (E_{\mu,\nu}\cap(\sigma\times\sigma)^{-n}(E_{\mu,\nu}))\\
 &\quad-\frac{1}{N}\sum_{n=1}^N(\mu\times\nu)(
       (E_{\mu,\nu}\cap(\sigma\times\sigma)^{-n}(E_{\mu,\nu}))
       \setminus D_n).
 \end{aligned}
\end{equation}

We first estimate the first term on the right-hand side
of \eqref{eq:exchange-average-comparison}.

Since \(\mu\) and \(\nu\) are \(\sigma\)-invariant,
\(\mu\times\nu\) is \((\sigma\times\sigma)\)-invariant.
Applying the \(L^p\) Ergodic Theorem of Von Neumann
(\cite[Corollary~1.14.1]{Walters1982}, with \(p=2\)) to
\(\mathbf{1}_{E_{\mu,\nu}}\)
for the measure-preserving transformation
\((\sigma\times\sigma)|_{\Sigma^2}\) on \(\Sigma^2\) with its Borel
\(\sigma\)-algebra and invariant probability measure \(\mu\times\nu\)
yields a real-valued \(f^*\in L^2(\mu\times\nu)\) such that
\begin{equation}\label{eq:ergodic-limit-invariance}
f^*((\sigma\times\sigma)(x,y))=f^*(x,y) \quad\text{for
\((\mu\times\nu)\)-a.e. } (x,y)\in\Sigma^2,
\end{equation}
and
\[
 \lim_{N\to\infty}\left\|
 \frac{1}{N}\sum_{n=1}^{N}
 \mathbf{1}_{E_{\mu,\nu}}((\sigma\times\sigma)^n(x,y))
 -f^*(x,y)
 \right\|_{L^2(\mu\times\nu)}=0.
\]
These, together with the Cauchy--Schwarz inequality, imply that, for each
\(g\in\{1,\mathbf{1}_{E_{\mu,\nu}},f^*\}\),
\[
 \begin{aligned}
 &\left|\int_{\Sigma^2}
       \left(\frac{1}{N}\sum_{n=1}^{N}
       \mathbf{1}_{E_{\mu,\nu}}((\sigma\times\sigma)^n(x,y))
       -f^*(x,y)\right)\cdot g(x,y)\,d(\mu\times\nu)(x,y)\right|\\
 &\quad\leq
       \left\|\frac{1}{N}\sum_{n=1}^{N}
       \mathbf{1}_{E_{\mu,\nu}}((\sigma\times\sigma)^n(x,y))
       -f^*(x,y)\right\|_{L^2(\mu\times\nu)}
       \cdot\|g\|_{L^2(\mu\times\nu)}
       \rightarrow0\quad(N\to\infty).
 \end{aligned}
\]
Taking \(g=1\) and using the \((\sigma\times\sigma)\)-invariance of \(\mu\times\nu\) gives
\[
 \int_{\Sigma^2}f^*(x,y)\,d(\mu\times\nu)(x,y)
 =\lim_{N\to\infty}\frac{1}{N}\sum_{n=1}^{N}
       \int_{\Sigma^2}\mathbf{1}_{E_{\mu,\nu}}
       ((\sigma\times\sigma)^n(x,y))
       \,d(\mu\times\nu)(x,y)
 =(\mu\times\nu)(E_{\mu,\nu}).
\]
Taking \(g=f^*\) and using the \((\sigma\times\sigma)\)-invariance of \(\mu\times\nu\),
together with \eqref{eq:ergodic-limit-invariance}, gives
\[
 \begin{aligned}
 &\int_{\Sigma^2}\mathbf{1}_{E_{\mu,\nu}}(x,y)
       \cdot f^*(x,y)\,d(\mu\times\nu)(x,y)\\
 &\quad=\lim_{N\to\infty}\frac{1}{N}\sum_{n=1}^{N}
       \int_{\Sigma^2}\mathbf{1}_{E_{\mu,\nu}}((\sigma\times\sigma)^n(x,y))
       \cdot f^*((\sigma\times\sigma)^n(x,y))\,d(\mu\times\nu)(x,y)\\
 &\quad=\lim_{N\to\infty}\int_{\Sigma^2}
       \left(\frac{1}{N}\sum_{n=1}^{N}
       \mathbf{1}_{E_{\mu,\nu}}((\sigma\times\sigma)^n(x,y))\right)
       \cdot f^*(x,y)\,d(\mu\times\nu)(x,y)\\
 &\quad=\int_{\Sigma^2}(f^*(x,y))^2\,d(\mu\times\nu)(x,y).
 \end{aligned}
\]
Taking \(g=\mathbf{1}_{E_{\mu,\nu}}\) and using these two identities,
the Cauchy--Schwarz inequality, \((\mu\times\nu)(\Sigma^2)=1\),
and Claim~\ref{clm:two-splices} gives
\begin{equation}
\label{eq:recurrence-product}
 \begin{split}
 &\lim_{N\rightarrow\infty}\frac{1}{N}\sum_{n=1}^{N}
       (\mu\times\nu)
       (E_{\mu,\nu}\cap(\sigma\times\sigma)^{-n}(E_{\mu,\nu}))\\
 &\quad=\lim_{N\to\infty}\int_{\Sigma^2}\mathbf{1}_{E_{\mu,\nu}}(x,y)
       \cdot\left(\frac{1}{N}\sum_{n=1}^{N}
       \mathbf{1}_{E_{\mu,\nu}}((\sigma\times\sigma)^n(x,y))\right)
       \,d(\mu\times\nu)(x,y)\\
 &\quad=\int_{\Sigma^2}\mathbf{1}_{E_{\mu,\nu}}(x,y)
       \cdot f^*(x,y)\,d(\mu\times\nu)(x,y)\\
 &\quad=\int_{\Sigma^2}(f^*(x,y))^2\,d(\mu\times\nu)(x,y)\\
 &\quad\geq\left[\int_{\Sigma^2}f^*(x,y)
       \,d(\mu\times\nu)(x,y)\right]^2\\
 &\quad=((\mu\times\nu)(E_{\mu,\nu}))^2\geq\frac{\kappa_{\ell^{*}}^2}{4}.
 \end{split}
\end{equation}

We next estimate the second term on the right-hand side of
\eqref{eq:exchange-average-comparison}.

Put \(\eta=\frac{1}{4}(\htop(\Sigma)-s^*)\).
The choices of \(s^*\) and \(t^*\) ensure that
\(\htop(\Sigma)<\htop(\Sigma)+\eta<2t^*-s^*\).
By \eqref{eq:language-limit}, there exists \(C_\eta\geq1\) such that
\(\#\cL(\Sigma,n)\leq C_\eta\cdot e^{(\htop(\Sigma)+\eta)n}\)
for all \(n\in\N\).

Fix any \(n\in\N\), and take any
\((x,y)\in(E_{\mu,\nu}\cap(\sigma\times\sigma)^{-n}(E_{\mu,\nu}))\setminus D_n\).
Since \((x,y)\notin D_n\), the definition of \(D_n\) gives
\(\Upsilon_n(x,y)\notin\Sigma^2\), and thus
\((\Upsilon_n(x,y))_1\notin\Sigma\) or
\((\Upsilon_n(x,y))_2\notin\Sigma.\)
Since \(\Sigma=\Sigma_{\mathcal{F}}\), there exist
\(m\in\N\), \(i\in\Z\), and \(f\in\mathcal{F}_m\) such that
\(f=((\Upsilon_n(x,y))_1)_{[i,i+m)}\)
or
\(f=((\Upsilon_n(x,y))_2)_{[i,i+m)}\), i.e.,
\(f\in\{((\Upsilon_n(x,y))_1)_{[i,i+m)},\allowbreak
((\Upsilon_n(x,y))_2)_{[i,i+m)}\}\).
Moreover, since
\((x,y)\in E_{\mu,\nu}\cap(\sigma\times\sigma)^{-n}(E_{\mu,\nu})\),
the definition of \(E_{\mu,\nu}\) and the
\(\sigma\)-invariance of \(\Sigma\) give
\(\Gamma(x,y),\allowbreak
\Gamma(y,x)\in\Sigma\) and
\(\sigma^{-n}(\Gamma(\sigma^n(y),\sigma^n(x))),\allowbreak
\sigma^{-n}(\Gamma(\sigma^n(x),\sigma^n(y)))\allowbreak\in\Sigma\),
implying that
\begin{itemize}
\setlength{\emergencystretch}{3em}
\item If \(i+m\leq n\), then, noting that \(((\Upsilon_{n}(x,y))_{1})_{(-\infty, n)}
\allowbreak=(\Gamma(x, y))_{(-\infty, n)}\) and \(((\Upsilon_{n}(x,y))_{2})_{(-\infty, n)}
\allowbreak=(\Gamma(y, x))_{(-\infty, n)}\) by \eqref{eq:block-exchange-map} and \eqref{eq:splice-map}, we have
\(f\in\{((\Upsilon_n(x,y))_1)_{[i,i+m)},\;\allowbreak
       ((\Upsilon_n(x,y))_2)_{[i,i+m)}\}
\allowbreak=\{\Gamma(x,y)_{[i,i+m)},\allowbreak
\Gamma(y,x)_{[i,i+m)}\}
       \allowbreak\subseteq\cL(\Sigma,m);\)
\item If \(i\geq0\), then, noting that
\(((\Upsilon_{n}(x,y))_{1})_{[0,\infty)}
\allowbreak=(\sigma^{-n}(\Gamma(\sigma^n(y),\sigma^n(x))))_{[0,\infty)}\)
and \(((\Upsilon_{n}(x,y))_{2})_{[0,\infty)}
\allowbreak=(\sigma^{-n}(\Gamma(\sigma^n(x),\sigma^n(y))))_{[0,\infty)}\)
by \eqref{eq:block-exchange-map} and \eqref{eq:splice-map}, we have
\(f\in\{((\Upsilon_n(x,y))_1)_{[i,i+m)},\;\allowbreak
       ((\Upsilon_n(x,y))_2)_{[i,i+m)}\}
  \allowbreak=\{(\sigma^{-n}(\Gamma(\sigma^n(y),\sigma^n(x))))_{[i,i+m)},\;\allowbreak
  (\sigma^{-n}(\Gamma(\sigma^n(x),\sigma^n(y))))_{[i,i+m)}\}\allowbreak\subseteq\cL(\Sigma,m).\)
\end{itemize}
These, together with \(f\in\mathcal{F}_m\subseteq\mathcal{A}^m\setminus\cL(\Sigma,m)\),
imply that \(i+m>n\) and \(i<0\), and thus \(m\geq n+2\) and
\(n-m+1\leq i\leq-1\).
Therefore, by \eqref{eq:block-exchange-map},
\(f\in\{
x_{[i,0)}y_{[0,n)}x_{[n,i+m)},\;
y_{[i,0)}x_{[0,n)}y_{[n,i+m)}\}\).
Since \((x,y)\) is arbitrary in
\((E_{\mu,\nu}\cap(\sigma\times\sigma)^{-n}(E_{\mu,\nu}))\setminus D_n\) and
\(E_{\mu,\nu}\subseteq G_\mu\times G_\nu\),
\begin{equation}\label{eq:exchange-forbidden-cover}
 \begin{aligned}
 &(E_{\mu,\nu}\cap(\sigma\times\sigma)^{-n}(E_{\mu,\nu}))
       \setminus D_n\\
 &\quad\subseteq
       \bigcup_{m\geq n+2}\ \bigcup_{f\in\mathcal{F}_m}\
       \bigcup_{i=n-m+1}^{-1}
       \Bigl\{(x,y)\in G_\mu\times G_\nu:\\
 &\hspace{5em}
       f\in\{x_{[i,0)}y_{[0,n)}x_{[n,i+m)},\;
                 y_{[i,0)}x_{[0,n)}y_{[n,i+m)}\}\Bigr\}.
 \end{aligned}
\end{equation}
For any \(m\geq n+2\), any \(n-m+1\leq i\leq -1\), and any
\(f\in\mathcal{F}_m\), \eqref{eq:good-cylinder} and
\eqref{eq:half-cylinder} give
\[
 \begin{aligned}
 &(\mu\times\nu)(\{(x,y)\in G_\mu\times G_\nu:x_{[i,0)}y_{[0,n)}x_{[n,i+m)}=f\})\\
 &\quad\leq\nu(G_{\nu}\cap[f_{[-i,n-i)}]_{\Sigma})
 \sum_{w\in\cL(\Sigma,n)}
          \mu(G_{\mu}\cap
          \{x\in\Sigma:x_{[i,i+m)}=f_{[0,-i)}wf_{[n-i,m)}\})\\
 &\quad\leq C_0 e^{-t^*n}\sum_{w\in\cL(\Sigma,n)}C_0 e^{-t^*m}
       =C_0^2 \#\cL(\Sigma,n)e^{-t^*(m+n)}
 \end{aligned}
\]
and
\[
 \begin{aligned}
 &(\mu\times\nu)(\{(x,y)\in G_\mu\times G_\nu:y_{[i,0)}x_{[0,n)}y_{[n,i+m)}=f\})\\
 &\quad\leq\mu(G_{\mu}\cap[f_{[-i,n-i)}]_{\Sigma})
 \sum_{w\in\cL(\Sigma,n)}
          \nu(G_{\nu}\cap
          \{y\in\Sigma:y_{[i,i+m)}=f_{[0,-i)}wf_{[n-i,m)}\})\\
 &\quad\leq C_0 e^{-t^*n}\sum_{w\in\cL(\Sigma,n)}C_0 e^{-t^*m}
       =C_0^2 \#\cL(\Sigma,n) e^{-t^*(m+n)}.
 \end{aligned}
\]
These, together with \eqref{eq:exchange-forbidden-cover},
\(\#\cL(\Sigma,n)\leq C_\eta\cdot e^{(\htop(\Sigma)+\eta)n}\),
\(\#\mathcal{F}_m\leq C_{\mathcal{F}} e^{s^*m}\),
and \(t^*>s^*\), imply
\[
 \begin{aligned}
 &0\leq(\mu\times\nu)
   ((E_{\mu,\nu}\cap(\sigma\times\sigma)^{-n}(E_{\mu,\nu}))\setminus D_n)\\
 &\quad\leq\sum_{m\geq n+2}\sum_{f\in\mathcal{F}_m}\sum_{i=n-m+1}^{-1}
 \left[\begin{aligned}
 &\ (\mu\times\nu)(\{(x,y)\in G_\mu\times G_\nu:x_{[i,0)}y_{[0,n)}x_{[n,i+m)}=f\})\\
 &{}+(\mu\times\nu)(\{(x,y)\in G_\mu\times G_\nu:y_{[i,0)}x_{[0,n)}y_{[n,i+m)}=f\})
 \end{aligned}\right]\\
 &\quad\leq\sum_{m\geq n+2}\sum_{f\in\mathcal{F}_m}\sum_{i=n-m+1}^{-1}
       2 C_0^2 \#\cL(\Sigma,n) e^{-t^*(m+n)}\\
 &\quad=2 C_0^2 \#\cL(\Sigma,n) e^{-t^*n}
       \sum_{m\geq n+2}(m-n-1) \#\mathcal{F}_m e^{-t^*m}\\
 &\quad\leq 2 C_0^2 (C_\eta e^{(\htop(\Sigma)+\eta)n}) e^{-t^*n}
 \sum_{m\geq n+2}(m-n-1)
       (C_{\mathcal{F}} e^{s^*m}) e^{-t^*m}\\
 &\quad=\frac{2 C_0^2 C_\eta C_{\mathcal{F}} e^{-2(t^*-s^*)}}
                  {(1-e^{-(t^*-s^*)})^2}
       e^{-(2t^*-\htop(\Sigma)-s^*-\eta)n},
 \end{aligned}
\]
and thus, for all \(N\in\N\),
\begin{equation}\label{eq:exchange-error-average}
 \begin{aligned}
 0&\leq\frac{1}{N}\sum_{n=1}^N(\mu\times\nu)
 ((E_{\mu,\nu}\cap(\sigma\times\sigma)^{-n}(E_{\mu,\nu}))
       \setminus D_n)\\
 &\leq\frac{2 C_0^2 C_\eta C_{\mathcal{F}} e^{-2(t^*-s^*)}}
             {N(1-e^{-(t^*-s^*)})^2}
       \sum_{n=1}^Ne^{-(2t^*-\htop(\Sigma)-s^*-\eta)n}\\
 &\leq\frac{2 C_0^2 C_\eta C_{\mathcal{F}} e^{-2(t^*-s^*)}}
             {N(1-e^{-(t^*-s^*)})^2
              (e^{2t^*-\htop(\Sigma)-s^*-\eta}-1)}
       \rightarrow0\quad(N\to\infty).
 \end{aligned}
\end{equation}

Combining \eqref{eq:exchange-average-comparison},
\eqref{eq:recurrence-product}, and \eqref{eq:exchange-error-average} yields
\begin{align*}
 &\liminf_{N\to\infty}\frac{1}{N}
       \sum_{n=1}^N(\mu\times\nu)(D_n)\\
 &\quad\geq\lim_{N\to\infty}\frac{1}{N}\sum_{n=1}^N
       (\mu\times\nu)
       (E_{\mu,\nu}\cap(\sigma\times\sigma)^{-n}(E_{\mu,\nu}))\\
 &\qquad-\lim_{N\to\infty}\frac{1}{N}\sum_{n=1}^N
       (\mu\times\nu)(
       (E_{\mu,\nu}\cap(\sigma\times\sigma)^{-n}(E_{\mu,\nu}))
       \setminus D_n)\\
 &\quad\geq\frac{\kappa_{\ell^{*}}^2}{4}>0.\qedhere
\end{align*}
\end{proof}

For any pair of ergodic measures of maximal entropy
\(\mu,\nu\) on \(\Sigma\), the definitions of \(\mathbb{P}^{(1)}_{\mu,\ell^{*}}\)
and \(\mathbb{P}^{(1)}_{\nu,\ell^{*}}\),
together with \eqref{eq:total-variation}, give
\[
 \frac{1}{2}\sum_{w\in\cL(\Sigma,2\ell^{*})}
       |\mu([w]_\Sigma)-\nu([w]_\Sigma)|
 =\frac{1}{2}\sum_{u,v\in\cL(\Sigma,\ell^{*})}
       |\mu([uv]_\Sigma)-\nu([uv]_\Sigma)|
       =\TV(\mathbb{P}^{(1)}_{\mu,\ell^{*}},\mathbb{P}^{(1)}_{\nu,\ell^{*}}).
\]
To complete the proof of
Lemma~\ref{lem:positive-entropy-separation}, suppose, to the contrary, that there exist distinct
ergodic measures of maximal entropy \(\mu_0,\nu_0\) on \(\Sigma\) with
\(\TV(\mathbb{P}^{(1)}_{\mu_0,\ell^{*}},\mathbb{P}^{(1)}_{\nu_0,\ell^{*}})\leq\frac{1}{16}\).

Define the probability measure \(\tilde{\Lambda}(\mu_0\times\nu_0)\)
on \(\mathcal{X}_{\Sigma}\) by
\((\tilde{\Lambda}(\mu_0\times\nu_0))(B)=(\mu_0\times\nu_0)(\Lambda^{-1}(B))\)
for any Borel set \(B\subseteq\mathcal{X}_{\Sigma}\).
By this definition, the measure-preserving systems
\((\Sigma^2,\mu_0\times\nu_0,(\sigma\times\sigma)|_{\Sigma^2})\) and
\((\mathcal{X}_{\Sigma},\tilde{\Lambda}(\mu_0\times\nu_0),\sigma|_{\mathcal{X}_{\Sigma}})\)
are isomorphic via the conjugacy \(\Lambda\).
Since \(\mu_0\) and \(\nu_0\) are measures of maximal entropy,
\cite[Theorem~4.23]{Walters1982}~gives
\begin{equation}\label{eq:product-pushforward-entropy}
 h_{\tilde{\Lambda}(\mu_0\times\nu_0)}(\sigma|_{\mathcal{X}_{\Sigma}})
 \mathrel{=}
 h_{\mu_0\times\nu_0}((\sigma\times\sigma)|_{\Sigma^2})
 =h_{\mu_0}(\sigma|_\Sigma)+h_{\nu_0}(\sigma|_\Sigma)
 =2\htop(\Sigma)=\htop(\mathcal{X}_{\Sigma}).
\end{equation}

Fix any \(n\in\N\).
For any \(u,v\in\cL(\Sigma,n)\), put
\(\Xi_{u,v}:=(u_0,v_0)\cdots(u_{n-1},v_{n-1})
\in\cL(\mathcal{X}_{\Sigma},n)\).
For distinct \(u,v\in\cL(\Sigma,n)\),
by the definitions of $\Upsilon_{n}$,
$\Lambda$, $D_{n}$, and $D_{\mathcal{X}_{\Sigma}}$,
it can be verified that
\begin{equation}\label{eq:product-replacement-domain}
 \Lambda(D_n\cap([u]_\Sigma\times[v]_\Sigma))
 \subseteq D_{\mathcal{X}_{\Sigma}}(\Xi_{u,v},\Xi_{v,u}),
\end{equation}
and that, for any \((x,y)\in D_n\cap([u]_\Sigma\times[v]_\Sigma)\),
\begin{equation}\label{eq:product-exchange-conjugacy}
 \Lambda(\Upsilon_n(x,y))
 =\Phi^{\mathcal{X}_{\Sigma}}_{\Xi_{u,v},\Xi_{v,u}}(\Lambda(x,y)).
\end{equation}

By Proposition~\ref{prop:product-forbidden-entropy}, we have
\(\hG(\mathcal{X}_{\Sigma})<\htop(\mathcal{X}_{\Sigma})\).
Since \(\htop(\mathcal{X}_{\Sigma})=2\htop(\Sigma)>0\),
the hypotheses of Lemma~\ref{lem:finite-replacement-invariance} are satisfied
for \(\mathcal{X}_{\Sigma}\) and \(\tilde{\Lambda}(\mu_0\times\nu_0)\).
For any \(u,v\in\cL(\Sigma,n)\) and any Borel set
\(B\subseteq D_n\cap([u]_\Sigma\times[v]_\Sigma)\),
we prove that
\((\mu_0\times\nu_0)(\Upsilon_n(B))=(\mu_0\times\nu_0)(B)\)
by considering the following two cases:
\begin{itemize}
\item If \(u=v\), then
\(\Upsilon_n(B)=B\), and hence
\((\mu_0\times\nu_0)(\Upsilon_n(B))=(\mu_0\times\nu_0)(B)\).
\item If \(u\ne v\), then
\(\Lambda(B)\) is a Borel subset of
\(D_{\mathcal{X}_{\Sigma}}(\Xi_{u,v},\Xi_{v,u})\)
by \eqref{eq:product-replacement-domain} and the fact that
\(\Lambda\) is a homeomorphism.
This, together with Lemma~\ref{lem:finite-replacement-invariance} and
\eqref{eq:product-exchange-conjugacy}, implies
\[
 \begin{aligned}
 &(\mu_0\times\nu_0)(\Upsilon_n(B))
 =(\tilde{\Lambda}(\mu_0\times\nu_0))(\Lambda(\Upsilon_n(B)))\\
 &\quad =(\tilde{\Lambda}(\mu_0\times\nu_0))
   (\Phi^{\mathcal{X}_{\Sigma}}_{\Xi_{u,v},\Xi_{v,u}}(\Lambda(B)))\\
 &\quad =(\tilde{\Lambda}(\mu_0\times\nu_0))(\Lambda(B))\\
 &\quad =(\mu_0\times\nu_0)(B).
 \end{aligned}
\]
\end{itemize}
Noting that
\(D_n=D_n\cap \Sigma^{2}=\bigsqcup_{u,v\in\cL(\Sigma,n)}
(D_n\cap([u]_\Sigma\times[v]_\Sigma))\), we obtain
by the injectivity of \(\Upsilon_n\), for any Borel set
\(B\subseteq D_n\),
\begin{equation}
\label{eq:swap-invariance}
 \begin{split}
 & (\mu_0\times\nu_0)(\Upsilon_{n}(B))
 =\sum_{u,v\in\cL(\Sigma,n)}
       (\mu_0\times\nu_0)(\Upsilon_{n}(B\cap([u]_{\Sigma}\times[v]_{\Sigma})))\\
 &\quad =\sum_{u,v\in\cL(\Sigma,n)}
       (\mu_0\times\nu_0)(B\cap([u]_{\Sigma}\times[v]_{\Sigma}))
       =(\mu_0\times\nu_0)(B).
 \end{split}
\end{equation}

Since \(\mu_0\neq\nu_0\), there exist \(r\in\N\) and
\(w\in\cL(\Sigma,r)\) such that
\(\mu_0([w]_\Sigma)\neq\nu_0([w]_\Sigma)\).
Choose disjoint open intervals \(I_{\mu_0},I_{\nu_0}\) containing
\(\mu_0([w]_{\Sigma})\) and \(\nu_0([w]_{\Sigma})\), respectively.
For any integer \(n\geq r\), let
\[
 \Omega_n=\biggl\{(x,y)\in\Sigma^2:\
 \frac{1}{n-r+1}\sum_{j=0}^{n-r}
       \mathbf{1}_{[w]_{\Sigma}}(\sigma^{j}(x))\in I_{\mu_0},\
 \frac{1}{n-r+1}\sum_{j=0}^{n-r}
       \mathbf{1}_{[w]_{\Sigma}}(\sigma^{j}(y))\in I_{\nu_0}
                                                     \biggr\},
 \]
and
\[
 \tilde{\Omega}_n=\biggl\{(x,y)\in\Sigma^2:\
 \frac{1}{n-r+1}\sum_{j=0}^{n-r}
       \mathbf{1}_{[w]_{\Sigma}}(\sigma^{j}(x))\in I_{\nu_0},\
 \frac{1}{n-r+1}\sum_{j=0}^{n-r}
       \mathbf{1}_{[w]_{\Sigma}}(\sigma^{j}(y))\in I_{\mu_0}
                                                     \biggr\}.
 \]
Then \(\tilde{\Omega}_n=\{(x,y)\in\Sigma^2:(y,x)\in\Omega_n\}.\)
Since \(I_{\mu_0}\cap I_{\nu_0}=\varnothing\), we have \(\Omega_n\cap\tilde{\Omega}_n=\varnothing\).
The Birkhoff Ergodic Theorem
(\cite[Theorem~1.14]{Walters1982}), applied to
\(\mathbf{1}_{[w]_\Sigma}\) with respect to \(\mu_0\) and \(\nu_0\) separately,
gives
\[
\lim_{n\to\infty}\frac{1}{n-r+1}\sum_{j=0}^{n-r}\mathbf{1}_{[w]_\Sigma}(\sigma^j(x))
=\mu_0([w]_\Sigma) \quad\text{for \(\mu_0\)-a.e. \(x\in\Sigma\)},
\]
 and
\[
\lim_{n\to\infty}\frac{1}{n-r+1}\sum_{j=0}^{n-r}\mathbf{1}_{[w]_\Sigma}(\sigma^j(y))
=\nu_0([w]_\Sigma) \quad\text{for \(\nu_0\)-a.e. \(y\in\Sigma\)}.
\]
Since \(I_{\mu_0}\) and \(I_{\nu_0}\) are open neighborhoods of these respective limits,
\[
\lim_{n\to\infty}\mathbf{1}_{\Omega_n}(x,y)=1
\quad\text{for \((\mu_0\times\nu_0)\)-a.e. \((x,y)\in\Sigma^2\)}.
\]
By \(\Omega_n\cap\tilde{\Omega}_n=\varnothing\), we have
\(0\leq\mathbf{1}_{\tilde{\Omega}_n}\leq1-\mathbf{1}_{\Omega_n}\),
and thus
\[
\lim_{n\to\infty}\mathbf{1}_{\tilde{\Omega}_n}(x,y)=0
\quad\text{for \((\mu_0\times\nu_0)\)-a.e. \((x,y)\in\Sigma^2\)}.
\]
The Dominated Convergence Theorem
(\cite[Corollary~0.9.1]{Walters1982}) gives
\begin{equation}\label{eq:swap-frequency-limits}
 \begin{cases}
 \displaystyle\lim_{n\to\infty}(\mu_0\times\nu_0)(\Omega_n)
 =\lim_{n\to\infty}\int_{\Sigma^2}
 \mathbf{1}_{\Omega_n}(x,y)\,d(\mu_0\times\nu_0)(x,y)
 =\int_{\Sigma^2}1\,d(\mu_0\times\nu_0)(x,y)=1,\\[8pt]
 \displaystyle\lim_{n\to\infty}(\mu_0
 \times\nu_0)(\tilde{\Omega}_n)
 =\lim_{n\to\infty}\int_{\Sigma^2}
 \mathbf{1}_{\tilde{\Omega}_n}(x,y)\,d(\mu_0\times\nu_0)(x,y)
 =\int_{\Sigma^2}0\,d(\mu_0\times\nu_0)(x,y)=0.
 \end{cases}
\end{equation}
For any integer \(n\geq r\) and any
\((x,y)\in D_n\cap\Omega_n\), we have
\(\Upsilon_n(x,y)\in\Sigma^2\).
By \eqref{eq:block-exchange-map},
\[
 ((\Upsilon_n(x,y))_1)_{[j,j+r)}=y_{[j,j+r)}\quad\text{and}\quad
 ((\Upsilon_n(x,y))_2)_{[j,j+r)}=x_{[j,j+r)}
\]
for \(0\leq j\leq n-r\).
Since \((x,y)\in\Omega_n\), these identities give
\[
\frac{1}{n-r+1}\sum_{j=0}^{n-r}
 \mathbf{1}_{[w]_\Sigma}(\sigma^j((\Upsilon_n(x,y))_2))
=\frac{1}{n-r+1}\sum_{j=0}^{n-r}
 \mathbf{1}_{[w]_\Sigma}(\sigma^j(x))\in I_{\mu_0},
\]
and
\[
\frac{1}{n-r+1}\sum_{j=0}^{n-r}
 \mathbf{1}_{[w]_\Sigma}(\sigma^j((\Upsilon_n(x,y))_1))
 =\frac{1}{n-r+1}\sum_{j=0}^{n-r}
 \mathbf{1}_{[w]_\Sigma}(\sigma^j(y))\in I_{\nu_0},
\]
implying that
\(((\Upsilon_n(x,y))_2,(\Upsilon_n(x,y))_1)\in\Omega_n\), and hence
\(\Upsilon_n(x,y)\in\tilde{\Omega}_n\).
Since \((x,y)\) is arbitrary in \(D_n\cap\Omega_n\), we have
\(\Upsilon_n(D_n\cap\Omega_n)\subseteq\tilde{\Omega}_n\).
This, together with \eqref{eq:swap-invariance} and
\eqref{eq:swap-frequency-limits}, gives
\[
 \begin{aligned}
 0&\leq(\mu_0\times\nu_0)(D_n)
 \leq(\mu_0\times\nu_0)(\Sigma^2\setminus\Omega_n)+(\mu_0\times\nu_0)(D_n\cap\Omega_n)\\
 &=(\mu_0\times\nu_0)(\Sigma^2\setminus\Omega_n)
    +(\mu_0\times\nu_0)(\Upsilon_n(D_n\cap\Omega_n))\\
 &\mathrel{\leq}1-(\mu_0\times\nu_0)(\Omega_n)
 +(\mu_0\times\nu_0)(\tilde{\Omega}_n)\rightarrow 0\quad(n\to \infty),
 \end{aligned}
\]
implying that
\(\lim_{N\to\infty}\frac{1}{N}\sum_{n=1}^{N}(\mu_0\times\nu_0)(D_n)=0\),
contradicting Claim~\ref{clm:positive-exchange}.
This proves
\eqref{eq:finite-word-separation} with \(\ell=\ell^{*}\).
\end{proof}

\begin{theorem}
\label{thm:positive-entropy}
Let \(\Sigma\) be a subshift over a finite alphabet
\(\mathcal{A}\) satisfying \(\htop(\Sigma)>0\) and
\(\hG(\Sigma)<\htop(\Sigma)\).
Then \(\Sigma\) has only finitely many ergodic measures of maximal entropy,
each of which is Bernoulli up to a finite period.
\end{theorem}

\begin{proof}
Choose
\(\ell^{*}\in\N\) as in the proof of
Lemma~\ref{lem:positive-entropy-separation}.
Write \(\cL(\Sigma,2\ell^{*})=\{w_1,\ldots,w_m\}\), where
\(m=\#\cL(\Sigma,2\ell^{*})\in\N\), and define
\[
 I_j=
 \begin{cases}
 [\frac{j}{16m}, \frac{j+1}{16m}),
       & j=0,\ldots,16m-2,\\
 [\frac{16m-1}{16m}, 1],& j=16m-1.
 \end{cases}
\]
These intervals form a partition of \([0,1]\).
Consequently, the sets \(I_{j_1}\times\cdots\times I_{j_m}\),
\(j_1,\ldots,j_m\in\{0,\ldots,16m-1\}\),
form a finite partition of \([0,1]^m\).
For any ergodic measure of maximal entropy \(\mu\) on \(\Sigma\), the vector
\((\mu([w_1]_\Sigma),\ldots,\mu([w_m]_\Sigma))\)
belongs to one of these sets.
Let \(\mu,\nu\) be any pair of ergodic measures of maximal entropy on \(\Sigma\)
whose vectors belong to the same set
\(I_{j_1}\times\cdots\times I_{j_m}\) in this partition.
Then \(\mu([w_i]_\Sigma),\nu([w_i]_\Sigma)\in I_{j_i}\) for all
\(i\in\{1,\ldots,m\}\). Since each \(I_{j_i}\) has length
\(\frac{1}{16m}\),
\[
 \frac{1}{2}\sum_{w\in\cL(\Sigma,2\ell^{*})}
 |\mu([w]_\Sigma)-\nu([w]_\Sigma)|
 =\frac{1}{2}\sum_{i=1}^{m}
 |\mu([w_i]_\Sigma)-\nu([w_i]_\Sigma)|
 \leq\frac{1}{2}\sum_{i=1}^{m}\frac{1}{16m}
 =\frac{1}{32}<\frac{1}{16}.
\]
This, together with Lemma~\ref{lem:positive-entropy-separation}, implies \(\mu=\nu\).
Thus each set in this finite partition contains the vector of at most
one ergodic measure of maximal entropy on \(\Sigma\). Therefore, \(\Sigma\) has only finitely
many ergodic measures of maximal entropy.

It remains to prove that every ergodic measure of maximal entropy is
Bernoulli up to a finite period.
Fix any ergodic measure of maximal entropy \(\mu\) on \(\Sigma\).
Recall \(\Upsilon_n\) from \eqref{eq:block-exchange-map} and
\(D_n=\{(x,y)\in\Sigma^2:\Upsilon_n(x,y)\in\Sigma^2\}\).
Applying Claim~\ref{clm:positive-exchange} with \(\nu=\mu\) gives
\(\liminf_{N\to\infty}\frac{1}{N}\sum_{n=1}^{N}(\mu\times\mu)(D_n)
\geq\frac{\kappa_{\ell^{*}}^2}{4}>0,\)
and thus
\[
\limsup_{n\to\infty}(\mu\times\mu)(D_n)\geq\frac{\kappa_{\ell^{*}}^2}{4}>0,
\]
where \(\kappa_{\ell^{*}}\) is as in Claim~\ref{clm:two-splices}.

Since \(\mu\) is \(\sigma\)-invariant, for any
\(r\in\N\), any integer \(n\geq r\), and any \(u,v\in\cL(\Sigma,r)\),
\begin{align*}
 &(\mu\times\mu)(\{(x,y)\in D_n:
 x_{[-r,0)}=u,\ y_{[0,r)}=v\})\\
 & \quad \leq
 \mu(\{x\in \Sigma: x_{[-r,0)}=u\})
 \mu(\{y\in \Sigma: y_{[0,r)}=v\})\\
 & \quad =\mu(\sigma^{r}([u]_{\Sigma}))
 \mu([v]_{\Sigma})=\mu([u]_{\Sigma}) \mu([v]_{\Sigma})=
 \mathbb{P}^{(2)}_{\mu,r}(u,v).
\end{align*}

On the other hand, by \eqref{eq:block-exchange-map} and
\(\Upsilon_n(D_n)=D_n\), we have
\(\Upsilon_n(\{(x,y)\in D_n:x_{[-r,0)}=u,\ y_{[0,r)}=v\})
=\{(x,y)\in D_n:x_{[-r,0)}=u,\ x_{[0,r)}=v\}.\)
The argument establishing \eqref{eq:swap-invariance}, applied to
\(\mu\times\mu\), shows that \(\Upsilon_n\) preserves
\((\mu\times\mu)|_{D_n}\). Hence
\[
 \begin{aligned}
 &(\mu\times\mu)(\{(x,y)\in D_n:
 x_{[-r,0)}=u,\ y_{[0,r)}=v\})\\
 &\quad=(\mu\times\mu)(\{(x,y)\in D_n:
                 x_{[-r,0)}=u,\ x_{[0,r)}=v\})\\
 &\quad\leq(\mu\times\mu)(\{x\in\Sigma:x_{[-r,r)}=uv\}\times\Sigma)\\
 &\quad=\mu(\{x\in\Sigma:x_{[-r,r)}=uv\})
 =\mu([uv]_\Sigma)=\mathbb{P}^{(1)}_{\mu,r}(u,v).
 \end{aligned}
\]
These estimates, together with \eqref{eq:total-variation}, give
\[
 \begin{aligned}
 &\TV(\mathbb{P}^{(1)}_{\mu,r},\mathbb{P}^{(2)}_{\mu,r})
 =1-\sum_{u,v\in\cL(\Sigma,r)}\min\left\{\mathbb{P}^{(1)}_{\mu,r}(u,v),\mathbb{P}^{(2)}_{\mu,r}(u,v)\right\}\\
 &\quad \leq1-\sum_{u,v\in\cL(\Sigma,r)}
 (\mu\times\mu)(\{(x,y)\in D_n:x_{[-r,0)}=u,\ y_{[0,r)}=v\})\\
 &\quad =1-(\mu\times\mu)(D_n),
 \end{aligned}
\]
implying that, for any \(r\in\N\),
\begin{equation}\label{eq:uniform-block-total-variation}
\TV(\mathbb{P}^{(1)}_{\mu,r},\mathbb{P}^{(2)}_{\mu,r})
\leq\liminf_{n\to\infty}\left(1-(\mu\times\mu)(D_n)\right)
=1-\limsup_{n\to\infty}(\mu\times\mu)(D_n)
\leq1-\frac{\kappa_{\ell^{*}}^2}{4}<1.
\end{equation}

Let \(\pi_-\), \(\pi_+\), \(\mathbb{P}\), and \(\mathbb{Q}\) be as in Lemma~\ref{lem:berbee-symbolic}.
For any \(r\in\N\), define
\(\rho_r:\pi_-(\Sigma)\times\pi_+(\Sigma)\to\cL(\Sigma,r)^2\) by
\[
 \rho_r(a,b)=(a_{[-r,0)},b_{[0,r)})
\]
for any \((a,b)\in\pi_-(\Sigma)\times\pi_+(\Sigma)\).
By the definitions of \(\mathbb{P}\),
\(\mathbb{Q}\), \(\rho_r\) and
the \(\sigma\)-invariance of \(\mu\), for any \(u,v\in\cL(\Sigma,r)\), we have
 \begin{align}
 &\mathbb{P}(\rho_r^{-1}(\{(u,v)\}))
 =\mu(\{x\in\Sigma:x_{[-r,0)}=u,\ x_{[0,r)}=v\})\notag\\
 &\quad =\mu(\sigma^r([uv]_\Sigma))
 =\mu([uv]_\Sigma)=\mathbb{P}^{(1)}_{\mu,r}(u,v),
 \label{eq:joint-coordinate-correspondence}
 \end{align}
and
\begin{align}
 &\mathbb{Q}(\rho_r^{-1}(\{(u,v)\}))
 =\mu(\{x\in\Sigma:x_{[-r,0)}=u\})\,
                  \mu(\{x\in\Sigma:x_{[0,r)}=v\})\notag\\
 &\quad =\mu(\sigma^r([u]_\Sigma))\,\mu([v]_\Sigma)
 =\mu([u]_\Sigma)\,\mu([v]_\Sigma)=\mathbb{P}^{(2)}_{\mu,r}(u,v).
 \label{eq:product-coordinate-correspondence}
 \end{align}
As in the proof of Lemma~\ref{lem:finite-replacement-invariance},
it can be verified that the sets \(\rho_r^{-1}(T)\), with \(r\in\N\) and
\(T\subseteq\cL(\Sigma,r)^2\), form an algebra generating the Borel
\(\sigma\)-algebra of \(\pi_-(\Sigma)\times\pi_+(\Sigma)\).
Since \(\frac{\mathbb{P}+\mathbb{Q}}{2}\) is a probability measure on
\(\pi_-(\Sigma)\times\pi_+(\Sigma)\) with its Borel \(\sigma\)-algebra,
\cite[Theorem~0.7]{Walters1982} implies that, for any Borel set
\(B\subseteq\pi_-(\Sigma)\times\pi_+(\Sigma)\) and any \(\varepsilon>0\),
there exist \(r\in\N\) and \(T_B\subseteq\cL(\Sigma,r)^2\) such that
\((\mathbb{P}+\mathbb{Q})(B\mathbin{\triangle}\rho_r^{-1}(T_B))<\varepsilon.\)
This, together with finite additivity,
\eqref{eq:joint-coordinate-correspondence}, \eqref{eq:product-coordinate-correspondence}, and
\eqref{eq:uniform-block-total-variation}, implies
\[
 \begin{aligned}
 &|\mathbb{P}(B)-\mathbb{Q}(B)|
\leq |\mathbb{P}(\rho_r^{-1}(T_B))-\mathbb{Q}(\rho_r^{-1}(T_B))|
 +|\mathbb{P}(B)-\mathbb{P}(\rho_r^{-1}(T_B))|
 +|\mathbb{Q}(B)-\mathbb{Q}(\rho_r^{-1}(T_B))|\\
 &\quad =|\mathbb{P}(\rho_r^{-1}(T_B))-\mathbb{Q}(\rho_r^{-1}(T_B))|
 +|\mathbb{P}(B\setminus\rho_r^{-1}(T_B))-\mathbb{P}(\rho_r^{-1}(T_B)\setminus B)|\\
 & \qquad +|\mathbb{Q}(B\setminus\rho_r^{-1}(T_B))-\mathbb{Q}(\rho_r^{-1}(T_B)\setminus B)|\\
 &\quad\leq |\mathbb{P}(\rho_r^{-1}(T_B))-\mathbb{Q}(\rho_r^{-1}(T_B))|
 +(\mathbb{P}+\mathbb{Q})(B\mathbin{\triangle}\rho_r^{-1}(T_B))\\
 &\quad <\left|\sum_{(u,v)\in T_B}
 \left(\mathbb{P}^{(1)}_{\mu,r}(u,v)-\mathbb{P}^{(2)}_{\mu,r}(u,v)\right)\right|
 +\varepsilon\\
 &\quad \leq\TV(\mathbb{P}^{(1)}_{\mu,r},\mathbb{P}^{(2)}_{\mu,r})
       +\varepsilon
       \leq 1-\frac{\kappa_{\ell^{*}}^2}{4}+\varepsilon.
 \end{aligned}
\]
Letting \(\varepsilon \searrow 0\) and taking the supremum
over all such \(B\), we get
\[
 \TV(\mathbb{P},\mathbb{Q})\leq1-\frac{\kappa_{\ell^{*}}^2}{4}<1.
\]

Since \(\mu\) is an ergodic invariant probability measure and
\(\TV(\mathbb{P},\mathbb{Q})<1\), it follows from
Lemma~\ref{lem:berbee-symbolic} that there exist
\(p\in\N\) and pairwise disjoint Borel sets
\(\Sigma_0,\ldots,\Sigma_{p-1}\subseteq\Sigma\) satisfying properties
(i)--(iii), with \(\mu_i=\text{\mbox{\(p\,\mu\)}}|_{\Sigma_i}\).

For any \(i\in\{0,\ldots,p-1\}\), let
\[
 \widetilde{\Sigma}_i
 =\bigcap_{n\in\Z}\sigma^{-n}(\Sigma_{(i+n)\bmod p}).
\]
These sets are Borel and pairwise disjoint, with
\(\widetilde{\Sigma}_i\subseteq\Sigma_i\).
By Lemma~\ref{lem:berbee-symbolic}~(i) and the
\(\sigma\)-invariance of \(\mu\), for any
\(i\in\{0,\ldots,p-1\}\) and any \(n\in\Z\),
\[
 \mu(\Sigma_i\mathbin{\triangle}
       \sigma^{-n}(\Sigma_{(i+n)\bmod p}))=0,
\]
implying that
\[
 \mu(\Sigma_i\setminus\widetilde{\Sigma}_i)
 \leq\sum_{n\in\Z}\mu(\Sigma_i\mathbin{\triangle}
       \sigma^{-n}(\Sigma_{(i+n)\bmod p}))=0,
\]
and consequently
\[
 \mu\left(\bigcup_{i=0}^{p-1}\widetilde{\Sigma}_i\right)=1,
 \quad \mu(\widetilde{\Sigma}_i)=\frac{1}{p},\quad\text{and}\quad
  \mu_i=\text{\mbox{\(p\,\mu\)}}|_{\widetilde{\Sigma}_i}.
\]
Since \(\sigma\) is invertible, we have
\[
 \sigma(\widetilde{\Sigma}_i)
 =\bigcap_{n\in\Z}\sigma^{1-n}(\Sigma_{(i+n)\bmod p})
 =\bigcap_{n\in\Z}\sigma^{-n}(\Sigma_{(i+n+1)\bmod p})
 =\widetilde{\Sigma}_{(i+1)\bmod p}.
\]

Since \(\widetilde{\Sigma}_i\) is a Borel subset of the compact metric
space \(\Sigma\) and \(\mu_i(\widetilde{\Sigma}_i)=1\),
\((\widetilde{\Sigma}_i,\mu_i|_{\widetilde{\Sigma}_i})\) is a
standard Borel probability space.
By Lemma~\ref{lem:berbee-symbolic}~(ii) and
\(\sigma^p(\widetilde{\Sigma}_i)=\widetilde{\Sigma}_i\), the homeomorphism \(\sigma^p|_{\widetilde{\Sigma}_i}\) preserves
\(\mu_i|_{\widetilde{\Sigma}_i}\).

Fix any \(i\in\{0,\ldots,p-1\}\). Let
\[
 \alpha_i=\{\widetilde{\Sigma}_i\cap[u]_\Sigma:
              u\in\cL(\Sigma,p)\}
              \quad\text{and}\quad
 \mathscr{P}_m=\bigvee_{j=-m}^{m}\sigma^{-jp}(\alpha_i)
 \quad(m\in\N_0).
\]
These are finite Borel partitions of \(\widetilde{\Sigma}_i\).
We verify conditions (i)--(iii) of Lemma~\ref{lem:friedman-ornstein}
for the sequence \((\mathscr{P}_m)_{m\in\N_0}\):
\begin{enumerate}
\item[(i)] For any \(m\in\N_0\),
\[
 \mathscr{P}_{m+1}
 =\sigma^{(m+1)p}(\alpha_i)\vee\mathscr{P}_m
     \vee\sigma^{-(m+1)p}(\alpha_i),
\]
so \(\mathscr{P}_{m+1}\) is a refinement of \(\mathscr{P}_m\).
\item[(ii)] For any \(n\in\N_0\) and any \(x\in\widetilde{\Sigma}_i\),
\[
\{y\in\widetilde{\Sigma}_i:
       y_{[-np,(n+1)p)}=x_{[-np,(n+1)p)}\}
 =\bigcap_{j=-n}^{n}
 \sigma^{-jp}\bigl(\widetilde{\Sigma}_i
       \cap[x_{[jp,(j+1)p)}]_\Sigma\bigr)
 \in\mathscr{P}_n.
\]
These cylinders generate the Borel \(\sigma\)-algebra of
\(\widetilde{\Sigma}_i\). Thus,
\(\mathscr{A}\!\left(\bigcup_{m\in\N_0}\mathscr{P}_m\right)
 =\{B\subseteq\widetilde{\Sigma}_i:B\text{ is Borel}\}.\)

\item[(iii)] Fix any \(m\in\N_0\). For any \mbox{\(n\in\N\)} and any integer \mbox{\(k>1\)},
\[
\bigvee_{j=-n}^{0}\sigma^{-jp}(\mathscr{P}_m)
 =\bigvee_{q=-n-m}^{m}\sigma^{-qp}(\alpha_i)
 \quad \text{and} \quad
 \bigvee_{j=k+2m}^{k+2m+n}\sigma^{-jp}(\mathscr{P}_m)
 =\bigvee_{q=k+m}^{k+n+3m}\sigma^{-qp}(\alpha_i).
\]
For any atoms
\(A\in\bigvee_{j=-n}^{0}\sigma^{-jp}(\mathscr{P}_m)\) and
\(B\in\bigvee_{j=k+2m}^{k+2m+n}\sigma^{-jp}(\mathscr{P}_m)\),
these imply that there exist words
\(u_{-n-m},\ldots,u_m\in\cL(\Sigma,p)\) and
\(v_{k+m},\ldots,v_{k+n+3m}\in\cL(\Sigma,p)\) such that
\(A=\bigcap_{q=-n-m}^{m}
       \sigma^{-qp}(\widetilde{\Sigma}_i\cap[u_q]_\Sigma)\) and
\(B=\bigcap_{q=k+m}^{k+n+3m}
       \sigma^{-qp}(\widetilde{\Sigma}_i\cap[v_q]_\Sigma)\).
Since \(A\) and \(B\) are nonempty, the concatenations
\(u=u_{-n-m}\cdots u_m\) and
\(v=v_{k+m}\cdots v_{k+n+3m}\)
belong to \(\cL(\Sigma,(n+2m+1)p)\) and satisfy
\[
 A=\widetilde{\Sigma}_i\cap\sigma^{(n+m)p}([u]_\Sigma)
 \quad\text{and}\quad
 B=\widetilde{\Sigma}_i\cap\sigma^{-(k+m)p}([v]_\Sigma).
\]
By the \(\sigma^{(m+1)p}\)-invariance of \(\mu_i\) and
\(\mu_i(\widetilde{\Sigma}_i)=1\),
\[
 \begin{aligned}
 &\mu_i(A\cap B)-\mu_i(A)\mu_i(B)\\
 &=\mu_i(\sigma^{(m+1)p}(A)\cap\sigma^{(m+1)p}(B))
   -\mu_i(\sigma^{(m+1)p}(A))\,\mu_i(\sigma^{(m+1)p}(B))\\
 &=\mu_i(\sigma^{(n+2m+1)p}([u]_\Sigma)
                  \cap\sigma^{-(k-1)p}([v]_\Sigma))
   -\mu_i(\sigma^{(n+2m+1)p}([u]_\Sigma))\,
    \mu_i(\sigma^{-(k-1)p}([v]_\Sigma)).
 \end{aligned}
\]
Thus, Lemma~\ref{lem:berbee-symbolic}~(iii) gives
\[
 \begin{aligned}
 &\sup_{n\in\N}
  \sum_{A\in\bigvee_{j=-n}^{0}\sigma^{-jp}(\mathscr{P}_m)}
  \sum_{B\in\bigvee_{j=k+2m}^{k+2m+n}\sigma^{-jp}(\mathscr{P}_m)}
  |\mu_i(A\cap B)-\mu_i(A)\mu_i(B)|\\
 &\quad\leq\sup_{r,s\in\N}
  \sum_{\substack{u\in\cL(\Sigma,r)\\v\in\cL(\Sigma,s)}}
  \left|\mu_i(\sigma^r([u]_\Sigma)
                   \cap\sigma^{-(k-1)p}([v]_\Sigma))
   -\mu_i(\sigma^r([u]_\Sigma))\,
        \mu_i(\sigma^{-(k-1)p}([v]_\Sigma))\right|\\
  &\quad\rightarrow0
  \quad(k\to\infty),
 \end{aligned}
\]
implying that, for any \(\varepsilon>0\), there exists an integer
\(k_0>1\) such that
\[
 \sum_{A\in\bigvee_{j=-n}^{0}\sigma^{-jp}(\mathscr{P}_m)}
 \sum_{B\in\bigvee_{j=k_0+2m}^{k_0+2m+n}\sigma^{-jp}(\mathscr{P}_m)}
 |\mu_i(A\cap B)-\mu_i(A)\mu_i(B)|<\varepsilon
 \quad(n\in\N).
\]
Hence, each \(\mathscr{P}_m\) is weak Bernoulli for
\(\sigma^p|_{\widetilde{\Sigma}_i}\).
\end{enumerate}

By Lemma~\ref{lem:friedman-ornstein},
\((\widetilde{\Sigma}_i,\text{\mbox{\(p\,\mu\)}}|_{\widetilde{\Sigma}_i},
\sigma^p|_{\widetilde{\Sigma}_i})\) is
measure-theoretically isomorphic to a Bernoulli scheme
for any \(i\in\{0,\ldots,p-1\}\).
Thus, \(\mu\) is Bernoulli up to a finite
period.
\end{proof}

\begin{proof}[Proof of Theorem~\ref{thm:finiteness}]
Since \(\htop(\Sigma)\geq0\), Theorems~\ref{thm:zero-entropy}
and~\ref{thm:positive-entropy} show that there are only finitely many ergodic
measures of maximal entropy, each of which is Bernoulli up to a finite period.

Denote all the ergodic measures of maximal entropy by
\(\mu_1,\ldots,\mu_N\), where \(N\geq1\), and let
\(\mu=\frac{1}{N}\sum_{n=1}^{N}\mu_n\).
Any automorphism of \(\Sigma\) permutes these measures, since it
preserves ergodicity and measure-theoretic entropy. Hence, \(\mu\)
is characteristic. By the affinity of entropy
(\mbox{\cite[Theorem~8.1]{Walters1982}}),
\[
h_\mu(\sigma|_\Sigma)
=\frac{1}{N}\sum_{n=1}^{N}h_{\mu_n}(\sigma|_\Sigma)
=\htop(\Sigma).
\qedhere
\]
\end{proof}

\begin{corollary}
\label{cor:minimal-forbidden-entropy}
If \(\Sigma\) is a subshift over a finite alphabet and
\(\hM(\Sigma)<\htop(\Sigma)\),
then \(\Sigma\) has only finitely many ergodic measures of maximal entropy,
each of which is Bernoulli up to a finite period.
\end{corollary}

\begin{proof}
It follows directly from~\mbox{\eqref{eq:entropy-comparison}}
and Theorem~\mbox{\ref{thm:finiteness}}.
\end{proof}

\begin{remark}
\label{rem:buzzi-comparison}
\begin{enumerate}
\renewcommand{\labelenumi}{(\roman{enumi})}
\item By~\mbox{\eqref{eq:entropy-comparison}} and
Example~\mbox{\ref{ex:strict-entropy-classes}}~(1),
Theorem~\mbox{\ref{thm:finiteness}} strictly extends Buzzi's result
\mbox{\cite[Theorem~1~(1)]{Buzzi2005}} beyond subshifts of quasi-finite type.
\item Corollary~\mbox{\ref{cor:minimal-forbidden-entropy}} gives a negative
answer to Question~\ref{ques:buzzi-maximal-measures}.
Our proof of Theorem~\ref{thm:finiteness} also provides a
direct proof of the conclusion on ergodic measures of maximal entropy in \cite[Theorem~1~(1)]{Buzzi2005}
without using Markov diagrams, thereby partially answering Buzzi's question
in \cite[Section~1.4, p.~380]{Buzzi2005}.
\end{enumerate}
\end{remark}

\section{Applications}
\label{sec:applications}

We apply Theorem~\ref{thm:finiteness} to power-free subshifts,
which arise from the study of repetition avoidance in
combinatorics on words
(see \cite[Section~2.2]{ClimenhagaPowerFree2025}).
Let \(d\in\N\) with \(d\geq2\), let \(\beta>1\), and put
\(\mathcal{A}=\{a_0,\ldots,a_{d-1}\}\).
Given \(n\in\N\), a word \(u\in\mathcal{A}^n\), and a rational number
\(\alpha>1\) with \(\alpha n\in\N\), the
\emph{\(\alpha\)-power} of \(u\) is the word \(u^\alpha\) formed by
the first \(\alpha n\) symbols of \(u^\infty=uuuu\cdots\).
Let
\[
 \mathcal{F}_\beta^{(d)}
 =\{u^\alpha:n\in\N,\ u\in\mathcal{A}^n,\
       \mathbb{Q}\ni\alpha\geq\beta,\ \alpha n\in\N\},
\]
and
\[
 \mathcal{F}_{\beta^+}^{(d)}
 =\{u^\alpha:n\in\N,\ u\in\mathcal{A}^n,\
       \mathbb{Q}\ni\alpha>\beta,\ \alpha n\in\N\}.
\]
Following \cite[Section~1.2]{ClimenhagaPowerFree2025},
we define the \(d\)-ary \emph{\(\beta\)-free} and
\emph{\(\beta^+\)-free} subshifts, respectively, by
\[
 \Sigma_\beta^{(d)}=\Sigma_{\mathcal{F}_\beta^{(d)}}
 \quad\text{and}\quad
 \Sigma_{\beta^+}^{(d)}=\Sigma_{\mathcal{F}_{\beta^+}^{(d)}},
\]
whenever these sets are nonempty.
A finite word over \(\mathcal{A}\) is \(\beta\)-free or
\(\beta^+\)-free if it contains no factor in
\(\mathcal{F}_\beta^{(d)}\) or \(\mathcal{F}_{\beta^+}^{(d)}\),
respectively.

\begin{proposition}
\label{prop:power-free-entropy}
Let \(d\in\N\) with \(d\geq2\) and \(\beta>1\).
If \(\Sigma\) is either \(\Sigma_\beta^{(d)}\) or
\(\Sigma_{\beta^+}^{(d)}\) and is nonempty, then
\[
 \hG(\Sigma)\leq\frac{1}{\beta}\htop(\Sigma).
\]
\end{proposition}

\begin{proof}
We first consider \(\Sigma=\Sigma_\beta^{(d)}\).
Let \(\mathscr{D}\) be the set of all finite \(\beta\)-free words.
For any \(n\in\N_0\), put
\(\mathscr{D}_n=\mathscr{D}\cap\mathcal{A}^n\).
By \cite[p.~303, Eq.~(1) and the following paragraph]{LouidorMarcusPavlov2013},
\begin{equation}\label{eq:power-free-word-growth}
 \lim_{n\to\infty}\frac{1}{n}\log\#\mathscr{D}_n
 =\htop(\Sigma).
\end{equation}

Let
\[
 \mathcal{F}
 =\left\{w\in\bigcup_{n\in\N}\mathcal{A}^n\setminus\mathscr{D}:
   \text{every proper factor of }w\text{ belongs to }\mathscr{D}\right\}.
\]
Since all finite factors of each point of \(\Sigma\) belong to
\(\mathscr{D}\) and \(\mathcal{F}\cap\mathscr{D}=\varnothing\),
we have \(\Sigma\subseteq\Sigma_{\mathcal{F}}\).
Conversely, for any \(x\in\Sigma_{\mathcal{F}}\),
if some finite factor of \(x\) does not belong to \(\mathscr{D}\),
choose a factor of minimum length with this property.
All its proper factors belong to \(\mathscr{D}\), so it belongs to
\(\mathcal{F}\). Thus \(x\notin\Sigma_{\mathcal{F}}\), a contradiction.
Hence \(x\in\Sigma\), proving
\(\Sigma_{\mathcal{F}}=\Sigma\).

Fix any \(n\in\N\).
For any \(w\in\mathcal{F}_n\), we have \(w\notin\mathscr{D}\),
so \(w\) contains a factor \(u^\alpha\) for some nonempty finite
word \(u\) over \(\mathcal{A}\) and some \(\alpha\in\mathbb{Q}\)
with \(\alpha\geq\beta\).
This factor cannot be proper, by the definition of \(\mathcal{F}\),
so \(w=u^\alpha\).
Then \(n=\alpha|u|\), and hence
\(1\leq |u|\leq\lfloor\frac{n}{\beta}\rfloor<n\).
The word \(u\) is a proper prefix of \(w\), so
\(u\in\mathscr{D}_{|u|}\).
The word \(w\) consists of the first \(n\) symbols of \(u^\infty\),
so it is uniquely determined by \(u\).
Consequently,
\[
 \#\mathcal{F}_n
 \leq\sum_{p=1}^{\lfloor\frac{n}{\beta}\rfloor}\#\mathscr{D}_p.
\]
For any \(\varepsilon>0\),
\eqref{eq:power-free-word-growth} gives a constant
\(C_\varepsilon\geq1\) such that
\(\#\mathscr{D}_p\leq C_\varepsilon e^{(\htop(\Sigma)+\varepsilon)p}\)
for all \(p\in\N\).
Since \(\htop(\Sigma)\geq0\), we obtain
\[
 \#\mathcal{F}_n
 \leq C_\varepsilon n e^{\frac{(\htop(\Sigma)+\varepsilon)n}{\beta}}
 \quad(n\in\N).
\]
This, together with \(\Sigma_{\mathcal{F}}=\Sigma\) and
\eqref{eq:generating-forbidden-entropy}, implies that
\[
 \hG(\Sigma)
 \leq\limsup_{n\to\infty}\frac{1}{n}\log\#\mathcal{F}_n
 \leq\frac{\htop(\Sigma)+\varepsilon}{\beta}.
\]
Letting \(\varepsilon\searrow0\) proves the assertion in this case.

For \(\Sigma=\Sigma_{\beta^+}^{(d)}\), the same proof applies with
\(\beta^+\)-free words in place of \(\beta\)-free words and
\(\alpha>\beta\) in place of \(\alpha\geq\beta\).
\end{proof}

\begin{corollary}
\label{cor:power-free-measures}
Let \(d\in\N\) with \(d\geq2\) and \(\beta>1\).
If \(\Sigma\) is either \(\Sigma_\beta^{(d)}\) or
\(\Sigma_{\beta^+}^{(d)}\) and \(\htop(\Sigma)>0\), then
\(\Sigma\) has only finitely many ergodic measures of maximal
entropy, each of which is Bernoulli up to a finite period.
Moreover, \(\Sigma\) admits a characteristic measure of maximal
entropy and is not of quasi-finite type.
\end{corollary}

\begin{proof}
Proposition~\ref{prop:power-free-entropy} gives
\[
 \hG(\Sigma)\leq\frac{1}{\beta}\htop(\Sigma)
 <\htop(\Sigma).
\]
The conclusions about measures follow from
Theorem~\ref{thm:finiteness}.
If \(x\in\Sigma\) satisfies \(\sigma^p(x)=x\) for some
\(p\in\N\), choose an integer \(k>\beta\).
Then \(x_{[0,kp)}=(x_{[0,p)})^k\), contradicting either
power-free condition.
Thus \(\Sigma\) has no periodic points and cannot be of
quasi-finite type, by \cite[Theorem~1~(3)]{Buzzi2005}.
\end{proof}

\begin{remark}
\label{rem:power-free-bernoulli}
Fix any \(\beta>3\).
By \cite[Theorem~2.7 and its proof]{ClimenhagaPowerFree2025},
there exists an integer \(d_0(\beta)\geq2\) such that,
for any integer \(d\geq d_0(\beta)\) and any
\(\Sigma\in\{\Sigma_\beta^{(d)},\Sigma_{\beta^+}^{(d)}\}\),
we have \(\htop(\Sigma)>0\), and \(\Sigma\) has a
unique measure of maximal entropy \(\mu\) with the \(K\)-property.
By Corollary~\ref{cor:power-free-measures}, \(\mu\) is Bernoulli up to
some finite period \(p\), with cyclic components
\(\Sigma_0,\ldots,\Sigma_{p-1}\).
Then
\(\mu(\Sigma_0)=\frac{1}{p}\) and \(\sigma^p(\Sigma_0)=\Sigma_0\).
Since the \(K\)-property implies mixing
(\cite[Section~2.3]{ClimenhagaPowerFree2025}),
we have \(\frac{1}{p}=\mu(\Sigma_0\cap\sigma^{-np}(\Sigma_0))
\to\mu(\Sigma_0)^2=\frac{1}{p^2}\) as \(n\to\infty\), so \(p=1\).
Thus \((\Sigma,\mu,\sigma|_\Sigma)\) is
measure-theoretically isomorphic to a Bernoulli scheme.
This gives a partial affirmative answer to
\cite[Question~1.4]{ClimenhagaPowerFree2025} for \(d\geq d_0(\beta)\).
\end{remark}

\begin{remark}
\label{rem:classical-power-free}
Consider subshifts \(\Sigma_2^{(3)}\) and \(\Sigma_3^{(2)}\).
By \cite[Theorems~6.10 and~6.5]{Kobayashi1986} and
\eqref{eq:power-free-word-growth}
(see also \cite[Section~2.2.2]{ClimenhagaPowerFree2025}),
both subshifts have positive topological entropy.
By Proposition~\ref{prop:power-free-entropy} and
Corollary~\ref{cor:power-free-measures}, we have
\[
 \Sigma_2^{(3)}
 \in\mathscr{C}_{\mathcal{G}}^{(\{a_0,a_1,a_2\})}
 \setminus\mathscr{C}_{\mathcal{SC}}^{(\{a_0,a_1,a_2\})}
 \quad\text{and}\quad
 \Sigma_3^{(2)}
 \in\mathscr{C}_{\mathcal{G}}^{(\{a_0,a_1\})}
 \setminus\mathscr{C}_{\mathcal{SC}}^{(\{a_0,a_1\})}.
\]
Thus, for \(\#\mathcal{A}\in\{2,3\}\),
\(\mathscr{C}_{\mathcal{G}}^{(\mathcal{A})}\) contains classical
subshifts that are not of quasi-finite type.
\end{remark}

\end{document}